\documentclass[11pt,a4paper]{article}

\usepackage[T1]{fontenc}
\usepackage[utf8]{inputenc}
\usepackage[font=small,labelfont={bf}]{caption}
\usepackage{epigraph}
\usepackage[explicit]{titlesec}
\usepackage{fancyhdr}
\usepackage{authblk}
\usepackage{xspace}	
\usepackage{amsmath}
\usepackage{amssymb}
\usepackage{mathtools}
\usepackage{amsthm}
\usepackage{enumerate}
\usepackage[shortlabels]{enumitem}
\usepackage{pinlabel}
\usepackage{subfig}
\usepackage{booktabs}
\usepackage{bbm}
\providecommand{\keywords}[1]{\par\smallskip\textbf{\textbf{Keywords. }} #1}
\newtheorem{theorem}{Theorem}[section]
\newtheorem{lemma}[theorem]{Lemma}

\theoremstyle{definition}
\newtheorem{definition}[theorem]{Definition}

\newtheorem{example*}{Example}
\newtheorem*{running_example*}{Running Example}

\numberwithin{figure}{section}
\usepackage{xcolor}

\definecolor{thm_color_red}{RGB}{245, 188, 174}
\definecolor{thm_color}{RGB}{238, 187, 240}
\definecolor{dark_green}{RGB}{0, 123, 2}
\definecolor{very_light_blue}{RGB}{211, 255, 255}
\definecolor{dark_blue}{RGB}{0, 80, 255}
\definecolor{light_orange}{RGB}{255, 165, 0}
\definecolor{red_one}{RGB}{183, 0, 0}
\definecolor{gold}{RGB}{255, 255, 240}
\definecolor{gradient_blue_red}{RGB}{91, 40, 128}
\definecolor{turquoise}{RGB}{29, 169, 255}
\definecolor{turquoise_darker}{RGB}{0, 174, 255}
\definecolor{block_bgcolor1}{RGB}{29, 169, 255}
\definecolor{red_alert}{RGB}{222, 20, 20}

\newcommand{\dd}{ \mathrm{d}}

\newcommand{\cD}{\mathcal{D}}
\newcommand{\cE}{\mathcal{E}}

\newcommand{\cH}{\mathcal{H}}
\newcommand{\cI}{\mathcal{I}}
\newcommand{\cJ}{\mathcal{J}}
\newcommand{\cK}{\mathcal{K}}

\newcommand{\cM}{\mathcal{M}}

\newcommand{\cP}{\mathcal{P}}

\newcommand{\cR}{\mathcal{R}}

\def\calM{\mathcal M}

\newcommand{\R}{\mathbb R}
\renewcommand{\P}{\mathbb P}

\def\e{\epsilon}

\usepackage[hidelinks,bookmarks=false]{hyperref}	

\usepackage{rotating}

\let\e\varepsilon
\let\wt\widetilde
\newcommand{\tmsp}{{Q_T}}
\newcommand{\tmspz}{{Q_{T}^0}}

\definecolor{blue_fig}{RGB}{0, 80, 255}
\definecolor{red_fig}{RGB}{183, 0, 0}

\newcommand{\labelsize}{\scriptsize}

\usepackage{stackrel}
\usepackage{mathtools} 
\usepackage{longtable} 

\usepackage{ifthen}
\newlength{\leftstackrelawd}
\newlength{\leftstackrelbwd}
\def\leftstackrel#1#2{\settowidth{\leftstackrelawd}%
{${{}^{#1}}$}\settowidth{\leftstackrelbwd}{$#2$}%
\addtolength{\leftstackrelawd}{-\leftstackrelbwd}%
\leavevmode\ifthenelse{\lengthtest{\leftstackrelawd>0pt}}%
{\kern-.5\leftstackrelawd}{}\mathrel{\mathop{#2}\limits^{#1}}}

\usepackage{xparse}
\ExplSyntaxOn
    \NewDocumentEnvironment{remark} { o }
     {
      \par\medbreak\refstepcounter{theorem}\noindent
        \IfNoValueTF{#1} { \textbf{Remark~\thetheorem.\ } }
                         { \textbf{Remark~\thetheorem\ } (#1)\textbf.\ }
       \ignorespaces
     }
     {\qed\par\medskip}
\ExplSyntaxOff
\newenvironment{assumption}%
  {\par\medbreak\refstepcounter{theorem}%
    \noindent\textbf{Assumption~\thetheorem. }}%
  {\qed\par\medskip}

\DeclareMathOperator{\dist}{dist}
\DeclareMathOperator{\supp}{supp}
\let\ds\displaystyle
\def\RelEnt{\mathcal H}

\DeclareMathOperator{\sign}{sign}

\usepackage{dsfont}
\newcommand{\bONE}{\mathds{1}}

\long\def\drop#1{}

 \def\sfE{{\mathsf E}}

  \def\sfR{{\mathsf R}}

\def\CE{\mathrm{CE}}
\newcommand\idhalf{\mathrm{id}_{1/2}}

\let\weakto\rightharpoonup

\let\longweakto\longrightharpoonup

\usepackage{fancyhdr}
\usepackage{datetime}

\title{Gamma-convergence of a gradient-flow structure to a non-gradient-flow structure}
\author[1,2]{Mark A. Peletier}
\author[1]{Mikola C. Schlottke}
\affil[1]{Department of Mathematics and Computer Science, Eindhoven University of Technology}
\affil[2]{Institute for Complex Molecular Systems (ICMS), Eindhoven University of Technology}

\begin{document}
\maketitle
\begin{abstract}
We study the asymptotic behaviour of a gradient system in a regime in which the driving energy becomes singular. For this system gradient-system convergence concepts are ineffective. We characterize the limiting behaviour in a different way, by proving  $\Gamma$-convergence of the so-called \emph{energy-dissipation functional}, which combines the gradient-system components of energy and dissipation in a single functional. The $\Gamma$-limit of these functionals again characterizes a variational evolution, but this limit functional is not the energy-dissipation functional of any gradient system.

The system in question describes the diffusion of a particle in a one-dimensional double-well energy landscape, in the limit of small noise. The wells have different depth, and in the small-noise limit the process converges to a Markov process on a two-state system, in which jumps only happen from the higher to the lower well. 

This transmutation of a gradient system into a variational evolution of non-gradient type is a model for how many one-directional chemical reactions emerge as limit of reversible ones. The $\Gamma$-convergence proved in this paper both identifies the `fate' of the gradient system for these reactions and the variational structure of the limiting irreversible reactions.
\end{abstract}

\medskip

\keywords{Kramers problem, irreversible limit, variational evolution, gradient system, EDP-convergence.}

\section{Introduction}

\subsection{Diffusion in an asymmetric potential landscape}
\begin{figure}[h!]
	\labellist
	\pinlabel  $x$ at 2750 570
	\pinlabel  $V(x)$ at 2550 1200
	\pinlabel  $x_a$ at 500 570
	\pinlabel  $x_0$ at 1080 570
	\pinlabel  $x_b$ at 1800 700
	\pinlabel  $x_{b-}$ at 1530 850
	\pinlabel  $x_{b+}$ at 2100 850
	\endlabellist
	\centering
	\includegraphics[height=4cm]{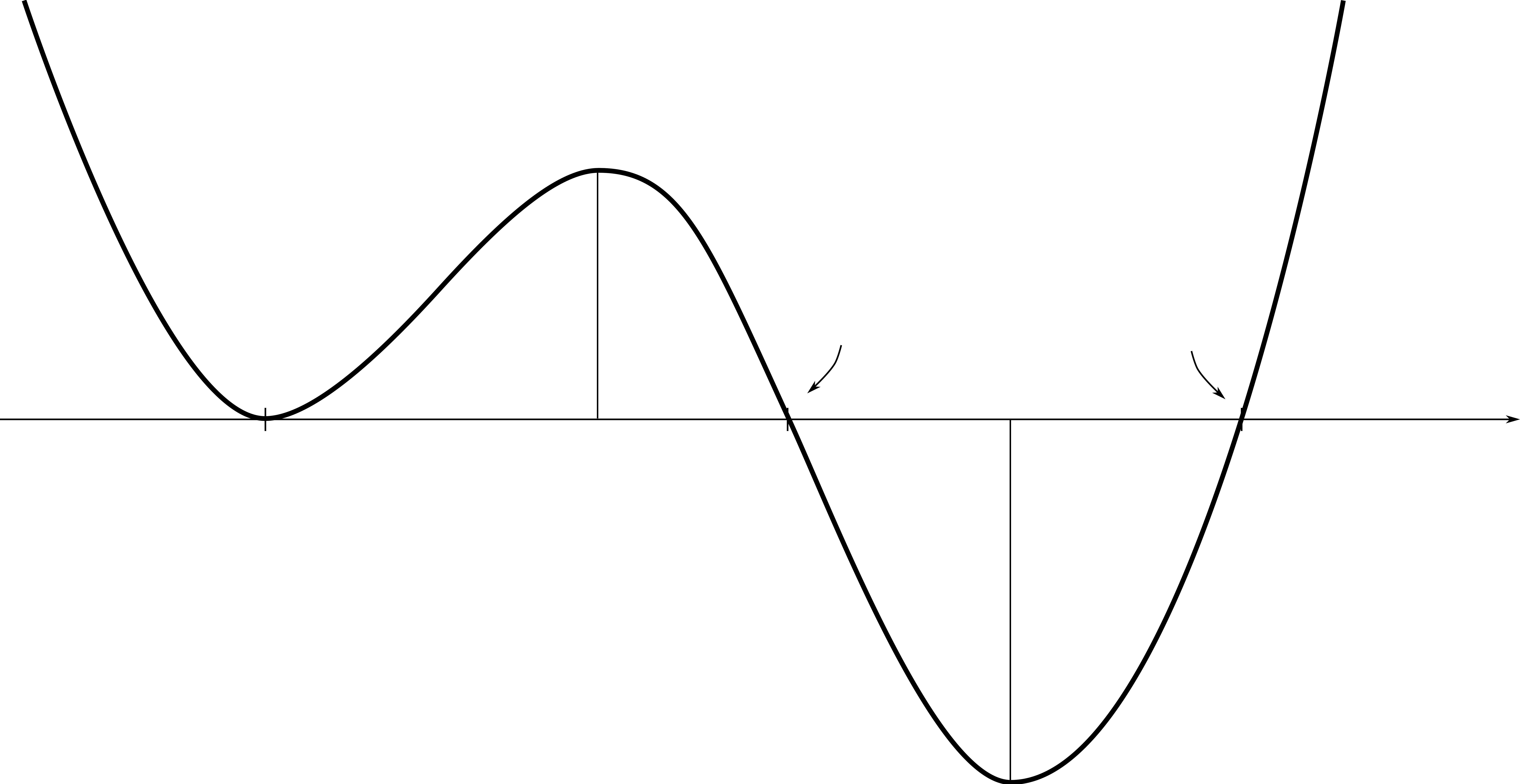}
	\caption{A typical asymmetric potential $V(x)$. }
	\label{fig:asymmetric-doublewell-potential}
\end{figure}
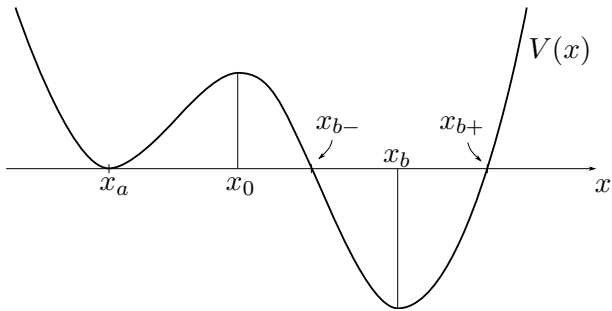

Our interest in this paper is the limit $\e\to0$ in the 
family of Fokker-Planck equations in one dimension defined by 
\begin{equation}\label{eq:intro:upscaled-FP}
\partial_t \rho_\e = \tau_\e\Bigl( \e\, \partial_{xx} \rho_\e + \partial_x \left( \rho_\e  V'\right)\Bigr), 
\quad \text{on }\R_+ \times \R.
\end{equation}
Here we take an \emph{asymmetric} double-well potential~$V:\mathbb{R}\to\mathbb{R}$ as depicted in Figure~\ref{fig:asymmetric-doublewell-potential}.
\smallskip

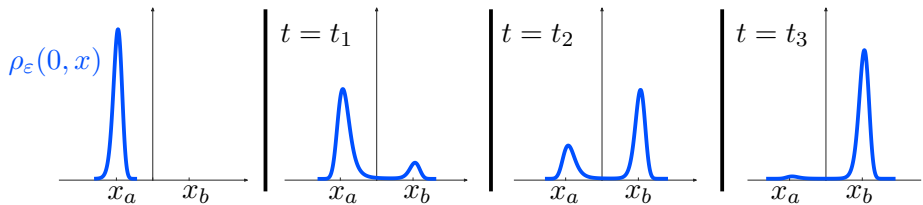
\begin{figure}[hb]
	\labellist
	\pinlabel $x_a$ at 130 0
	\pinlabel $x_b$ at 240 0
	\pinlabel $x_a$ at 460 0
	\pinlabel $x_b$ at 570 0
	\pinlabel $x_a$ at 790 0
	\pinlabel $x_b$ at 900 0
	\pinlabel $x_a$ at 1120 0
	\pinlabel $x_b$ at 1230 0
	\pinlabel $t=t_1$ at 420 240
	\pinlabel $t=t_2$ at 750 240
	\pinlabel $t=t_3$ at 1100 240
	\pinlabel {\color{dark_blue}{$\rho_\e(0,x)$}} at 30 200
	\endlabellist
	\centering
	\includegraphics[scale=.25]{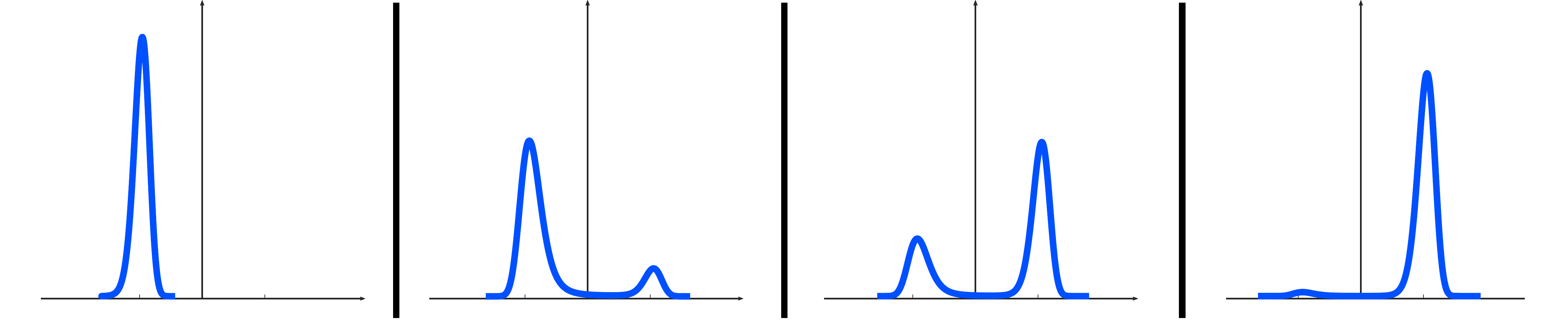}
	\caption{The time evolution of a solution~$\rho_\e(t,x)$ to~\eqref{eq:intro:upscaled-FP} whose initial distribution is supported on the left. Time increases from left to right. At the final time, the solution is close to the equilibrium distribution, which is proportional to~$\exp\{-V(x)/\e\}$. The smaller the value of~$\e$, the sharper the equilibrium distribution concentrates around the global minimum~$x_b$.}
	\label{fig:intro:evolution-FP}
\end{figure}

A typical solution~$\rho_\e(t,x)$ is displayed in Figure~\ref{fig:intro:evolution-FP}, showing mass flowing from left to right. There are two parameters, $\e>0$ and~$\tau_\e>0$.  The first parameter~$\e$ controls how fast mass can move between the potential wells, where  smaller values of~$\e$ correspond to larger transition times. The second parameter~$\tau_\e$ sets the global time scale, and is chosen  such that typical transition times from the local minimum~$x_a$ to the global minimum~$x_b$ are of order one as $\e\to0$ (see equation~\eqref{intro:eq:def-time-scale-parameter} below).

\smallskip

The small-$\e$ limit in the PDE~\eqref{eq:intro:upscaled-FP} is known as the \emph{high activation energy} limit in the context of chemical reactions. In this setting, the PDE can be derived from  the stochastic evolution of a chemical system, modelled  by a one-dimensional diffusion process $Y^\e_t = Y^\e(t)$ in $\mathbb{R}$, satisfying
\begin{equation*}
\dd Y^\e_t = - V'(Y^\e_t)\,\dd t + \sqrt{2\e}\;\dd B_t,
\end{equation*}
where~$B_t$ is a  standard Brownian motion. For example, consider a particle  starting in the left minimum~$x_a$ and propagating from left to right.  This propagation  models a reaction event in which a molecule's state changes from a low-energy state~$x_a$ via a high-energy state~$x_0$ to another low-energy state~$x_b$. The assumption of asymmetry of the potential~$V$ corresponds to modelling a reaction in which the  final energy is lower than the initial energy. The energy barrier that the particle has to overcome,~$V(x_0)-V(x_a)$, is the \emph{activation energy} of the reaction. 

Hendrik Antony Kramers was the first to translate the question of determining the rate of a chemical reaction into properties of PDEs such as~\eqref{eq:intro:upscaled-FP}~\cite{Kramers40,HanggiTalknerBorkovec90}. Decreasing~$\e$ reduces the noise level in comparison to the potential energy barrier, and a transition from~$x_a$ to~$x_b$ becomes more unlikely, and hence the average time until a transition~$x_a\rightsquigarrow x_b$  increases. Kramers derived an asymptotic expression for this average time:
\begin{multline}
\label{eq:Kramers-time}
\mathbb{E}\Bigl[\inf\{ t>0: Y_t^\e = x_b\} \Big| Y_0^\e= x_a\Bigr]\\
 =\left[1 + o(1)_{\e\to 0}\right] \frac{2\pi}{\sqrt{V''(x_a) |V''(x_0)|}} \exp\{\e^{-1}(V(x_0) - V(x_a))\},
\end{multline}
which now is  known as the \emph{Kramers formula}. It shows that the average transition time scales exponentially with respect to  the ratio of the energy barrier~$V(x_0)-V(x_a)$ to the diffusion coefficient~$\e$. For further details and background on this model, we refer to the monographs of Bovier and den Hollander~\cite{BovierDenHollander2016}, and of Berglund and Gentz~\cite{BerglundGentz2005}.
\smallskip
 
We are interested in the limit~$\e\to 0$ in the equation~\eqref{eq:intro:upscaled-FP}. In this limit we expect the solution $\rho_\e$ to concentrate at the minima $x_a$ and~$x_b$. Furthermore, transitions from left to right face a lower energy barrier than from right to left, and because of the exponential scaling in the energy barrier in~\eqref{eq:Kramers-time}, we expect that in the limit $\e\to0$ transitions occur much more often from left to right than from right to left. 

\smallskip

Since we want to follow left-to-right transitions, we choose the global  time-scale parameter~$\tau_\e$ approximately equal to the left-to-right transition time:
\begin{equation}\label{intro:eq:def-time-scale-parameter}
\tau_\e := \frac{2\pi}{\sqrt{V''(x_a) |V''(x_0)|}} \exp\{\e^{-1}(V(x_0) - V(x_a))\}.
\end{equation}
Speeding up the process~$Y^\e(t)$ by $\tau_\e$ as $X^\e(t):= Y^\e(\tau_\e t)$, the accelerated process $X^\e$ 
satisfies the SDE
\begin{equation}\label{eq:intro:upscaled-diffusion-process}
\dd X^\e_t = - \tau_\e V'(X^\e_t)\,\dd t + \sqrt{2\e \tau_\e}  \; \dd B_t,
\end{equation}
and the equation~\eqref{eq:intro:upscaled-FP} is the Fokker-Planck equation for  the transition probabilities $\rho_\e(t,\dd x) := \mathbb{P}\left(X^\e_t \in \dd x\right)$.
\smallskip

\begin{figure}[ht!]
	\labellist
	\pinlabel $x_a$ at 130 0
	\pinlabel $x_b$ at 240 0
	\pinlabel $x_a$ at 460 0
	\pinlabel $x_b$ at 570 0
	\pinlabel $x_a$ at 790 0
	\pinlabel $x_b$ at 900 0
	\pinlabel $x_a$ at 1120 0
	\pinlabel $x_b$ at 1230 0
	\pinlabel $t=t_1$ at 420 240
	\pinlabel $t=t_2$ at 750 240
	\pinlabel $t=T$ at 1080 240
	\pinlabel {\color{red_one}{$\rho_0(0,x)$}} at 30 200
	\endlabellist
	\centering
	\includegraphics[scale=.25]{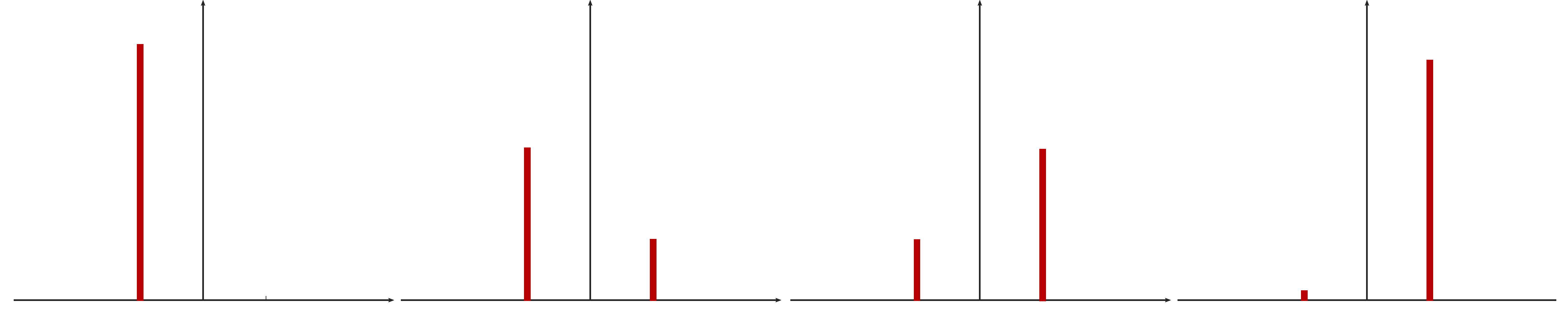}
	\caption{The time evolution of~$\rho_0$, defined as the~$\e\to 0$ limit of the solution~$\rho_\e(t,x)$ to~\eqref{eq:intro:upscaled-FP}. The initial distribution is supported solely on the left. Mass flows only from left to right, with rate one.}
	\label{fig:intro:limit-evolution-FP}
\end{figure}
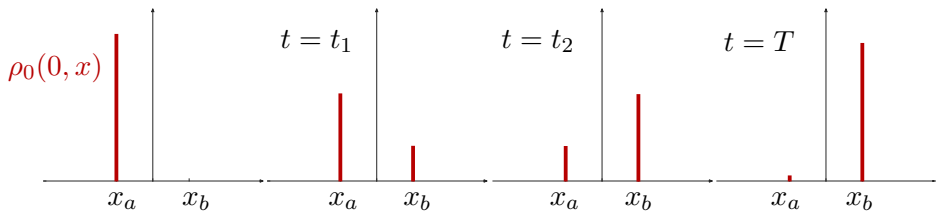

In the rescaled equation~\eqref{eq:intro:upscaled-FP} we therefore expect the limiting dynamics to be characterized by mass being transferred at rate one from the local minimum~$x_a$ to the global minimum~$x_b$, and to see no mass move in the opposite direction. In terms of the solution $\rho_\e$, we expect that 
\begin{equation}
\label{eq:expected-limit}	
\rho_\e \to \rho_0 = z \delta_{x_a} + (1-z) \delta_{x_b},
\end{equation}
where the density~$z=z(t)$ of particles at $x_a$ satisfies~$\partial_t z = - z$, corresponding to left-to-right transitions happening at rate $1$. The time evolution of the limiting density is depicted in Figure~\ref{fig:intro:limit-evolution-FP}.

The main results of this paper  imply the convergence~\eqref{eq:expected-limit}, but they provide more information: they describe the fate of the gradient-system, vari\-a\-tion\-al-evolutionary structure satisfied by~\eqref{eq:intro:upscaled-FP}. We describe this next.

\subsection{Gradient systems and convergence}

Both the convergence of stochastic processes and the convergence of PDEs are classical problems, and the particular case of the small-noise or high-activation-energy limit is very well studied; see the monographs of Berglund--Gentz  and Bovier--Den Hollander that we already mentioned for much more on this topic~\cite{BerglundGentz2005,BovierDenHollander2016}. 

In this paper, however, our main interest  in the $\e\to0$ limit of equation~\eqref{eq:intro:upscaled-FP} is the relation with \emph{convergence of gradient systems}. 
One of the main points of this paper is that while the $\e>0$ systems are of gradient type, there is no reasonable convergence that remains within the class of gradient systems. Instead we prove a convergence result to a more general variational evolution that is not of gradient type.

\medskip

In this paper we focus on gradient systems in the space of probability measures on $\R$ with a continuity-equation structure. Equation~\eqref{eq:intro:upscaled-FP} is of this form; it can be written as  the triplet of equations
\begin{subequations}
\label{eq:upscaled-pair}
\begin{alignat}2
\label{eq:upscaled-CE}
&\partial_t\rho_\e + \partial_x \,j_\e = 0, &\qquad\text{(continuity equation),}\\
&j_\e =  J_\e^\rho,&\qquad\text{(specification of flux),}\label{eq:upscaled-flux}
\\
&J_\e^\rho := -\tau_\e\left[\e\, \partial_x \rho_\e + \rho_\e\partial_x V\right],&\qquad\text{(definition of $J_\e^\rho$ in terms of $\rho_\e$)}.
\end{alignat}
\end{subequations}

For pairs $(\rho_\e,j_\e)$ satisfying~\eqref{eq:upscaled-CE}, the second equation~\eqref{eq:upscaled-flux} can formally be written as
\begin{equation}
\label{defeq:I_e}
\cI_\e(\rho_\e,j_\e) \leq 0, 
\end{equation}
in terms of the trivially nonnegative functional $\cI_\e$, 
\begin{equation}
\label{defeq:I_e2}
\cI_\e(\rho_\e,j_\e) := \frac{1}2 \int_0^T\int_\mathbb{R} \frac{1}{\e\, \tau_\e} \frac{1}{\rho(t,x)} \big|j_\e(t,x) -J_\e^\rho(t,x)\big|^2\,\dd x\dd t.
\end{equation}
By expanding the square in $\cI_\e$ (see Lemma~\ref{l:reg-rho-j} for details) one finds the equivalent  form of~\eqref{defeq:I_e},
\begin{subequations}
\label{eq:gf:Wass}	
\begin{equation}
\label{eq:EDP-formulation}
\bigl( \;\cI_\e(\rho_\e,j_\e) = \ \bigr)\quad 
E_\e(\rho)\Big|_{t=0}^{t=T} 
+ \underbrace{\int_0^T \bigl[ R_\e(\rho,j) + R_\e^*\bigl(\rho, -\mathrm D E_\e\bigl(\rho)\bigr)\Bigr]\, dt}
      _{\textstyle =: \ \cD_\e^T(\rho,j)}
\leq 0.
\end{equation}
In~\eqref{eq:EDP-formulation} the functional $E_\e$ is given as 
\begin{equation}
\label{eq:gf:Wass-E}	
E_\e(\rho) := \cH(\rho|\gamma_\e)
\quad\text{where}\quad
\gamma_\e(\dd x) := \frac1{Z_\e}e^{-V(x)/\e}\, \dd x,
\end{equation}
and $\cH(\mu|\nu)$ is the relative entropy of $\mu$ with respect to $\nu$. The dual pair $(R_\e,R_\e^*)$ of \emph{dissipation potentials}
is formally defined as 
\begin{equation}
\label{eq:gf:Wass-R}
R_\e(\rho,j) := \frac1{2\e\tau_\e} \int_\R \frac{j^2} \rho
\quad\text{and}\quad
R_\e^*(\rho,\xi) := \frac{\e\tau_\e}2 \int_\R \xi^2 \,\rho.
\end{equation}
\end{subequations}


The inequality~\eqref{eq:EDP-formulation} is known as the \emph{EDP-formulation} of the gradient system defined by $E_\e$, $R_\e$, and the continuity equation; see e.g.~\cite{AmbrosioGigliSavare2008,Peletier2014, Mielke16a} for a general discussion of gradient systems, and~\cite{PeletierRossiSavareTse20TR} for a specific treatment of gradient systems with continuity-equation structure.
The dissipation potential $R_\e^*$ in~\eqref{eq:gf:Wass-R} and its dual $R_\e$ can be interpreted as infinitesimal versions of the \emph{Wasserstein metric}, and for this reason system~\eqref{eq:upscaled-pair} or equivalently equation~\eqref{eq:intro:upscaled-FP} is known as a Wasserstein gradient flow~\cite{AmbrosioGigliSavare2008,Peletier2014,Santambrogio15}.


\medskip

The EDP-formulation~\eqref{eq:gf:Wass} can be used not only to define gradient-system solutions, but also to define \emph{convergence} of a sequence of gradient systems  to a limiting gradient system. Although this method will not be directly of use to us for the proofs in this paper, since the limiting system of this paper will not be of gradient-system type, we will use a number of elements of this method. In addition, it is useful to contrast the method of this paper with this convergence concept.

\begin{definition}[EDP-convergence]
A sequence $(E_\e,R_\e)$ \emph{EDP-converges} to a limiting gradient system $(E_0,R_0)$ if 
\begin{enumerate}
\item $E_\e \stackrel \Gamma \longrightarrow E_0$,
\item $\cD_\e^T \stackrel \Gamma \longrightarrow \cD_0^T$ for all $T$, and 
\item the limit functional $\cD_0^T$ can again be written in terms of the limiting functional $E_0$ and a dissipation potential $R_0$ as 
\begin{equation}
\label{eq:D0-dual-form}
\cD_0^T(\rho,j) = \int_0^T \bigl[R_0(\rho,j) + R_0^*(\rho,-\mathrm DE_0(\rho))\bigr]\, dt.
\end{equation}
\end{enumerate}
\end{definition}

EDP-convergence implies convergence of solutions: If $(\rho_\e,j_\e)$ is a sequence of solutions of~\eqref{eq:intro:upscaled-FP} or equivalently of~\eqref{eq:gf:Wass} that converges to a limit $(\rho_0,j_0)$, and if the initial state $\rho_\e(0)$ satisfies the well-preparedness condition 
\begin{equation}
\label{eq:well-preparedness}
E_\e(\rho_\e(t=0) )\to E_0(\rho_0(t=0)),
\end{equation}
then the limit $(\rho_0,j_0)$ is a solution of the gradient flow associated with $(E_0,R_0)$. 
See~\cite{Mielke16a,MielkeMontefuscoPeletier2020} for an in-depth discussion of EDP-conver\-gence.

\subsection{(Non-)convergence as $\e\to0$ in the Kramers problem}
\label{ss:conv-Kramers}

For symmetric potentials $V$, EDP-convergence of the gradient systems $(E_\e,R_\e)$ of~(\ref{eq:gf:Wass-E}--\ref{eq:gf:Wass-R}) has been proved in~\cite{ArnrichMielkePeletierSavareVeneroni2012,LieroMielkePeletierRenger17}.
For non-symmetric potentials as in this paper, however, we claim that the sequence $(E_\e,R_\e)$  \emph{can not} converge in this sense, and we now explain this.

\medskip
\emph{1. The functional $E_\e$ blows up.} The first argument for non-convergence follows from the singular behaviour of the driving functional $E_\e$. We can rewrite this functional as
\begin{equation}
\label{eq:rewrite-Ee-blowup}
E_\e(\rho) = \cH(\rho|Z_\e^{-1}e^{-V/\e})
= \int_\R \rho(x)\log \rho(x) \, dx + \int_\R \rho \Bigl(\frac 1\e V + \log Z_\e\Bigr).
\end{equation}
Since the normalization constant $Z_\e$ is chosen such that $\gamma_\e$ has mass one, the term in parentheses converges to $+\infty$ at all $x$ except for the global minimizer $x=x_b$ (this follows from Lemma~\ref{l:basic-props-gamma}). Therefore $E_\e$ $\Gamma$-converges to the singular limit functional 
\[
E_0(\rho) := \begin{cases}
0 & \text{if }\rho = \delta_{x_b}\\
+\infty & \text{otherwise}.
\end{cases}	
\]

This  implies that if $\rho_\e(0)$ retains any mass in the higher well around~$x_a$ as $\e\to0$, then $E_\e(\rho_\e(0))\to\infty$. The `well-preparedness condition'~\eqref{eq:well-preparedness} therefore can only be satisfied in a trivial way, with the initial mass being `already' in the lower of the two wells. Indeed, a gradient system driven by~$E_0$ admits only constants as solutions, and does not allow us to follow transitions from $x_a$ to $x_b$. 

\medskip

\emph{2. Other scalings of $E_\e$ also fail.} One could mitigate the blow-up of $E_\e$ by choosing a different scaling of $E_\e$, 
\[
\wt E_\e(\rho) := \e E_\e(\rho) = \e\int_\R \rho(x)\log \rho(x) \, dx + \int_\R \rho \bigl( V + \e\log Z_\e\bigr),
\]
which $\Gamma$-converges to the  functional $\rho \mapsto \int\rho V$. With this scaling the well-preparedness condition~\eqref{eq:well-preparedness} is simple to satisfy, and by general compactness arguments (e.g.~\cite[Ch.~10]{DalMaso93}) the correspondingly rescaled functionals $\wt\cD_\e^T := \e\cD_\e^T$ also $\Gamma$-converge to a limit $\wt\cD_0^T$. However, this limit functional~$\wt\cD_0^T$ fails to characterize an evolution; we prove this in Section~\ref{ss:renormalized} below. Other rescaling choices suffer from similar problems.

\medskip

\emph{3. EDP-convergence \textbf{should} fail.}
There also is a more abstract argument why EDP-convergence should fail, and in fact why any gradient-system convergence should fail. In the limit $\e\to0$ the ratio of forward to reverse transitions diverges, leading to a situation in which motion becomes one-directional. On the other hand, in gradient systems motion can be reversed by appropriate tilting of the driving functional. Therefore the one-directionality is incompatible with a gradient structure. 

Note that the limiting equation itself, $\dot z = -z$ (see Section~\ref{ss:limiting-functional}), \emph{can} be given a gradient structure, even many different gradient structures; one example is 
\[
E(z) := \frac12 z^2, \qquad R(\dot z) := \frac12 {\dot z}^2.
\]
Our claim here is the following: although the limiting equation can in fact be given a multitude of gradient structures, none of these structures can be found as the limit of the Wasserstein gradient structure of equation~\eqref{eq:intro:upscaled-FP}. The simplest proof of this statement is the $\Gamma$-convergence theorem that we prove in this paper (Theorem~\ref{thm:intro-main-result}), which identifies the limit functional; this functional does not generate a gradient structure. %

\bigskip
Summarizing, although for each $\e>0$ the equation~\eqref{eq:intro:upscaled-FP} is a Wasserstein gradient flow with components $E_\e$ and $R_\e$, these components diverge in  the limit $\e\to0$, and only trivial gradient-system convergence is possible. 

On the other hand, the functional $\cI_\e$ combines the components $E_\e$,  $R_\e$, and $R_\e^*$ in such a way that their divergences compensate each other; in the case of solutions of~\eqref{eq:intro:upscaled-FP}, $\cI_\e$ even is zero for all $\e$. This suggests that $\cI_\e$ is a better candidate for a variational convergence analysis, and the rest of this paper is devoted to this. Indeed we find below that the limit of $\cI_\e$ is not of gradient-flow structure, confirming the earlier suggestion that the sequence leaves the class of gradient systems. 

\begin{remark}
In~\cite{PeletierSavareVeneroni2010,PeletierSavareVeneroni12} one of us developed convergence results for this same limit $\e\to0$ for the case of a symmetric potential $V$, using a functional framework based on $L^2$-spaces that are weighted with the invariant measure~$\gamma_\e$. This approach suffers from a similar problem as the Wasserstein-based approach above. The limiting state space is the space $L^2$, weighted by the limiting invariant measure $\delta_b$, which is a one-dimensional function space; in combination with the constraint of unit mass, the effective state space is a singleton. Consequently the limiting evolution would be trivial. 
\end{remark}

\subsection{Main result---$\Gamma$-convergence of $\cI_\e$}
\label{sec:main-result}

In the previous section we introduced the functional $\mathcal I_\e$ of a  pair $(\rho,j)$ with the property that solutions of the equation~\eqref{eq:intro:upscaled-FP} are minimizers of $\mathcal I_\e$ at value zero. As for gradient structures, we can therefore reformulate the question of convergence as $\e\to0$ in terms of $\Gamma$-convergence of these functionals.  The main questions then are:
\begin{enumerate}[label=(\roman*)]
	\item \emph{Compactness:} For a family of pairs~$(\rho_\e',j_\e')$ depending on~$\e$, does boundedness of~$\mathcal{I}_\e(\rho_\e',j_\e')$ imply the existence of a subsequence of~$(\rho_\e',j_\e')$ that converges in a certain topology~$\mathcal{T}$?
\item \emph{Convergence along sequences:} Is there a limit functional~$\mathcal{I}_0$ such that 
	\begin{equation*}
	    \Gamma-\lim_{\varepsilon\to 0}\mathcal{I}_\varepsilon = \mathcal{I}_0\, ?
	\end{equation*}
	\item \emph{Limit equation:} Does the equation $\mathcal I_0(\rho,j)=0$  characterize the evolution of $(\rho,j)$? 
\end{enumerate}
We answer the first question in Theorem~\ref{th:compactness-lower-bound-untransformed}, which establishes that sequences $(\rho_\e',j_\e')$ such that~$\mathcal{I}_\e(\rho_\e',j_\e')$ remains bounded are compact in a certain topology. 
\smallskip

The second question is answered by Theorems~\ref{th:compactness-lower-bound-untransformed} (liminf bound) and Theorem~\ref{thm:upper-bound} (limsup bound), which together establish a limit of~$\mathcal{I}_\e$ in the sense of~$\Gamma$-convergence. Here, we give a short version that combines these theorems into one statement. For convenience we collect pairs $(\rho,j)$ that satisfy the continuity equation~\eqref{eq:upscaled-CE} in a set $\CE(0,T;\mathbb{R})$; convergence in this set is defined in a distributional sense (see Definitions~\ref{def:continuity-equation} and~\ref{def:converge-in-CE}). The following theorem summarizes Theorems~\ref{th:compactness-lower-bound-untransformed} and~\ref{thm:upper-bound}.

\begin{theorem}
	[Main result]
	\label{thm:intro-main-result}
	Let $V$ satisfy Assumption~\ref{ass:V}.
	Then 
	\begin{enumerate}
	\item Sequences $(\rho_\e,j_\e)$ for which there exists a constant $C$ such that
	\begin{equation}\label{bounds:t:intro-main-result}
	\cI_\e(\rho_\e,j_\e) \leq C \qquad \text{and}\qquad E_\e(\rho_\e(0))\leq \frac C\e
	\end{equation}
	are sequentially compact in $\CE(0,T)$;
	\item Along sequences  $(\rho_\e,j_\e)$ satisfying 
\begin{equation}
\label{eq:conv-initial-data-intro}
\rho_\e(t=0) \longweakto \rho_0^\circ(\dd x) := z^\circ\delta_{x_a}(\dd x) + (1-z^\circ)\delta_{x_b}(\dd x)\qquad \text{as } \e\to0,
\end{equation}
the functional $\cI_\e$ $\Gamma$-converges to a limit $\cI_0$.
	\end{enumerate}
\end{theorem}
\noindent
In the next section we define the limit functional $\cI_0$ and show that it characterizes the limit evolution as $ z' = -z$. 

\begin{remark}
The condition~\eqref{eq:conv-initial-data-intro} can be interpreted as a well-preparedness property: it states that the initial datum converges to a measure of the same structure as the subsequent evolution (see~\eqref{structure:rho-zero-j-zero-rho} below). The bound~\eqref{bounds:t:intro-main-result} on the initial energy provides a second type of control on the initial data. 	
\end{remark}

\subsection{The limiting functional $\cI_0$}
\label{ss:limiting-functional}

Introduce  the function $S:\R^2 \to [0,\infty]$,
\begin{equation}\label{eq:S-fct}
S(a|b):=
\begin{cases}
\ds a\log\frac ab -a+b,	& 	a,b>0,\\
b,					&	a=0,\ b>0,\\
+\infty,			&	\text{otherwise}.
\end{cases}
\end{equation}
The map~$\mathcal{I}_0:\CE(0,T)\to[0,\infty]$ is defined by
\begin{equation}
\label{eqdef:I0}	
\mathcal{I}_0(\rho,j) := 2\int_0^T S(j(t)|z(t))\, \dd t,
\end{equation}
whenever
\begin{subequations}
\label{structure:rho-zero-j-zero}
\begin{align}
\label{structure:rho-zero-j-zero-rho}
& \rho(t,\dd x) = z(t)\delta_{x_a}(\dd x) + (1-z(t))\delta_{x_b}(\dd x) \text{for almost all $t$, }\\
& z(0) = z^\circ \quad\text{(see~\eqref{eq:conv-initial-data-intro}), and}\\
& \text{$j$ is piecewise constant in $x$ and given by }
 j(t,x) = j(t)\bONE_{(x_a,x_b)}(x). 
 \label{structure:rho-zero-j-zero-j}
\end{align}
\end{subequations}
Otherwise, we set $\mathcal{I}_0(\rho,j) = +\infty$.
\smallskip

\begin{lemma}[See Lemma~\ref{l:var-formulation-I0}]
 If $\cI_0(\rho,j)<\infty$, then 
\begin{enumerate}
\item the function $z$ in~\eqref{structure:rho-zero-j-zero-rho} is absolutely continuous and non-increasing, 
\item the function $j(t)$ in~\eqref{structure:rho-zero-j-zero-j} satisfies $j(t) = -z'(t)$ for almost all $t$.
\end{enumerate}
For all $(\rho,j)$, $\cI_0(\rho,j)\geq0$; if  $\cI_0(\rho,j) =0$, then $z$ satisfies $z'(t) = -z(t)$ for  all $t$. 
\end{lemma}
\noindent
The final part of this lemma allows us to characterize any limit of solutions $(\rho_\e,j_\e)$ of~\eqref{eq:upscaled-pair}. Such solutions satisfy $\cI_\e(\rho_\e,j_\e)=0$; therefore any limit  $(\rho_0,j_0)$ along a subsequence $\e_k\to0$ satisfies 
\[
0\leq \cI_0(\rho_0,j_0)\leq \liminf_{\e_k\to0} \cI_{\e_k}(\rho_\e,j_\e) = 0,
\]
and therefore $\rho_0$ has the structure~\eqref{structure:rho-zero-j-zero-rho} and the corresponding function $z$ satisfies $z'=-z$. Since the limit is unique, in fact any sequence $(\rho_{\e_\ell},j_{\e_\ell})$ converges. 
The evolution of such a function $\rho_0$  is depicted in Figure~\ref{fig:intro:limit-evolution-FP}.

\begin{remark}[$\cI_0$ does not define a gradient structure]
While the limiting \emph{equation} $z'=-z$ has multiple gradient structures (see Section~\ref{ss:conv-Kramers}), the limiting \emph{functional} $\cI_0$ does not define any gradient structure. This  example therefore is another illustration of how convergence of gradient structures is a stronger property than convergence of the equations (see~\cite{Mielke16a} for more discussion on this topic). 

To see that $\cI_0$ does not define a gradient system, at least formally, assume for the moment that  there exist $\sfE$ and $\sfR$ such that
\begin{equation}
\label{eq:tentative-GS}	
2\int_0^T S(-z'|z)\, \dd t
= \int_0^T \Bigl[ \sfR(z,z') + \sfR^*(z,-E'(z))\Bigr]\, \dd t + \sfE(z)\Big|_0^T.
\end{equation}
By taking a short-time limit we deduce that 
\[
2 S(-v|z) = \sfR(z,v) + \sfR^*(z,-E'(z)) + \sfE'(z)\cdot v
\qquad\text{for all }z,v,
\]
and by differentiating with respect to $v$ we find
\[
\mathrm D_v R(z,v) = -2\log \frac {-v}z - \sfE'(z)
\qquad\text{for all }z,v.
\]
Part of the definition of a gradient system is the requirement that $\sfR(z,\cdot)$ is minimal at $v=0$ for each $z$ (see the discussion in~\cite[p.~1296]{MielkePeletierRenger2014}), and the expression for the derivative $\mathrm D \sfR(z,v)$ above shows that this can not be the case. 
This mathematical argument backs up the more philosophical arguments in Section~\ref{ss:conv-Kramers} that that $\cI_0$ does not define a gradient system.
\end{remark}

\subsection{Discussion}
\label{ss:discussion}

\subsubsection{Main conclusions}
The main mathematical question in this paper is to understand the `fate' of a gradient structure in a limit in which this gradient structure itself must break down. What we find can be summarized as follows:
\begin{enumerate}
\item Although the energy $E_\e$ and the dissipation potentials $R_\e$ and $R_\e^*$ diverge, the single functional $\cI_\e$ that captures the Energy-Dissipation-Principle persists;
\item This functional $\cI_\e$ provides sufficient control for a proof of  compactness and $\Gamma$-convergence;
\item The limiting functional $\cI_0$ defines a `variational-evolution' system, but not a gradient system;
\item Both the EDP functional $\cI_\e$ and its limit $\cI_0$ have a clear connection to large deviations (see below). 
\end{enumerate}

Although the convergence proved in Theorem~\ref{thm:intro-main-result} is not a gradient-system convergence and  the energies $E_\e$ do not converge,  we do use  a small component of the typical gradient-system evolutionary-convergence proof. We need some control on the initial data; this is visible in the bound on $E_\e$ in~\eqref{bounds:t:intro-main-result}, which stipulates that  $E_\e(\rho_\e(t=0))$ is allowed to diverge, but not too fast. In fact, the requirement in the proof of Theorem~\ref{th:compactness-lower-bound-untransformed} is that $E_\e(\rho_\e(t=0))$ diverges more slowly than exponentially.

\subsubsection{Connection to Large-Deviation Principles}
\label{ss:connection-to-LDPs}

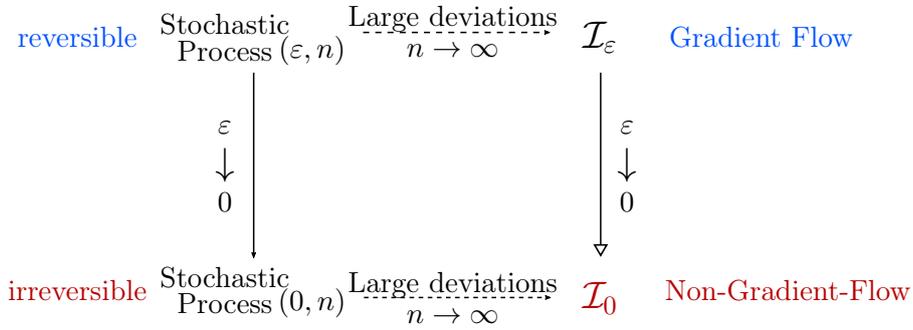
\begin{figure}[h!]
	\labellist
	\pinlabel \Large \color{black}{$\mathcal{I}_\varepsilon$} at 1300 1100
	\pinlabel \color{dark_blue}{$\text{reversible}$} at -650 1100
	\pinlabel \color{black}{$\text{Stochastic}$} at -100 1150
	\pinlabel \color{black}{$\text{Process}$} at -100 1050
	\pinlabel \color{black}{$(\varepsilon,n)$} at 225 1050
	\pinlabel \Large \color{red_one}{$\mathcal{I}_0$} at 1300 120
	\pinlabel \color{dark_blue}{$\text{Gradient Flow}$} at 1900 1100
	\pinlabel \color{red_one}{$\text{irreversible}$} at -650 150
	\pinlabel \color{black}{$\text{Stochastic}$} at -100 200
	\pinlabel \color{black}{$\text{Process}$} at -100 100
	\pinlabel \color{black}{$(0,n)$} at 225 100
	\pinlabel \color{red_one}{$\text{Non-Gradient-Flow}$} at 2000 150
	\pinlabel \color{black}{$\text{Large deviations}$} at 750 1170
	\pinlabel \color{black}{$n\to\infty$} at 750 1050
	\pinlabel \color{black}{$\text{Large deviations}$} at 750 170
	\pinlabel \color{black}{$n\to\infty$} at 750 50 
	\pinlabel \color{black}{$\varepsilon$} at -100 760
	\pinlabel \Large \color{black}{$\downarrow$} at -100 620
	\pinlabel \color{black}{$0$} at -100 480
	\pinlabel \color{black}{$\varepsilon$} at 1400 760
	\pinlabel \Large \color{black}{$\downarrow$} at 1400 620
	\pinlabel \color{black}{$0$} at 1400 480
	\endlabellist
	\centering
	\vskip5mm
	\includegraphics[scale=.1]{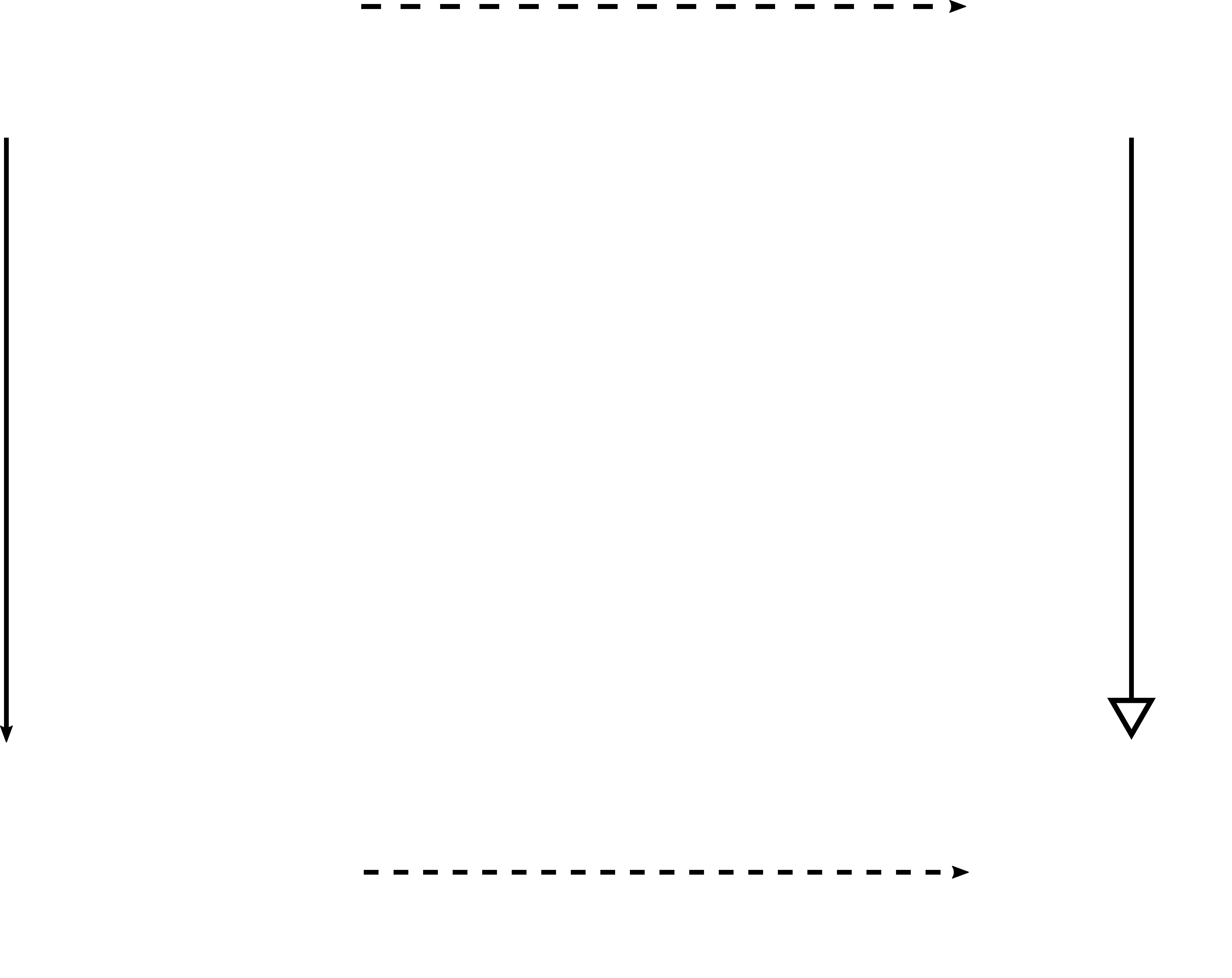}
	\caption{The top row corresponds to the empirical flux-density pairs~\eqref{eq:empirical-flux-density-pairs} stemming from i.i.d.\ copies of the reversible diffusion process~$X^\varepsilon_i(t)$ from~\eqref{eq:intro:upscaled-diffusion-process}, whose Fokker-Planck equation is~\eqref{eq:intro:upscaled-FP}. The bottom row corresponds similarly to a jump process defined on two states~$\{a,b\}$, with jumps only from~$a$ to~$b$. We prove the right arrow by Theorem~\ref{thm:intro-main-result}, and the left arrow also follows from this result. More explanation is given in the text.}
	\label{fig:commuting-diagram}
\end{figure}

Both the pre-limit functionals $\cI_\e$ and the limit functional $\cI_0$ have a clear interpretation as large-deviation rate functions of stochastic processes. In addition, the main result of this paper makes the diagram in Figure~\ref{fig:commuting-diagram} into a commuting diagram. We now explain this.

Let $X_i^\varepsilon(t)$ be independent copies of the upscaled diffusion process satisfying~\eqref{eq:intro:upscaled-diffusion-process}, and define formally 
the \emph{empirical flux-density pair}~$(\rho_{\varepsilon,n},j_{\varepsilon,n})$  by
\begin{equation}\label{eq:empirical-flux-density-pairs}
\rho_{\varepsilon,n} = \frac{1}{n}\sum_{i=1}^n \delta_{X_i^\varepsilon(t)}\quad\text{and}\quad j_{\varepsilon,n}\approx\frac{1}{n}\sum_{i=1}^n\delta_{X_i^\varepsilon(t)} \partial_t X_i^\varepsilon(t),
\end{equation}
The functional $\cI_\e$  characterizes the large deviations of $(\rho_{\varepsilon,n},j_{\varepsilon,n})$ in the limit $n\to\infty$ for fixed $\e$~\cite{DawsonGartner87,FengKurtz06} (see also~\cite[(1.3) and (2.8)]{BertiniDeSoleGabrielliJonaLasinioLandim2015})
\[
\P\Bigl[(\rho_{\varepsilon,n},j_{\varepsilon,n})\big|_{t\in[0,T]} \approx (\rho_\e,j_\e)\big|_{t\in[0,T]}\Bigr]
\ \stackrel{n\to\infty}\sim \ 
\exp \bigl(-n\cI_\e(\rho_\e,j_\e)\bigr).
\]
This is the top arrow in Figure~\ref{fig:commuting-diagram}.

The limit functional~$\mathcal{I}_0$, on the other hand, similarly characterizes the $n\to\infty$ large deviations  of flux-density pairs of $n$ independent particles jumping between two points~$x_a$ and $x_b$, with jump rates given by~$r_{a\to b}=1$ and~$r_{b\to a}=0$ (see e.g.~\cite{Renger2017,Kr17}). This is the bottom arrow in Figure~\ref{fig:commuting-diagram}.

The right-hand arrow in Figure~\ref{fig:commuting-diagram} is the main result of this paper, Theorem~\ref{thm:intro-main-result}, which establishes the $\Gamma$-convergence of $\cI_\e$ to $\cI_0$ in the limit $\e\to0$. 

In the case at hand, in which the $n$ particles constituting the stochastic processes on the left-hand side of the diagram are independent, the left-hand arrow \emph{also} follows from the results of this paper: The zero sets of $\cI_\e$ and $\cI_0$ are the forward Kolmogorov equations for the corresponding single-particle stochastic processes, for which the $\Gamma$-convergence implies convergence of solutions to solutions; in turn, this implies that the stochastic processes converge. 

In conclusion, with the results of this paper we see that the diagram of Figure~\ref{fig:commuting-diagram} commutes.

\medskip

A close connection between Gamma-convergence of rate functions and convergence of processes is observed more broadly; see~\cite{DawsonGartner94,Mariani12,BonaschiPeletier16}.

\subsubsection{Connections to chemical reactions}

There is a strong connection between the philosophy of this paper and results in the chemical literature on the appearance of irreversible chemical reactions as limits of reversible reactions, for instance using mass-action laws to describe the dynamics. Gorban, Mirkes, and Yablonksy~\cite{GorbanMirkesYablonsky13} perform an extensive analysis of such limits and the corresponding behaviour of thermodynamic potentials. Although the gradient-system description of~\eqref{eq:intro:upscaled-FP} has a clear thermodynamic interpretation (see e.g.~\cite[Ch.4--5]{Peletier2014}), the current paper is different in that the starting point is a diffusion problem, not a discrete reaction system. However, the connections between these two approaches do merit deeper study.

\subsubsection{The renormalized gradient system $(\wt E_\e,\wt R_\e)$ also does not converge}
	\label{ss:renormalized}

As we remarked in Section~\ref{ss:conv-Kramers}, the functionals $E_\e$ diverge as $\e\to0$, but the rescaled functionals $\wt E_\e := \e E_\e$ $\Gamma$-converge to a well-defined limit $\wt E_0(\rho) := \int \rho V$. It is a natural question whether switching to the rescaled gradient system $(\wt E_\e,\wt R_\e)$ might solve the singularity problems described in Section~\ref{ss:conv-Kramers}. Here the rescaled potentials are defined by
\[
\wt R_\e(\rho,j) := \e R_\e(\rho,j) 
\qquad\text{and}\qquad
\wt R_\e^* (\rho,\xi) := \e R_\e^*\Bigl(\rho,\frac1\e \xi\Bigr),
\]
and EDP-convergence of $(\wt E_\e,\wt R_\e)$ would follow from the $\Gamma$-convergence of 
\[
\wt \cD^T_\e(\rho,j) := \e \cD^T_\e(\rho,j) = \int_0^T \Bigl[\wt R_\e(\rho,j) + \wt R_\e^*\bigl(\rho,-\mathrm D \wt E_\e(\rho)\bigr)\Bigr]\, \dd t.
\]

Even if $(\wt E_\e,\wt R_\e)$ does converge in the EDP sense to $(\wt E_0,\wt R_0)$, for some dissipation potential $\wt R_0$, then this limiting gradient system $(\wt E_0,\wt R_0)$ admits a very wide class of curves as `solutions'. This can be recognized as follows. 

Let $z\in C^1([0,T])$ with $z'\leq 0$, and define $(\rho_0,j_0)$ according to~\eqref{structure:rho-zero-j-zero}; since $z'$ is bounded we have $\cI_0(\rho_0,j_0)<\infty$. By the recovery-sequence Theorem~\ref{t:recovery-sequence-untransformed} there exists a sequence $(\rho_\e,j_\e)$ converging to $(\rho_0,j_0)$ such that $\cI_\e(\rho_\e,j_\e)\to\cI_0(\rho_0,j_0)$ and $\liminf_{\e\to0} \wt E_\e(\rho_\e(t=0))\geq \wt E_0(\rho_0(t=0))$. We then calculate 
\begin{align*}
\wt\cD^T_0(\rho_0,j_0) &+ \wt E_0(\rho_0(T)) \\
&\leq \liminf_{\e\to0} \wt\cD^T_\e(\rho_\e,j_\e) + \wt E_\e(\rho_\e(T))\\
&= \liminf_{\e\to0} \e\cI_\e(\rho_\e,j_\e) + \wt E_\e(\rho_\e(0))\leq \wt E_0(\rho_0(0)).
\end{align*}
It follows that $(\rho_0,j_0)$ is a solution of the gradient system $(\wt E_0,\wt R_0)$. This shows that any decreasing function $z$ generates a solution of the gradient system $(\wt E_0,\wt R_0)$; this explains our claim that it is too degenerate to be of any use. 

\subsection{Notation}

\begin{small}
\begin{longtable}{lll}
$\CE(0,T)$	& set of pairs~$(\rho,j)$ satisfying the continuity equation &Def.~\ref{def:continuity-equation} \\
$C^{n,m}(X{\times} Y)$ & space of functions that are $n$ times differentiable on $X$\\
&\qquad  and $m$ times on $Y$ \\
$\gamma_\e$	&	invariant measure normalized to one	&	Eq.~\eqref{eq:gf:Wass-E} \\
$\gamma_\e^\ell$	&	left-normalized invariant measure	&	Sec.~\ref{ss:inv-measures}\\
$E_\e$	&	energy	&	Eq.~\eqref{eq:gf:Wass-E}\\
$\hat E_\e, \hat \cI_\e, \hat \cI_0$	&	rescaled functionals	&	Def.~\ref{def:rescaled-functionals}\\
$\cE(\mu|\nu,A)$	&	localized relative entropy	&	Sec.~\ref{ss:log-Sob-ineq}\\
$\cI_\e$	&	functional for pre-limit variational formulation&	Def.~\ref{eqdef:Ie} \\
$\cI_0$	&	functional for limit variational formulation	&	Eq.~\eqref{eqdef:I0} \\
$\hat\jmath_\e$	&	flux transformed under~$y_\e$	&	Eq.~\eqref{def:transformation-j-hat}\eqref{def:transformation-j-hat-ac}\\
$\calM(\Omega)$, $\cP(\Omega)$ & signed Borel and probability measures &Sec.~\ref{ss:prelim} \\
$\calM_{\geq 0}(\Omega)$		& non-negative Borel measures &Sec.~\ref{ss:prelim} \\
$\tmsp$, $\tmspz$ & $\tmsp = [0,T]\times \R$ and $\tmspz = [0,T]\times [-1/2,1/2]$. \\
$\cR(\mu|\nu,A)$	&	localized relative Fisher-information	&	Sec.~\ref{ss:log-Sob-ineq}\\
$\hat \rho_\e, \hat \gamma_\e^\ell$	& measures transformed under~$y_\e$	&	Eq.~\eqref{def:transformation-rho-hat}\\
$S(a|b)$	&	function in limit functional~$\cI_0$	&	Eq.~\eqref{eq:S-fct}\\
$\tau_\e$	& exponential time-scale parameter	& Eq.~\eqref{intro:eq:def-time-scale-parameter} \\
$\hat u_\e^\ell$	&	density transformed under~$y_\e$	&	Eq.~\eqref{def:transformation-u-hat}\\
$\hat u_0$	&	limit density	&	Eq.~\eqref{eq:limit-density-u0}\\
$V$	& potential/energy landscape &	Ass.~\ref{ass:V} \\
$y_\e,\phi_\e$	&	auxiliary functions	&	Sec.~\ref{ss:aux-fcns} 
\end{longtable}
\end{small}

\section{Elements of the proofs}
\label{s:elements-of-the-proof}

The proofs of compactness and $\Gamma$-convergence hinge on a number of ingredients.

\smallskip

\textit{Dual form of the functional $\cI_\e$.} The definition of $\cI_\e$ given in~\eqref{defeq:I_e2} is formal, since it only makes sense for sufficiently smooth measures $\rho_\e$ and $j_\e$. The dual formulation that arises naturally from the large-deviation context (see Section~\ref{ss:connection-to-LDPs})  solves this definition problem:
\begin{definition}
\label{def:Ie}
The functional $\cI_\e:\CE(0,T)\to [0,\infty]$ is defined by
\begin{multline}
\label{eqdef:Ie}	
	\cI_\e(\rho,j) := \sup\biggl\{ \int_0^T\int_\R\Bigl[jb - \e\tau_\e\rho \Bigl(\partial_xb - \frac1\e b V'+  \frac{1}2 b^2\Bigr)\Bigr]:\\
 b\in C^{0,1}_c([0,T]\times \R)\biggr\}.
\end{multline}
\end{definition}
\noindent
Note how  this dual form of $\cI_\e$ remains  singular in multiple ways: the factor $\e\tau_\e\rho$ is exponentially large in any region where $\rho$ has $O(1)$ mass, and it is small near the saddle $x_0$ where $\rho$ is expected to behave as $\gamma_\e$. 

The following lemma makes the connection rigorous between $\cI_\e$ and the $(E_\e,R_\e)$ gradient system. 
\begin{lemma}
\label{l:reg-rho-j}
Let $(\rho,j)\in \CE(0,T)$ satisfy $E_\e(\rho(0))<\infty$. Then
\begin{multline}
\label{eq:Ie=E+R}
\cI_\e(\rho,j) =  E_\e(\rho(T)) -  E_\e(\rho(0)) \\ 
+ 	\int_0^T \underbrace{\int_\R \biggl[ \frac1{2\e\tau_\e} \bigl|v(t,x) \bigr|^2 \rho_\e(t,\dd x)}_{R_\e(\rho_\e(t),j_\e(t))}
 + \underbrace{\frac{\e\tau_\e}2\Bigl|\partial_x \sqrt{u(t,x)} \Bigr|^2\gamma_\e(\dd x)}_{``R^*_\e\bigl(\rho_\e(t),-\mathrm DE_\e(\rho_\e(t))\bigr)\text{''}}\biggr]\, \dd t.
\end{multline}
Here the integral in~\eqref{eq:Ie=E+R} should be considered equal to $+\infty$ unless the following are satisfied:
\begin{enumerate}
\item\label{l:reg-rho-j:j} $j$ is absolutely continuous with respect to $\rho$ on $\tmsp$, with density $v := \dd j/\dd \rho$;
\item\label{l:reg-rho-j:rho} 	$\rho$ is Lebesgue-absolutely continuous on $\tmsp$, with density $u := \dd\rho/\dd \gamma_\e$;
\item \label{l:reg-rho-j:rho_x} $\partial_x u\in L^1_{\mathrm{loc}}(\tmsp)$.
\end{enumerate}
\end{lemma}
This type of reformulation is fairly standard, but we did not find an explicit proof for this case; we provide a proof in Appendix~\ref{app:pf:l:reg-rho}. In~\eqref{eq:Ie=E+R} we place   $R^*_\e\bigl(\rho_\e(t),-\mathrm DE(\rho_\e(t))\bigr)$ between quotes, since this expression  is only formal; in fact, the expression above the brace could be considered a rigorous interpretation of $R^*_\e\bigl(\rho_\e(t),-\mathrm DE(\rho_\e(t))\bigr)$.

\smallskip

\textit{Forcing concentration onto the two points $x_a$ and $x_b$.}
The starting point of the proofs of compactness and the lower bound in Theorem~\ref{th:compactness-lower-bound-untransformed} is the `fundamental estimate' of every $\Gamma$-convergence and compactness proof, 
\[
\cI_\e(\rho_\e,j_\e) \leq C.
\]
Restricting in~\eqref{eqdef:Ie} to functions~$b$ supported in~$[0,t]$, and taking into account the divergence of $E_\e(\rho_\e(0))$ as $C/\e$ (see Section~\ref{ss:conv-Kramers}) and the bound on~$\mathcal{I}_\varepsilon$, we obtain for each~$t\in[0,T]$ the estimate 
\begin{equation*}
\int_0^t \int_\R 
\frac{\e\tau_\e}2\Bigl|\partial_x \sqrt{u_\e(t,x)} \Bigr|^2\gamma_\e(\dd x)\, \dd t
+ E_\e(\rho_\e(t)) \leq \frac C\e
\end{equation*}
Since the integral is non-negative and the constant~$C$ is independent of~$t$, there are constants~$C_1,C_2$ such that for every~$t\in[0,T]$,
\begin{equation*}
  \int_0^T \int_\R 
\frac{\e\tau_\e}2\Bigl|\partial_x \sqrt{u_\e(t,x)} \Bigr|^2\gamma_\e(\dd x)\, \dd t \leq \frac{C_1}{\e} ,\; E_\e(\rho_\e(t)) \leq \frac{C_2}{\e}.
\end{equation*}
Hence
\begin{equation}
\label{ineq:fundamental-estimate-intro}
\int_0^T \int_\R 
\frac{\e\tau_\e}2\Bigl|\partial_x \sqrt{u_\e(t,x)} \Bigr|^2\gamma_\e(\dd x)\, \dd t
+ \sup_{t\in [0,T]}E_\e(\rho_\e(t)) \leq \frac C\e.
\end{equation}

The divergence of the right-hand side in~\eqref{ineq:fundamental-estimate-intro} has consequences for compactness:
\begin{enumerate}
\item Because of the growth of $V$ at $\pm\infty$, the divergence at rate $C/\e$ of $E_\e(\rho_\e(t))$ suffices to prove  tightness of $\rho_\e(t)$;
\item However, to prove concentration onto the two points $x_a$ and $x_b$, we need to use the polynomial divergence of the `Fisher information' integral that is guaranteed by~\eqref{ineq:fundamental-estimate-intro}. By applying Logarithmic Sobolev inequalities localized to each of the wells, this divergence is sufficiently slow to force concentration onto $\{x_a,x_b\}$. This does require us to assume uniform convexity of each of the two wells separately. 
\end{enumerate}
The details are given in Section~\ref{s:compactness-lower-bound}.

\smallskip

\emph{The form of the limit  functional $\cI_0$.}
One can understand how the limiting  functional $\cI_0$ appears in at least three different ways. The first is by observing that $\cI_0$ is the rate function for the Sanov large-deviation principle of a two-point jump process; see Section~\ref{ss:connection-to-LDPs} above.

The second understanding of the structure of $\cI_0$ follows from the proof of the lower bound. This bound follows from making a specific choice for the function~$b$ in the dual formulation~\eqref{eqdef:Ie}, of the form $b(t,x) = -2f(t)\delta^\e_{x_0}(x)$, where $ \delta^\e_{x_0}$ indicates an appropriately rescaled derivative of the classical committor function (see Section~\ref{ss:aux-fcns}); in the limit $\delta^\e_{x_0}$ converges to a Dirac measure at the saddle~$x_0$. With this choice we find the lower bound (Theorem~\ref{th:compactness-lower-bound-untransformed})
\begin{multline*}
\liminf_{\e\to0} \cI_\e(\rho_\e,j_\e)
\geq 2\int_0^T z(t) \Bigl[\underbrace{f'(t)}_{\text{from }jb} - \underbrace{(e^{f(t)} - 1)}_{\text{from }\e\tau_\e\rho(\dots)}\Bigr]\, \dd t + z^\circ f(0)\\ \text{for any }f\in C_b^1([0,T])\text{ with }f(T)=0.
\end{multline*}
The supremum of the right-hand side over functions $f$ equals the functional~$\cI_0$, expressed in terms of $z$. This argument is explained in detail in Section~\ref{s:compactness-lower-bound}.

The third way to understand the form of $\cI_0$ is through the construction of the recovery sequence. This sequence is obtained by first applying a spatial transformation $x 
\mapsto y = y_\e(x)$, where the mapping $y_\e$ is similar to the mapping $\hat s$ used in~\cite[Sec.~2.1]{ArnrichMielkePeletierSavareVeneroni2012}. The choice of $y_\e$ and $\tau_\e$ leads to a desingularization of $\cI_\e$, which takes the formal form
\begin{equation}
\label{eq:hatcIe-intro}
\hat \cI_\e(\hat \rho,\hat \jmath) 
 = \frac{1}{2}\int_0^T\int_\mathbb{R} \frac{1}{\hat{u}^\ell(t,y)}\big|\hat{\jmath}(t,y) + \partial_y \hat{u}^\ell(t,y)\big|^2\,\dd y\dd t.
\end{equation}
Here $\hat \rho$ and $\hat\jmath$ are transformed versions of $\rho$ and $j$ that again satisfy the continuity equation, and $\hat u^\ell$ is the density of $\hat\rho$ with respect to the `left-rescaled invariant measure'; see Section~\ref{s:upper-bound} for details. 

The remarkable aspect of this rescaling is that the expression~\eqref{eq:hatcIe-intro} no longer contains any singular parameters. The recovery sequence is constructed by solving an auxiliary PDE for $\hat u^\ell$, based on~\eqref{eq:hatcIe-intro}, which then is transformed back to a pair $(\rho_\e,j_\e)$. 

After transformation to the coordinate $y$, the left well at $x_a$ and the right well interval $[x_{b-},x_{b+}]$ (see Figure~\ref{fig:asymmetric-doublewell-potential}) are mapped to $-1/2$ and $1/2$.  From~\eqref{eq:hatcIe-intro} one then finds an alternative expression for the function $S$ of~\eqref{eq:S-fct} in terms of functions $\hat u(y)$ (see Lemma~\ref{lemma:variational-problem}):
\begin{multline*}
		S(j|z) = \frac{1}{4} \inf_u \biggl\{\int_{-1/2}^{+1/2} \frac{1}{\hat u(y)}\big|j + \hat u'(y)\big|^2\,\dd y: \quad \hat u:[-1/2,1/2]\to (0,\infty),\\
		\hat u(-1/2) = z, \text{ and } \hat u(+1/2) = 0.\biggr\}
\end{multline*}
This formula is closely related to the expression for the limiting rate functional in~\cite[Eq.~(1.30)]{ArnrichMielkePeletierSavareVeneroni2012}; see also~\cite[App.~A]{LieroMielkePeletierRenger17}.

%
%
%
%
%
\section{Rigorous setup}

\subsection{Preliminary remarks}
\label{ss:prelim}

Throughout this paper we use the following conventions and notation. We write $\tmsp$ for the time-space domain $[0,T]\times \R$. $C^{n,m}_b(\tmsp)$ is the space of functions $f:\tmsp\to\R$ that are $n$ times differentiable in $t$ and $m$ times differentiable in $x$, and these derivatives are continuous and bounded. (In the uses below we will require no mixed derivatives). $\mathcal M(\tmsp)$ and $\mathcal M(\R)$  are the sets of finite signed Borel measures on $\tmsp$ and $\R$. We will use two topologies for measures: 
\begin{itemize}
\item the \emph{narrow} topology, generated by duality with continuous and bounded functions; and
\item the \emph{wide} topology, generated by duality with continuous functions with compact support.
\end{itemize}
The sets~$\cM_{\geq0}(\R)$ and~$\mathcal P(\R)$ are the subsets of non-negative measures and probability measures with the same topology. 

For a measure $\mu\in \mathcal M(\R)$ that is absolutely continuous with respect to the Lebesgue measure, we write $\mu(dx)$ for the measure and $\mu(x)$ for the density, so that $\mu(dx) = \mu(x) dx$. The \emph{push-forward measure} of a measure $\mu\in \calM(\R)$ under a  map $T:\R\to\R$ is given by 
\[
(T_\#\mu) (A) := \mu(T^{-1}(A))\qquad\text{for all Borel sets }A\subset \Omega, 
\]
or equivalently
\[
\int_\R \varphi(y) (T_\#\mu)(\dd y) = \int _\R \varphi(T(x)) \mu(\dd x)
\qquad \text{for all Borel measurable }\varphi:\R\to\R.
\]

\subsection{Full definition of the continuity equation}
The  functionals $\cI_\e$ are defined on pairs of measures~$(\rho,j)$ satisfying the continuity equation $\partial_t\rho + \partial_x j = 0$ in the following sense.

\begin{definition}[Continuity Equation]\label{def:continuity-equation}
	We say that a pair~$(\rho(t,\cdot),j(t,\cdot))$ of time-dependent Borel measures on~$\mathbb{R}$ satisfies the continuity equation if:
	\begin{enumerate}[label=(\roman*)]
		\item For each~$t\in[0,T]$, $\rho(t,\cdot)$ is a probability measure on~$\mathbb{R}$. The map~$t\mapsto \rho(t,\cdot)\in\mathcal{P}(\mathbb{R})$ is continuous with respect to the narrow topology on~$\mathcal{P}(\mathbb{R})$.
		\item For each~$t\in[0,T]$, $j(t,\cdot)$ is a locally finite Borel measure on~$\mathbb{R}$. The map~$t\mapsto j(t,\cdot)\in \mathcal{M}(\mathbb{R})$ is measurable with respect to the wide topology on~$\mathcal{M}(\mathbb{R})$, and the joint measure on $\tmsp=[0,T]\times \R$ given by 
		\[
		\int_{t\in A} |j(t,B)|\, \dd t
		\qquad \text{for }A\subset [0,T], \ B\subset \R \text{ bounded,}
		\] is locally finite on $\tmsp$. 
		\item The pair solves $\partial_t\rho + \partial_x j = 0$ in the sense  that  for any test function $\varphi\in C_c^1(\tmsp)$ with $\varphi = 0$ at $t=T$, we have
		\begin{multline}
		\label{eq:weak-form-CE}
		 \int_0^T\int_\mathbb{R} \left[\rho(t,\dd x)\, \partial_t \varphi(t,x) 
		  +j(t,\dd x)\, \partial_x \varphi(t,x) \right]\,\dd t \\
		  + \int_\R \rho(0,\dd x) \varphi(0,x)= 0.
		\end{multline}
	\end{enumerate}
	We denote by~$\CE(0,T)$ the set of all pairs~$(\rho,j)$ satisfying the continuity equation.\qed
\end{definition}

This definition gives rise to a corresponding concept of convergence.

\begin{definition}[Convergence in $\CE$]
\label{def:converge-in-CE} 
We say that $(\rho_\e,j_\e)$ converges in $\CE(0,T)$ to $(\rho_0,j_0)\in \CE(0,T)$ if
\begin{enumerate}
\item $\rho_\e(0,\cdot)$ converges narrowly to $\rho_0(0,\cdot)$ on $\R$;
\item $\rho_\e$ converges narrowly to $\rho_0$ on $\tmsp$;
\item for all $\varphi\in C_c^1(\tmsp)$ with $\varphi = 0$ at $t=T$, 
\begin{equation}
\label{eq:conv-weak-form-j}
 \lim_{\e\to0}
 \int_0^T\int_\mathbb{R} j_\e(t,\dd x)\, \partial_x \varphi(t,x)\,\dd t =\int_0^T\int_\mathbb{R} j_0(t,\dd x)\, \partial_x \varphi(t,x)\,\dd t.
 \end{equation}
\end{enumerate}
Note that then the identity~\eqref{eq:weak-form-CE} for $(\rho_\e,j_\e)$ passes to the limit. 
\end{definition}

\begin{remark}[The convergence arises from a metric]
The narrow convergence of $\rho_\e$ is generated by well-known  metrics such as the L\'evy-Prokhorov or Bounded-Lipschitz metrics~\cite[Sec.~11.3]{Dudley04}. Since $C_c(\R)$ is separable, a metric can also be constructed for the wide topology in the usual way. 
\end{remark}

\begin{remark}[Other definitions of the continuity equation]
Definition~\ref{def:continuity-equation} is weaker than the common continuity-equation concept for Wasserstein-continuous curves~\cite[Sec.~8.1]{AmbrosioGigliSavare2008}, in which $j$ is of the form $j=v\rho$ with $\iint\rho|v|^2<\infty$. While for curves $(\rho_\e,j_\e)$  with $\cI_\e(\rho_\e,j_\e)<\infty$ the flux $j_\e$ indeed has this structure  (see~\eqref{eq:Ie=E+R}), in the limit $j$ no longer is absolutely continuous with respect to $\rho$ (see the characterization of finite $\cI_0$ in~\eqref{structure:rho-zero-j-zero-j}). 
	
In addition, we choose to incorporate the initial datum in the distributional definition of the continuity equation~\eqref{eq:weak-form-CE}, as is common in the theory of parabolic equations with weak time regularity (see e.g.~\cite[Sec.~I.3]{LSU}). The explicit initial datum is used below in proving that the limit of $\rho_\e$  connects continuously to the limiting initial datum; see steps 3 and 4 of the proof of Theorem~\ref{th:compactness-lower-bound-untransformed}. 
\end{remark}

\begin{remark}[Different topologies for $\rho$ and $j$]
\label{r:rho-widely-ct}
It may seem odd that for~$\rho$ we require narrow continuity  in Definition~\ref{def:continuity-equation} and narrow convergence in Definition~\ref{def:converge-in-CE}, but for~$j$ we require only wide convergence  in Definition~\ref{def:converge-in-CE}.  

This difference arises from the following considerations. For $j_\e$, convergence of the weak form~\eqref{eq:conv-weak-form-j} is what we obtain in the proof of the compactness (Theorem~\ref{th:compactness-lower-bound-untransformed}) and of the convergence of the recovery sequence (Theorem~\ref{t:recovery-sequence-untransformed}). In both cases it is not clear whether $j_\e$ converges in a stronger manner than widely.

For $\rho_\e$, however, it is important that in the limit no mass is lost at infinity; this requires narrow convergence. In the setup above, this narrow convergence follows from the wide convergence of $j_\e$ on $[0,T]\times \R$, which also implies wide convergence for $\rho_\e$ on the same space; since the limit $\rho(t,\cdot)$ is again required to be a probability measure for all $t$, no mass escapes to infinity, and the convergence of $\rho_\e$ in fact is narrow.

The narrow continuity of $t\mapsto \rho(t,\cdot)$ in Definition~\ref{def:continuity-equation} follows from the the conditions on $j$: the local bounds on $j$ imply wide continuity of $\rho$, and the requirement that $\rho$ is a probability measure at all $t$ upgrades this continuity to narrow.
\end{remark}

\section{Compactness}
\label{s:compactness-lower-bound}
%
%
%

The limit~$\e\to 0$ is accompanied by 
the concentration of $\rho_\e$ onto the two minima of the wells, at $x_a$ and $x_b$.
This concentration is essential for the further analysis of the functionals $\cI_\e$ and their $\Gamma$-limits; if $\rho_\e$ would maintain mass at other points in $\R$, then the main statement and the corresponding analysis of the functionals~$\cI_\e$ both would fail. 

In the case of a potential $V$ with wells of equal depth (as in~\cite{ArnrichMielkePeletierSavareVeneroni2012,LieroMielkePeletierRenger17}), a constant bound on the initial energy $E_\e(\rho_\e(t=0))$ leads to a similar bound on later energies $E_\e(\rho_e(t))$, which in turn leads to concentration onto $\{x_a,x_b\}$. In the unequal-well case of this paper, as we discussed in the introduction, we are forced to allow for divergent $E_\e$; consequently the concentration onto $\{x_a,x_b\}$  has to come from different arguments. 

Here we choose to obtain this concentration from the `Fisher-information' or `local-slope term'; this is the second term in $\cD_\e^T$ in~\eqref{eq:EDP-formulation}, or equivalently the second half of the integral in~\eqref{eq:Ie=E+R}. This requires imposing conditions on the convexity of the wells, which we do in part~\ref{ass:V:5} of the following  set of assumptions on~$V$. 




\begin{figure}[ht]
	\labellist
	\pinlabel \labelsize $x$ at 2750 690
	\pinlabel \labelsize $V(x)$ at 2530 1380
	\pinlabel \labelsize $x_a$ at 500 690
	\pinlabel \labelsize $x_{c\ell}$ at 770 690
	\pinlabel \labelsize $x_0$ at 1080 690
	\pinlabel \labelsize $x_{cr}$ at 1320 690
	\pinlabel \labelsize $x_{b-}$ at 1530 925
	\pinlabel \labelsize $x_b$ at 1800 800
	\pinlabel \labelsize $x_{b+}$ at 2130 925
	\pinlabel \labelsize {\color{blue_fig}{$V''\geq \alpha >0$}} at 350 1660
	\pinlabel \labelsize {\color{blue_fig}{$V''\geq \alpha >0$}} at 1900 1660
	\pinlabel \labelsize {\color{red_fig}{$B_a$}} at 600 380
	\pinlabel \labelsize {\color{red_fig}{$B_0$}} at 1020 470
	\pinlabel \labelsize {\color{red_fig}{$B_b$}} at 1650 -50
	\endlabellist
	\centering
	\includegraphics[scale=.1]{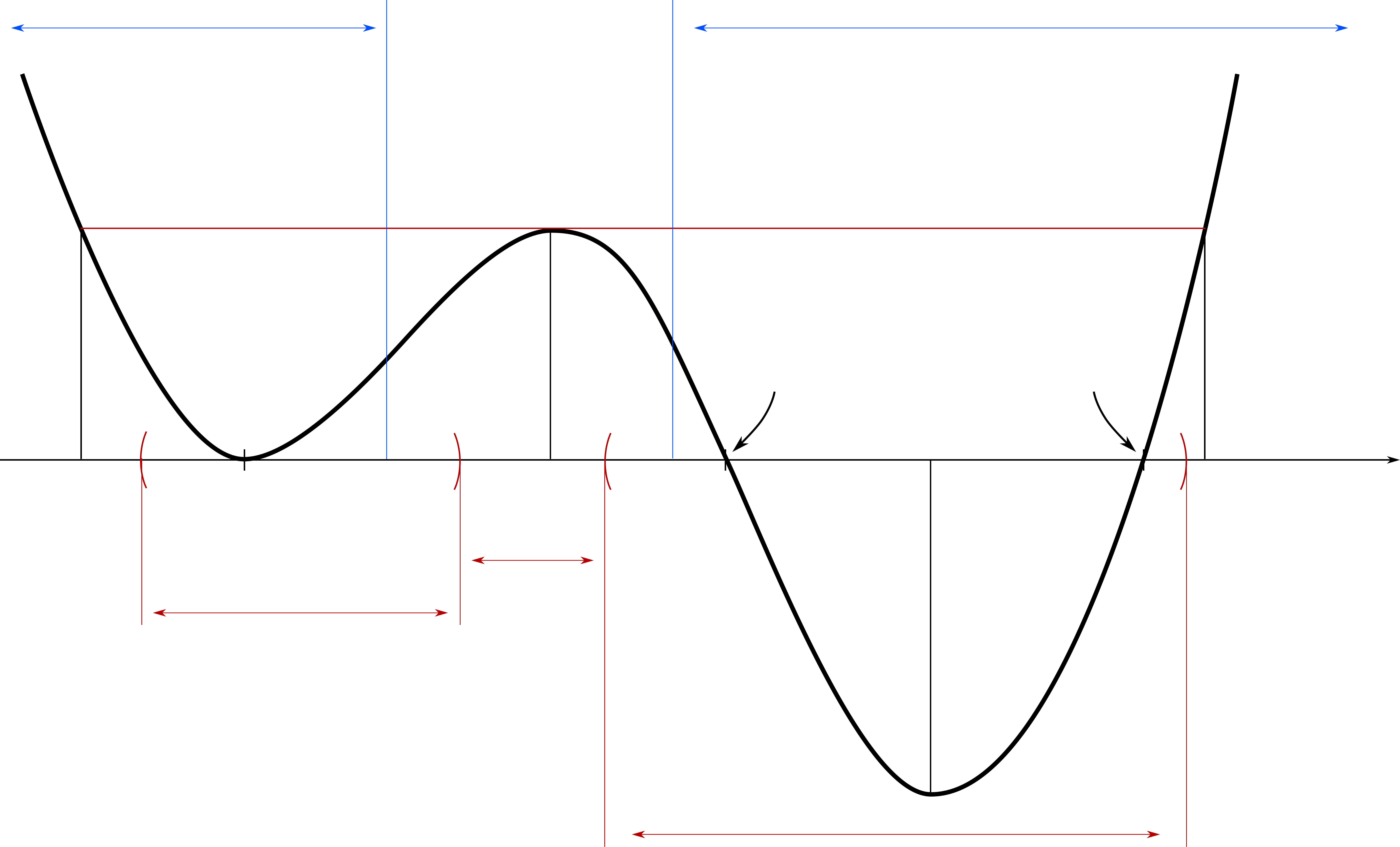}
\caption{Illustration of Assumption~\ref{ass:V}. }
\label{fig:detailed-V}
\end{figure}
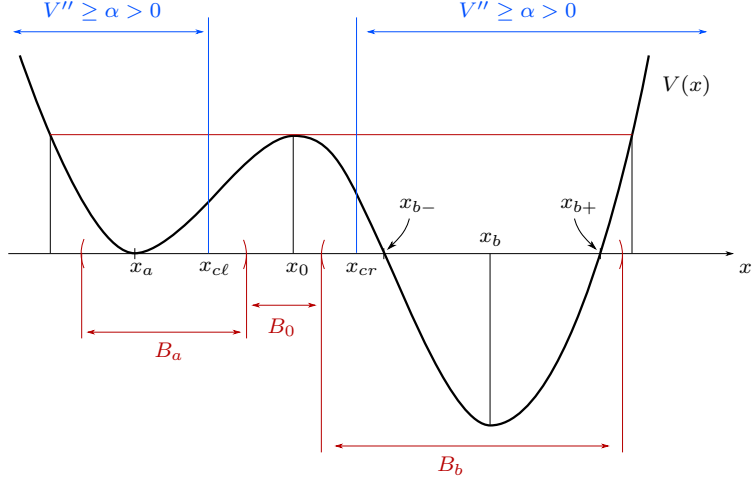

\begin{assumption}
\label{ass:V}
Let $V\in C^2(\R)$ and let the special $x$-values
\[
-\infty < x_a<x_{c\ell}<x_0<x_{cr}< x_{b-}<x_b<x_{b+}<\infty
\]
satisfy the following:
\begin{enumerate}
\item \label{ass:V:1} \emph{Two wells, the left well at value zero:} $\{V\leq0\} = \{x_a\}\cup [x_{b-},x_{b+}]$; 
	\item \label{ass:V:2}\emph{$x_b$ is the bottom of the right well:} $V(x_b) = \min_\R V < 0$;
	\item \label{ass:V:3}\emph{$x_0$ is the saddle, and the intermediate range lies below it:} $V(x)\leq V(x_0)$ for $x_a<x<x_b$, with $V(x)< V(x_0)$ unless $x=x_0$;
	\item \label{ass:V:4}\emph{The saddle is non-degenerate:} $V''(x_0)<0$;
	\item \label{ass:V:5}\emph{Uniform convexity away from the saddle:} there exist $A, \alpha>0$ such that $A \geq V''\geq \alpha> 0$ on $(-\infty, x_{c\ell}]$ and $[x_{cr},\infty)$.
\end{enumerate}
We also choose two open intervals $B_a$ and $B_b$ containing $x_a$ and $[x_{b-},x_{b+}]$, respectively, and such that $\sup_{B_a\cup B_b}V<V(x_0)$. The set $B_0$ is defined as the set separating $B_a$ and $B_b$. Figure~\ref{fig:detailed-V} illustrates these features. 
\end{assumption}


Assumptions~\ref{ass:V:1}--\ref{ass:V:4} encode the basic geometry of a two-well potential with unequal wells. Condition number~\ref{ass:V:5} is added to rule out concentration at different points than $x_a$ and $x_b$. The following two examples illustrate how concentration at different points may happen if this convexity condition is not imposed. 
%
%
%
%
%
%
\medskip

\begin{figure}[ht]
	\labellist
	\pinlabel \labelsize $x$ at 2670 1320
	\pinlabel \labelsize $V(x)$ at 2530 2000
	\pinlabel \labelsize $x_a$ at 500 1320
	\pinlabel \labelsize $x_0$ at 750 1320
	\pinlabel \labelsize $x_d$ at 1180 1450
	\pinlabel \labelsize $x_b$ at 2000 1450
	\pinlabel \labelsize {\color{red_fig}{$\text{barrier}$}} at -50 1600
	\pinlabel \labelsize {\color{red_fig}{$x_a \rightsquigarrow x_d$}} at -50 1500
	\pinlabel \labelsize {\color{red_fig}{$\text{barrier}$}} at 2600 1000
	\pinlabel \labelsize {\color{red_fig}{$x_d \rightsquigarrow x_b$}} at 2600 900
	\endlabellist
	\centering
	\includegraphics[scale=.1]{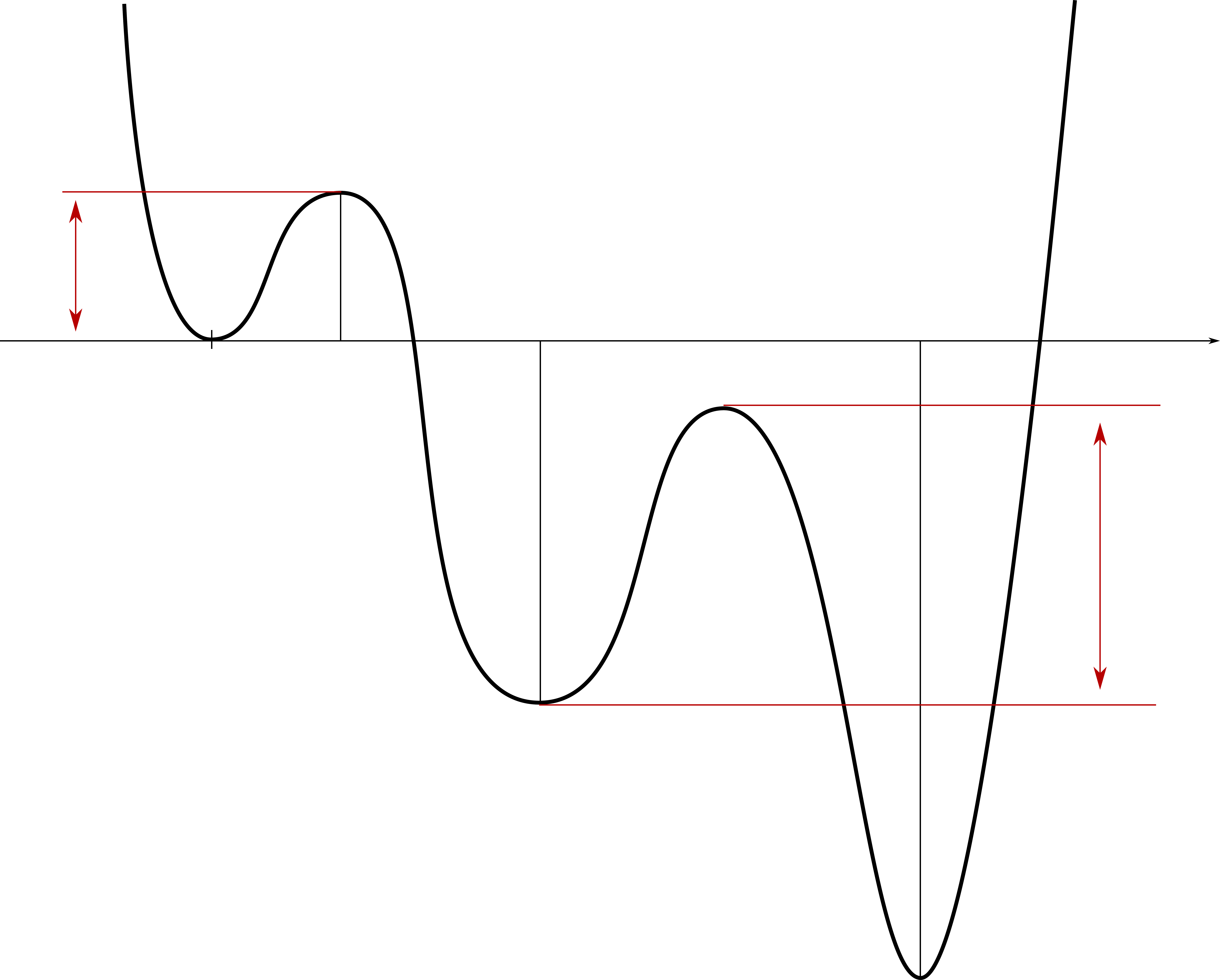}
	\caption{Example of a potential $V$ that is excluded by Assumption~\ref{ass:V} (failure of `type I' in the text).  If the deeper well has internal energy barriers that are larger than the barrier $V(x_0)-V(x_a)$ for escape from $x_a$, then mass may collect in  intermediary valleys instead of at $x_b$. In this example the mass will concentrate onto $\{ x_a, x_d, x_b\}$, with no mass moving from $x_d$ to $x_b$.}
	\label{fig:counterexample-I} 
\end{figure}
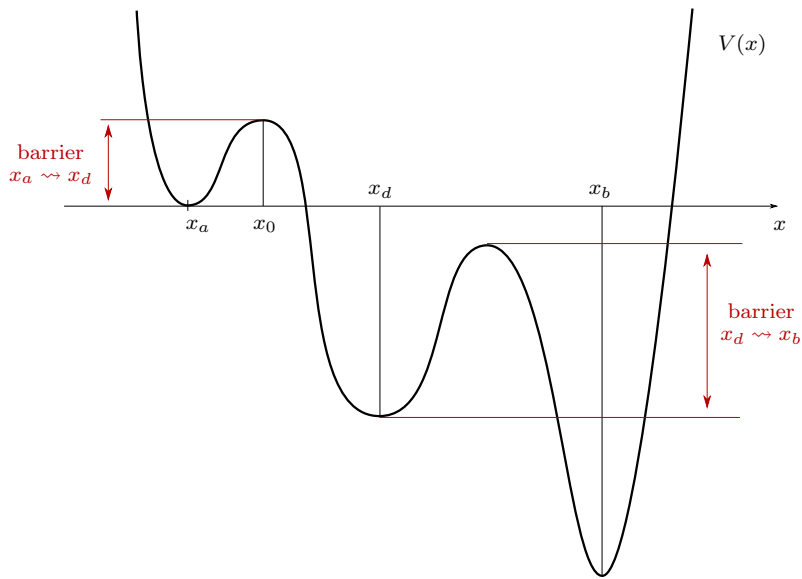

\textit{Failure type I: A hilly right well.} 
Since the energy barrier is lower for transitions from left to right than vice versa, it is natural to assume that in the limit all mass travels from left to right. Indeed, this is true under weak assumptions, but the mass that arrives in the right well $[x_{b-},x_{b+}]$ need not all end up in $x_b$. Figure~\ref{fig:counterexample-I} shows why: if the right well has a `sub-well' (say $x_d$) such that the transition $x_d\rightsquigarrow x_b$ has a higher energy barrier than the transition $x_a\rightsquigarrow x_d$, then the mass leaving $x_a$ will be held back at $x_d$, with further transitions to $x_b$ happening at an exponentially longer time scale. If we start with all mass concentrated at $x_a$, then the limiting evolution will be concentrated on $\{x_a,x_d\}$ instead of on $\{x_a,x_b\}$.

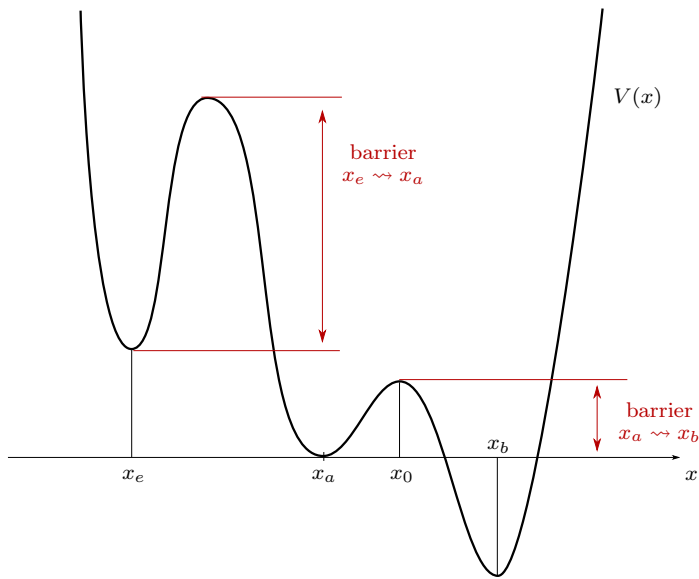
\begin{figure}[htb]
	\labellist
	\pinlabel \labelsize $x$ at 2550 380
	\pinlabel \labelsize $V(x)$ at 2350 1800
	\pinlabel \labelsize $x_e$ at 470 380
	\pinlabel \labelsize $x_a$ at 1180 380
	\pinlabel \labelsize $x_0$ at 1470 380
	\pinlabel \labelsize $x_b$ at 1830 490
	\pinlabel \labelsize {\color{red_fig}{$\text{barrier}$}} at 1400 1600
	\pinlabel \labelsize {\color{red_fig}{$x_e \rightsquigarrow x_a$}} at 1400 1500
	\pinlabel \labelsize {\color{red_fig}{$\text{barrier}$}} at 2430 630
	\pinlabel \labelsize {\color{red_fig}{$x_a \rightsquigarrow x_b$}} at 2430 530
	\endlabellist
	\centering
	\includegraphics[scale=.1]{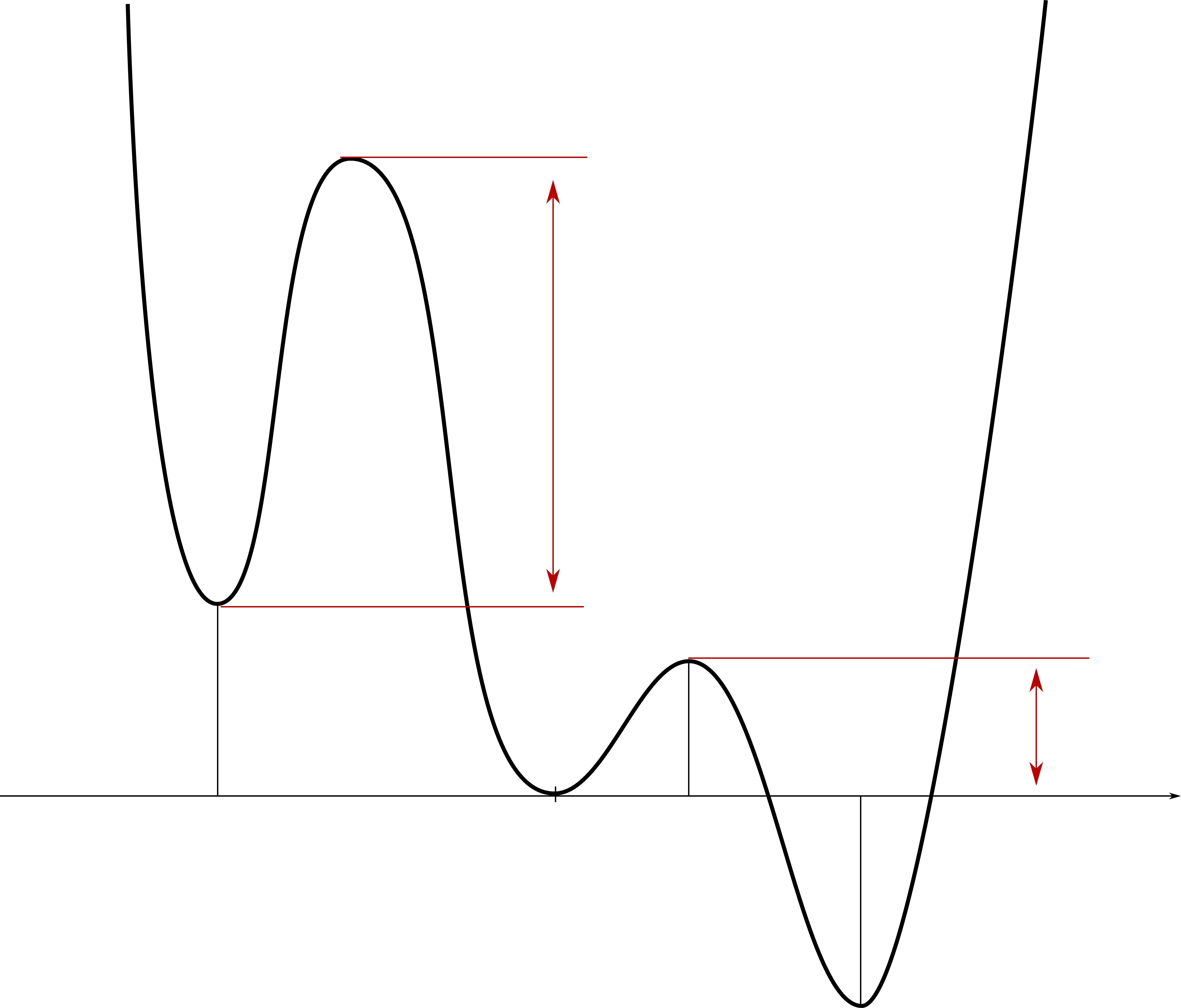}
	\caption{Second example of a potential $V$ that is excluded by Assumption~\ref{ass:V} (failure of `type II' in the text). If energy barriers exist outside of the range $[x_a,x_b]$ that are larger than the barrier $V(x_0)-V(x_a)$ of the escape from $x_a$, then these will prevent mass from moving to $x_a$ (and then also from transitioning to $x_b$). In this example the mass will concentrate onto $\{x_e, x_a, x_b\}$, with no mass moving from $x_e$ to $x_a$ or $x_b$.}
	\label{fig:counterexample-II} 
\end{figure}

\medskip
\textit{Failure type II: Hills at high energy levels.}
Something similar can happen in the `wings' of the energy landscape, as illustrated by Figure~\ref{fig:counterexample-II}. If valleys exist outside of the region $\{x:V(x)<V(x_0)\}$ with energy barriers larger than the $x_a\rightsquigarrow x_b$ barrier, then the slowness of transitions between such valleys  again will prevent concentration into the sub-zero zone $\{x:V(x)\leq 0\}$. 

\subsection{Logarithmic Sobolev inequalties}
\label{ss:log-Sob-ineq}
We use logarithmic Sobolev inequalities to capitalize on the uniform convexity bounds in part~\ref{ass:V:5} of Assumption~\ref{ass:V}. Such inequalities are usually formulated for reference measures 	 with unit mass, but in our case it will be convenient to generalize  to all finite positive measures, and also allow for localization to subsets of~$\R$. 

For $A\subset \R$ and $\mu,\nu\in \cM_{\geq0}(\R)$, we set 
\begin{align*}
\cE(\mu|\nu,A) &:= \begin{cases}
\ds\int_A f\log f \;\dd \nu  
  - \mu(A) \log \frac {\mu(A)}{\nu(A)}
  & \text{if }\mu\ll \nu,\ \mu = f\nu\\
+\infty & \text{otherwise},
\end{cases}\\
\cR(\mu|\nu,A) &:= \begin{cases} 
\ds2 \int_A \bigl|\partial_x  \sqrt {f\,}\,\bigr|^2 \, \dd \nu,
    & \text{if }\mu\ll \nu,\ \mu = f\nu\\
+\infty&\text{otherwise}.
 \end{cases}
\end{align*}

With these definitions, the energy $E_\e$ and the `slope' $R_\e^*(\rho,-\mathrm D E_\e(\rho))$ (see~\eqref{eq:gf:Wass-E}) can be written as 
\begin{equation}
\label{reform:ER-in-cER}	
E_\e(\rho) = \cE(\rho|\gamma_\e,\R) 
\qquad\text{and}\qquad
R_\e^*(\rho,-\mathrm D E_\e(\rho)) \stackrel{(*)}= \e\tau_\e \cR(\rho|\gamma_\e,\R).
\end{equation}
The identity $(*)$ can also be seen as a rigorous definition of the left-hand side $R_\e^*(\rho,-\mathrm D E_\e(\rho))$ in terms of the right-hand side $\cR(\rho|\gamma_\e,\R)$: this right-hand side is well defined for all $\rho$, and in addition Lemma~\ref{l:reg-rho-j} shows that this is the term that appears in the reformulation of $\cI_\e$ in gradient-system form.

Note that the functions $\cE$ and $\cR$ are  $(1,0)$-homogeneous in the pair $(\mu,\nu)$, i.e.\ for each $\mu,\nu\in \cM_{\geq0}(\R)$ and $a,b>0$, 
\[
\cE(a\mu|b\nu,A) = a\cE(\mu|\nu,A)\quad \text{and}\quad\cR(a\mu|b\nu,A) = a\cR(\mu|\nu,A).
\]
The following Lemma generalizes classical Logarithmic Sobolev inequalities based on uniform convexity bounds  to the homogeneous functionals~$\cE$ and~$\cR$ and the restriction to subsets $A\subset\R$.

\begin{lemma}[Logarithmic Sobolev inequality]
\label{l:LSI}
Let $A\subset \R$ be an interval. 
If $W\in C^2(A)$ with $W''\geq \alpha>0$ on $A$, then
\begin{equation}
\label{ineq:LSI}
\alpha \cE(\mu|e^{-W}\dd x, A) \leq \cR(\mu|e^{-W}\dd x, A)
\qquad\text{for all }\mu\in \cM_{\geq0}(A).
\end{equation}
\end{lemma}

\begin{proof}
By e.g.~\cite[Cor.~5.7.2]{BakryGentilLedoux2013} or~\cite[Cor.~1]{Cordero-Erausquin02}, if $W\in C^2(\R)$ with $W''\geq \alpha>0$ on $\R$, then the inequality~\eqref{ineq:LSI} holds for $A=\R$ and for all $\mu\in \mathcal P(\R)$. By the homogeneity of $\cE$ and $\cR$ the same applies to all $\mu\in\cM_{\geq0}(\R)$.

To generalize to the case of $A\subsetneqq \R$ and a given potential $W\in C^2(A)$ with $W''\geq \alpha$ on $A$, first smoothly extend $W$ to the whole of $\R$ in such a way that $W''\geq \alpha$ on $\R$ and  $\int_\R e^{-W}<\infty$. Next define the sequence of $C^2$ potentials
\[
W_k(x) := W(x) + k\dist(x,A)^4 \qquad\text{for }x\in \R.
\]
As $k\to\infty$ the measures $e^{-W_k(x)}\dd x$ converge narrowly on $\R$ to $e^{-W(x)}\bONE_{A}(x)\dd x$. Each $W_k$ satisfies $W_k''\geq \alpha$ on $\R$, and it follows that for any $\mu\in \cM_{\geq0}(\R)$ with $\mu(\R\setminus A) = 0$, 
\begin{align*}
\alpha \cE(\mu|e^{-W}\dd x, A) &= \alpha \cE(\mu|e^{-W}\bONE_{A}\dd x, \R)
= \lim_{k\to\infty} \alpha \cE(\mu|e^{-W_k}\dd x, \R) \\
&\leftstackrel{\eqref{ineq:LSI} \text{ on }\R}\leq \lim_{k\to\infty} \cR(\mu|e^{-W_k}\dd x, \R) =  \cR(\mu|e^{-W}\dd x, A).
\end{align*}
This proves the claim~\eqref{ineq:LSI}.
\end{proof}

Bounds on the entropy give rise to concentration estimates of the underlying measure.

\begin{lemma}[Concentration estimates based on $\cE$]
\label{l:concentration-LSI}
Let $A_1\subset A_2\subset \R$, and let $\mu,\nu\in \cM_{\geq0}(A_2)$ with $\nu(A_1)>0$. Then
\begin{equation}
\label{ineq:concentration-LSI}	
\mu(A_1) \leq \frac{\cE(\mu|\nu,A_2) + \mu(A_2)}{\log \bigl( \nu(A_2)/\nu(A_1)\bigr)}.
\end{equation}
\end{lemma}

\begin{proof}
By homogeneity of $\cE$ it is sufficient to prove the inequality for the case $\mu(A_2) = 	\nu(A_2) = 1$. We can also assume that $\mu\ll \nu$, and we set $\mu = f\nu$. 

Applying  Young's inequality with the dual pair $\eta(s) = s\log s -s + 1$ and $\eta^*(t) = e^t-1$, we find for any $a>0$ that
\begin{align*}
\mu(A_1) = \frac1a \int _{A_1}  f\, a \, \dd \nu
&\leq \frac1a \int_{A_1} \eta(f)\, \dd \nu 
  + \frac1a \int_{A_1} \bigl(e^a-1\bigr) \, \dd \nu\\
&\leq\frac1a \int_{A_2} \eta(f)\, \dd \nu 
  + \frac{e^a}a \nu(A_1).
\end{align*} 
Choosing $a = |\log \nu(A_1)|  = -\log\nu(A_1)$ we find 
\[
\mu(A_1) \leq \frac1{|\log \nu(A_1)|} \bigl(\cE(\mu|\nu,A_2) + 1\bigr),
\]
which is~\eqref{ineq:concentration-LSI} for the case $\mu(A_2) = 	\nu(A_2) = 1$.
\end{proof}

\subsection{Invariant measures and their normalizations}
\label{ss:inv-measures}
In the introduction we defined the invariant measure 
\[
\gamma_\e(\dd x) := \frac1{Z_\e} e^{-V(x)/\e}\, \dd x,
\qquad\text{with}\qquad
Z_\e := \int_\R e^{-V(x)/\e}\, \dd x.
\]
The measure $\gamma_\e$ is normalized in the usual manner, and is therefore a probability measure on $\R$. Since $V$ has a single global minimum at $x_b$, the measures $\gamma_\e$ converge to $\delta_{x_b}$; therefore the mass of $\gamma_\e$ around $x_a$ vanishes. It will also be useful to have a differently normalized measure $\gamma_\e^\ell$ in which the mass around $x_a$ does not vanish. For this reason we also define the \emph{left-normalized} measures $\gamma_\e^\ell$ by 
\begin{equation*}
	\gamma_\e^\ell(\dd x) := \frac1{Z_\e^\ell} e^{-V(x)/\e}\,\dd x, \quad\text{with}\quad  
	Z^\ell_\e := \int_{-\infty}^{x_0} e^{-V(x)/\e}\, \dd x.
	\end{equation*}

\begin{figure}[h!]
	\labellist
	\pinlabel \labelsize {\begin{turn}{54.52}left-normalized\end{turn}} at 1300 1500
	\pinlabel \labelsize {\begin{turn}{54.52}$\e\to 0$\end{turn}} at 1550 1300
	\pinlabel \labelsize {\begin{turn}{-54.52}fully normalized\end{turn}} at 3200 1500
	\pinlabel \labelsize {\begin{turn}{-54.52}$\e\to 0$\end{turn}} at 3000 1300
	\pinlabel \labelsize $x$ at 3500 3400
	\pinlabel \labelsize $V(x)$ at 3400 4000
	\pinlabel \labelsize $x_a$ at 1300 3450
	\pinlabel \labelsize $x_0$ at 1900 3450
	\pinlabel \labelsize $x_b$ at 2630 3630
	\pinlabel \labelsize $x_{b-}$ at 2300 3750
	\pinlabel \labelsize $x_{b+}$ at 2960 3750
	\pinlabel \labelsize $\gamma_\varepsilon(x)\sim\exp\{-V(x)/\e\}$ at 3500 2400
	\pinlabel \labelsize $+\infty$ at 1690 670
	\pinlabel \labelsize $\delta_{x_a}$ at 300 670
	\pinlabel \labelsize $x_a$ at 380 -60
	\pinlabel \labelsize $x_0$ at 1100 -60
	\pinlabel \labelsize $x_b$ at 1730 -60
	\pinlabel \labelsize $x_{b-}$ at 1530 -60
	\pinlabel \labelsize $x_{b+}$ at 1950 -60
	\pinlabel \labelsize $\delta_{x_b}$ at 3940 670
	\pinlabel \labelsize $x_b$ at 4050 -60 
	\endlabellist
	\centering
	\includegraphics[scale=.08]{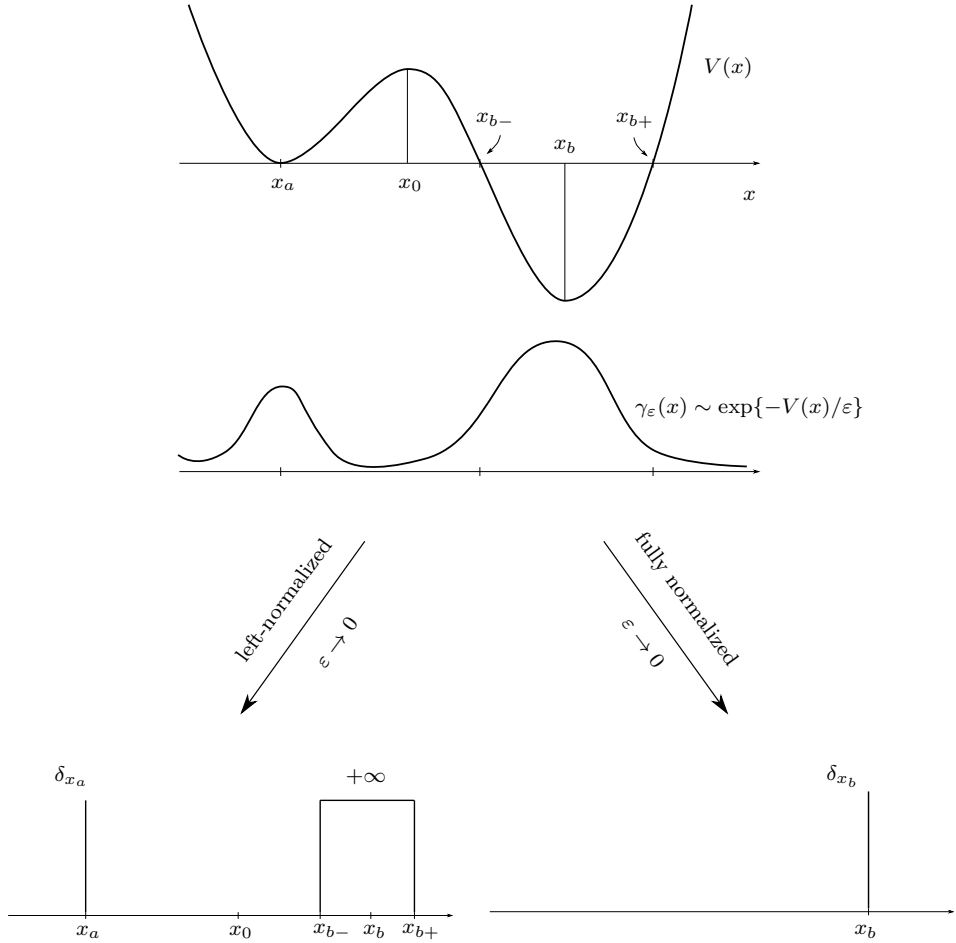}
	\vspace*{5mm}
	\caption{Behavior of the left-normalized invariant measure $\gamma_\e^\ell$ and the fully-normalized measure $\gamma_\e$ for small values of~$\varepsilon$.}
	\label{fig:left-and-fully-normalized}
\end{figure}

Figure~\ref{fig:left-and-fully-normalized} illustrates the behaviour of $\gamma_\e^\ell$ and $\gamma_\e$ as $\e\to0$. 
The following lemma characterizes some of their behaviour in precise form.

\begin{lemma}
\label{l:basic-props-gamma}
Let $V$ satisfy Assumption~\ref{ass:V}.
\begin{enumerate}
\item \label{l:basic-props-gamma:scaling-Z}
$\gamma_\e$ and $\gamma^\ell_\e$ are well-defined, and in the limit  $\e\to0$,
\begin{equation}
\label{char:asymp-Ze}
Z_\e=[1+o(1)] \sqrt{\frac{2\pi\e}{V''(x_b)}} e^{-V(x_b)/\e}, \qquad 
Z_\e^\ell=[1+o(1)] \sqrt{\frac{2\pi\e}{V''(x_a)}} .
\end{equation}
\item \label{l:basic-props-gamma:conv-prop-2-Z}
If $\tilde x>x_0$ and $V<V(x_0)$ on $(x_0,\tilde x]$, then 
\[
\frac{Z_\e^\ell}{\e\tau_\e} \int_{x_0}^{\tilde x} e^{V/\e} \longrightarrow \frac12 \qquad\text{as }\e\to0.
\]
\item \label{l:basic-props-gamma:concentration}
For any~$\delta > 0$, $	\lim_{\e\to 0}\gamma_\e^\ell(\{V>\delta\}) = 0$.
\item \label{l:basic-props-gamma:limit}
For any $x_a<c<x_0< x_{b-}< d$, the sequence $\gamma^\ell_\e\lfloor(-\infty,c)$ converges as measures to $\delta_{x_a}$, and $\gamma^\ell_\e((c,d)) \to \infty$. 
\end{enumerate}
\end{lemma}

Part~\ref{l:basic-props-gamma:concentration} above expresses the property that  the left-norm\-ali\-zed measures concentrate in the limit~$\e\to 0$ onto the set $\{V\leq 0\}=\{x_a\}\cup [x_{b-},x_{b+}]$.
Part~\ref{l:basic-props-gamma:limit} expresses the fact that the `left-hand' part of $\gamma^\ell_\e$ has a well-behaved limit $\delta_{x_a}$, while the right-hand part of $\gamma^\ell_\e$ has unbounded mass.

\begin{proof}
For part~\ref{l:basic-props-gamma:scaling-Z}, the superquadratic growth of $V $ towards $\pm\infty$ that follows from uniform convexity implies that $Z_\e$ and $Z^\ell_\e$ are finite for each $\e$; the scaling of $Z_\e$ and $Z_\e^\ell$ then follow directly from Laplace's method (Lemma~\ref{lemma:watson}). The same holds for part~\ref{l:basic-props-gamma:conv-prop-2-Z}, and the convergence of $\gamma^\ell_\e\lfloor(-\infty,c)$  to $\delta_{x_a}$ (part~\ref{l:basic-props-gamma:limit}).

For part~\ref{l:basic-props-gamma:concentration}, we estimate using the superquadratic growth of $V$ that 
\begin{align*}
\frac1{Z_\e^\ell} \int_{V> \delta} e^{-V/\e} \leq 
 \frac1{Z_\e^\ell} \int_\R \exp\Bigl(-\frac1\e\big(\delta\vee C(x^2-1)\bigr)\Bigr) 
 \stackrel{\e\to0}\longrightarrow0.
\end{align*}
This proves the claim.

Finally, to show that   $\gamma^\ell_\e((c,d))\to\infty$ (part~\ref{l:basic-props-gamma:limit}), note that $V(x)\leq -\mu<0$ for some constant $\mu>0$ on an open interval $(x_{b-}+\delta,x_{b-}+2\delta)\subset (c,d)$; from this the divergence follows. 
\end{proof}

\subsection{Auxiliary functions $\phi_\e$ and $y_\e$}
\label{ss:aux-fcns}

To desingularize the functional $\cI_\e$  we will need an auxiliary function $\phi_\e$ that is adapted to the singular structure of this system and distinguishes the two wells, in the sense of having constant, but different, values there. For the recovery sequence we will need a related function $y_\e$, and we define it here at the same time, and  study the properties of $\phi_\e$ and $y_\e$ together. 

 Fix two smooth functions $\chi_a,\chi_b\in C^\infty_c(\R)$ with $\chi_{a,b}\geq0$, $\supp \chi_a\subset B_a$ and $\supp \chi_b\subset B_b$, and $\chi_a(x_a) = 1 = \chi_b(x_b)$. Set 
\[
\mu_\e 
:= \frac{e^{-V/\e}\chi_a}{\int e^{-V/\e}\chi_a}
- \frac{e^{-V/\e}\chi_b}{\int e^{-V/\e}\chi_b} 
\qquad\text{and}\qquad
M_\e(x) := \int_{-\infty}^x \mu_\e.
\]
The function $M_\e$  has the following properties:
\begin{enumerate}
	\item $0\leq M_\e\leq 1$;
	\item $M_\e$ is equal to $1$ on $B_0$ and equal to zero outside of $B_a\cup B_0\cup B_b$, and converges in $L^1$ to $\bONE_{[x_a,x_b]}$;
\end{enumerate}
Define $\phi_\e\in C^2_b(\R)$ and $y_\e\in C^2(\R)$ by
\begin{align}
\label{eqdef:phi_e}	
\phi_\e(x) &:= \frac{Z^\ell_\e}{\e\tau_\e} \int_{x_0}^x e^{V(\xi)/\e}M_\e(\xi)\, \dd \xi\\
y_\e(x) &:= \frac{Z^\ell_\e}{\e\tau_\e} \int_{x_0}^x e^{V(\xi)/\e}\, \dd \xi.
\label{eqdef:y_e}	
\end{align}

\begin{figure}
\centering
\labellist
\pinlabel \labelsize $x$ at 1900 480
\pinlabel \labelsize $x_a$ at 570 460
\pinlabel \labelsize $x_b$ at 1300 460
\pinlabel \labelsize $+\frac{1}{2}$ at 830 870
\pinlabel \labelsize $-\frac{1}{2}$ at 830 200
\pinlabel \labelsize $B_a$ at 540 660
\pinlabel \labelsize $B_b$ at 1280 660
\pinlabel \labelsize {\color{red_fig}{$\phi_\varepsilon$}} at 1700 750
\pinlabel \labelsize {\color{blue_fig}{$y_\varepsilon$}} at 1600 1000
\endlabellist
\centering
\includegraphics[scale=.15]{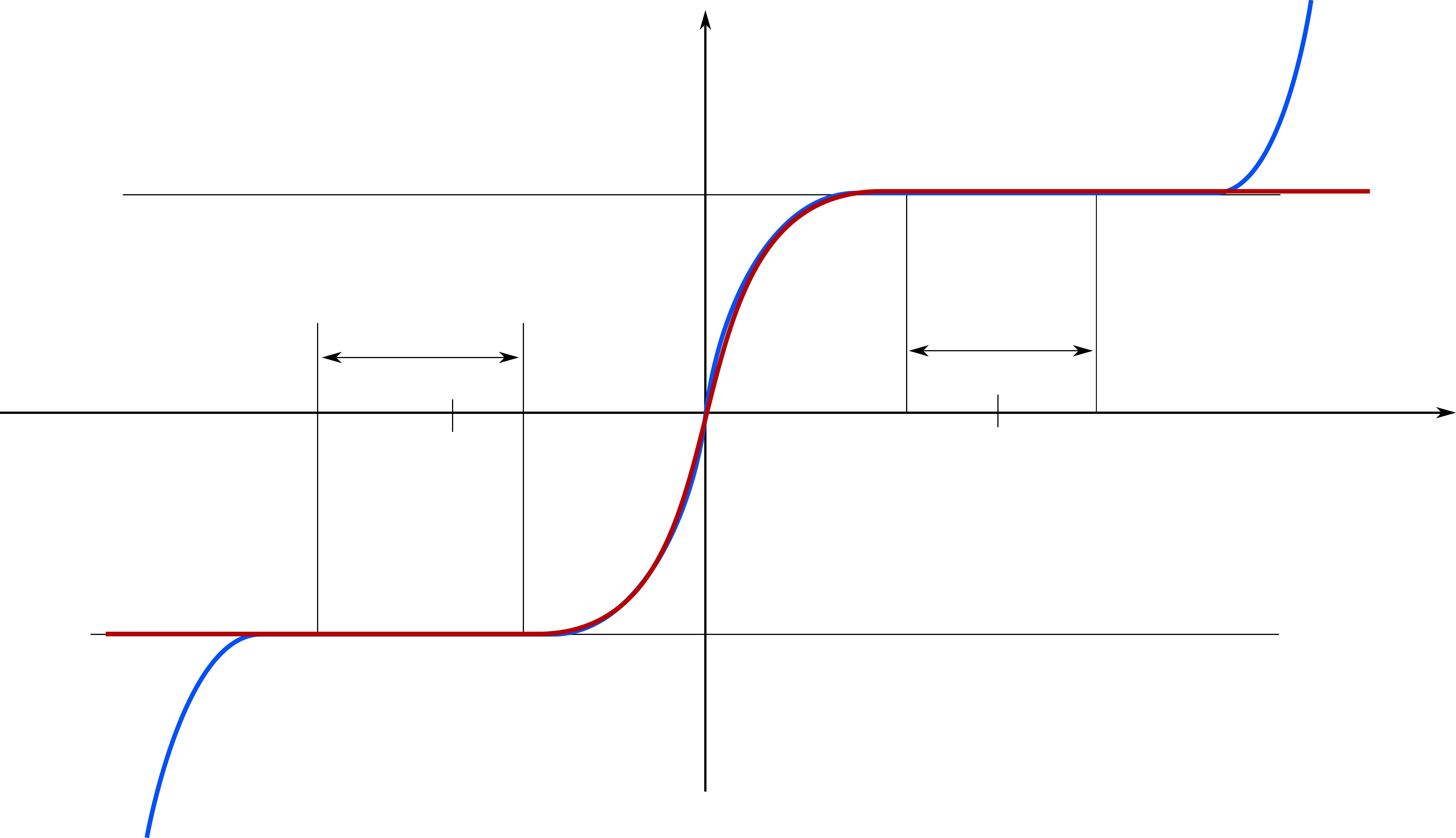}
\caption{Comparison of the functions $\phi_\e$ and $y_\e$. Note how the two functions are very similar in the region between and around the two wells; towards $\pm\infty$, however, $y_\e$ is unbounded, while the range of $\phi_\e$ is bounded.}
\end{figure}
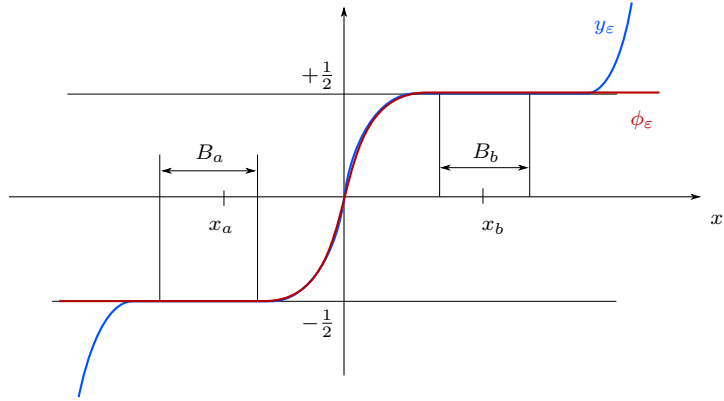
The definition of $\phi_\e$ is a minor modification of~\cite[Lemma~3.6]{EvansTabrizian16} and is nearly the same as the \emph{committor function}, known from potential theory~\cite{BovierDenHollander2016} and Transition-Path Theory~\cite{EVanden-Eijnden04}; see also~\cite{LuVanden-Eijnden14} for a discussion of its use in coarse-graining, which is similar to its function here. The following lemma describes in different ways how  $\phi_\e$ approximates the function $x\mapsto \sign(x)/2$.

\begin{lemma}
\label{l:props-phi_e}
The function $\phi_\e$  satisfies
\begin{enumerate}
	\item $\phi_\e$ is non-decreasing on $\R$;
	\item There exists $C>0$ such that 
	$|\phi_\e|\leq C$
	for sufficiently small~$\e$,  $\lim_{\e\to0} \phi_\e(-\infty) = -1/2$, and $\lim_{\e\to0} \phi_\e(+\infty) = 1/2$;
	\item $\phi_\e$ converges uniformly to $-1/2$ on $B_a$ and to $1/2$ on $B_b$.
	\item \label{l:props-phi_e:mue-conv}
$Z_\e^\ell e^{V/\e} \mu_\e$ converges uniformly on $\R$ to $-\chi_a$.
\end{enumerate}	
\end{lemma}

\begin{proof}
%
The non-negativity of $M_\e$ proves the monotonicity of $\phi_\e$. The bound on $\phi_\e$ and the convergence of the limit values follow from remarking that 
\[
\sup_\R \phi_\e = \phi_\e(+\infty) = \frac{Z^\ell_\e}{\e\tau_\e} \int_{x_0}^\infty e^{V/\e} M_\e.
\]
Since on $\supp M_\e\subset B_a\cup B_0\cup B_b$ the potential $V$ takes its maximum at the saddle $x_0$, and since $M_\e$ is equal to one around the saddle, the  integral converges to $1/2$ by part~\ref{l:basic-props-gamma:conv-prop-2-Z} of Lemma~\ref{l:basic-props-gamma}. The behaviour at $-\infty$ is proved in the same way.

Since the expression $Z_\e^\ell e^{V/\e}/\e\tau_\e $ converges to zero uniformly on $B_a$ and~$B_b$, equation~\eqref{eqdef:phi_e} implies that $\phi_\e$ becomes constant on $B_a$ and $B_b$ and converges uniformly on those sets to its limit values, which are $-1/2$ and $1/2$, respectively. 

Finally, 
\[
Z_\e^\ell e^{V/\e} \mu_\e = Z_\e^\ell\frac{\chi_b}{\int e^{-V/\e}\chi_b} - Z_\e^\ell\frac{\chi_a}{\int e^{-V/\e}\chi_a} =: \alpha_{\e,b} \chi_b - \alpha_{\e,a} \chi_a.
\]

The first term vanishes uniformly since the scalar $\alpha_{\e,b}$ equals $Z_\e^\ell /\int e^{-V/\e}\chi_b \sim e^{V(x_b)/\e} \to 0$. The second term converges to $-\chi_a$, since  
\[
\alpha_{\e,a}^{-1} = \frac1{Z_\e^\ell}\int e^{-V/\e}\chi_a = \int \chi_a \gamma^\ell_\e \longrightarrow \chi_a(x_a) = 1.
\qedhere
\]
\end{proof}

\bigskip

The function $y_\e$ is very similar to $\phi_\e$, but differs in the tails, and will be used as a coordinate transformation in Section~\ref{s:upper-bound}.

\begin{lemma}
\label{l:props-phie-ye}
	\begin{enumerate}
		\item The function $y_\e$ is strictly increasing and bijective.
		\item For any $x<x_0$ such that $V(x) < V(x_0)$, we have $y_\e(x) \to -\frac{1}{2}$ as $\e \to 0$.
		\item For any $x>x_0$ such that $V(x) < V(x_0)$, we have $y_\e(x) \to +\frac{1}{2}$ as $\e \to 0$.
		\item $\phi_\e\circ y_\e^{-1}$ converges uniformly on $\R$ to the truncated identity function $\idhalf$, defined by
\[
\idhalf(x) := \begin{cases}
-1/2 & \text{if }x\leq -1/2\\
x& \text{if } -1/2\leq x\leq 1/2\\
1/2 & \text{if }x\geq 1/2.	
\end{cases}
\]
\end{enumerate}
\end{lemma}

\begin{figure}[h!]
	\labellist
	\pinlabel $V(x)$ at 1730 1250
	\pinlabel $x_0$ at 800 600
	\pinlabel $x_a$ at 400 780
	\pinlabel $x_b$ at 1350 750
	\pinlabel $x$ at 1600 750
	\pinlabel $-\frac{1}{2}$ at 380 -75
	\pinlabel $0$ at 800 -60
	\pinlabel $+\frac{1}{2}$ at 1150 -75
	\pinlabel $y_\e(x)$ at 1750 -50
	\endlabellist
	\centering
	\includegraphics[scale=.1]{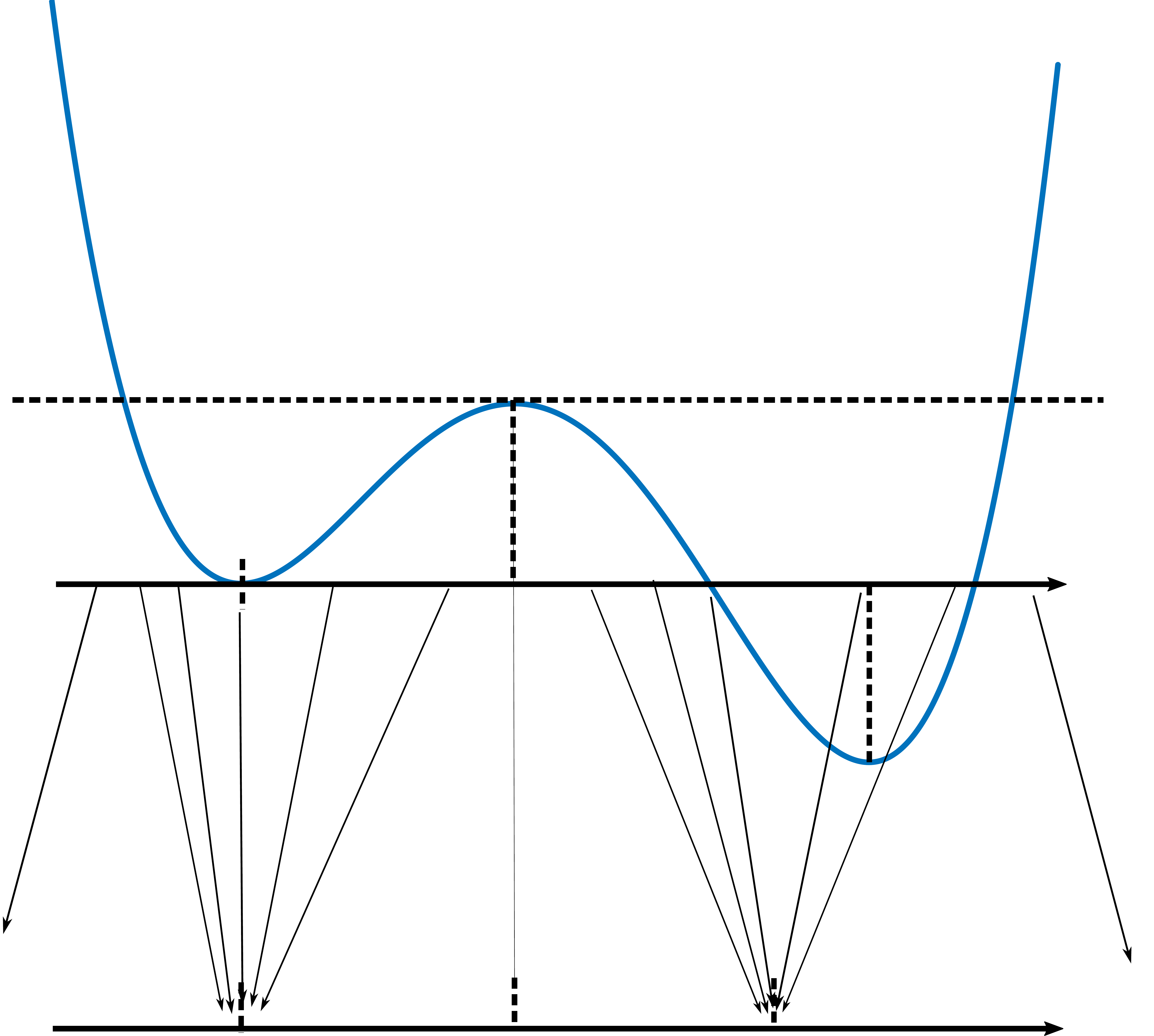}
	
	\kern2ex
	\caption{The structure of the map $y_\e$ of~\eqref{eqdef:y_e}. Points to the left of~$x_0$ with $V(x)<V(x_0)$ are mapped to~$-1/2$, and similarly, points to the right of~$x_0$ are mapped to~$+1/2$. The smaller the value of~$\e$, the sharper is the concentration effect. As~$\e\to 0$, points far to the left of~$x_a$ and far to the right of~$x_b$ are mapped to~$\mp\infty$, respectively.}
	\label{fig:coordinate-transformation}
\end{figure}
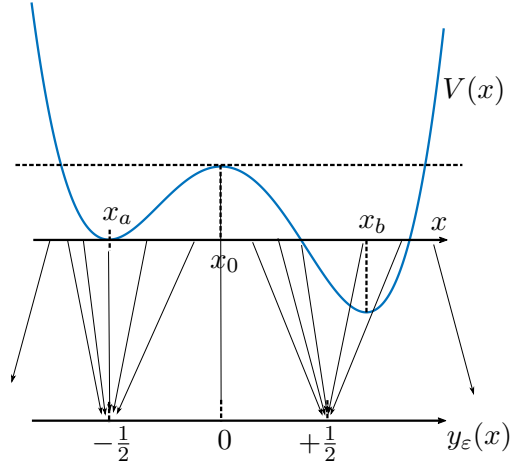

\begin{proof}
Since $y_\e'(x) > 0$ for any $x \in \mathbb{R}$ and $y_\e(x) \to \pm \infty$ as $x \to \pm\infty$, the map $y_\e$ is strictly increasing and bijective. For $x < x_0$ satisfying $V(x) < V(x_0)$, we obtain 
\begin{align*}
y_\e(x) &= \frac{1}{\e \tau_\e} \cdot Z^\ell_\e \cdot  \int_{x_0}^x e^{V(z)/\e}\, dz	\\
&=[1+o(1)] \frac{1}{\e \tau_\e} \cdot e^{-V(x_a)/\e} \sqrt{\frac{2\pi\e}{V''(x_a)}} \cdot \frac{1}{2}e^{V(x_0)/\e} \sqrt{\frac{2\pi\e}{|V''(x_0)|}} (-1)\\
&\longrightarrow -\frac{1}{2},
\end{align*}
by using~\eqref{char:asymp-Ze} and applying Lemma~\ref{lemma:watson}\ref{char:Laplace:boundary} to the integral.  The argument for the case~$x>x_0$ is similar.

To show that $\phi_\e\circ y_\e^{-1}$ converges uniformly on $\R$ to $\idhalf$, first note that 
\[
\frac \dd {\dd y} \phi_\e\bigl(y_\e^{-1}(y)\bigr)= \frac{\phi'_\e\bigl(y_\e^{-1}(y)\bigr)}{y_\e'\bigl(y_\e^{-1}(y)\bigr)} = M_\e\bigl(y_\e^{-1}(y)\bigr) \qquad\text{for any }y\in \R.
\]
The function $M_\e\circ y_\e^{-1}$ converges in $L^1(\R)$ to $\bONE_{[-1/2,1/2]}$; this can be recognized from the fact that $y_\e^{-1}(y)$ converges to $0$ for any $-1/2<y<1/2$, to $+\infty$ for $y>1/2$, and to $-\infty$ for $y< -1/2$. The uniform convergence of $\phi_\e\circ y_\e^{-1}$ then follows by integration.
\end{proof}

\subsection{Compactness and lower bound}

Having defined the auxiliary function $\phi_\e$ we can state and prove the main compactness theorem, which includes a lower bound  on~$\cI_\e$.

\begin{theorem}[Compactness and lower bound]
\label{th:compactness-lower-bound-untransformed}
Let $V$ satisfy Assumption~\ref{ass:V}.
Let $(\rho_\e,j_\e)\in \CE(0,T)$	satisfy 
\begin{equation}
\label{ass:bdd-Ie-Ee}
\sup_\e \cI_\e(\rho_\e,j_\e) + \e E_\e(\rho_\e(0)) \leq C< \infty, 
\end{equation}
and assume that $\rho_\e(0)$ satisfies the narrow convergence
\begin{equation}
\label{ass:initial-data-rho_e}	
\rho_\e(0) \longweakto \rho_0^\circ(\dd x) := z^\circ\delta_{x_a}(\dd x) + (1-z^\circ)\delta_{x_b}(\dd x)\qquad \text{as } \e\to0.
\end{equation}
Then there exists a $(\rho_0,j_0)\in \CE(0,T)$ and a subsequence along which
\begin{enumerate}
\item \label{th:compactness-lower-bound-untransformed:char-rhozero}
$\rho_\e\longweakto\rho_0$ narrowly in $\cM([0,T]\times\R)$, where $\rho_0\in \cM([0,T]\times \R)$ has the structure
\begin{equation}
\label{eq:structure-rhozero}
\rho_0(\dd t\dd x) = \rho_0(t,\dd x)\dd t := z(t)\delta_{x_a}(\dd x) \dd t + (1-z(t))\delta_{x_b}(\dd x)\dd t,
\end{equation}
and $z:[0,T]\to[0,1]$ is absolutely continuous.
\item \label{th:compactness-lower-bound-untransformed:conv-j} 
$j_\e$ converges in duality with $C_c^{1,0}([0,T)\times \R)$ to 
\[
j_0(\dd t \dd x) := j(t)\bONE_{[x_a,x_b]}(x) \dd x\dd t, 
\]
where $j(t) = -z'(t)$ for almost all $t\in[0,T]$.
\item $\liminf_{\e\to0} \cI_\e(\rho_\e,j_\e)\geq \cI_0(\rho_0,j_0)$.
\end{enumerate}
\end{theorem}

\begin{remark}
Note that the two assumptions on the initial data, the convergence~\eqref{ass:initial-data-rho_e} and the boundedness $E_\e(\rho_\e(0))\leq C/\e$ of~\eqref{ass:bdd-Ie-Ee}, are closely related, but independent: it is possible to satisfy one but not the other. 
\end{remark}

\begin{proof}
%
%
%
Recall from the discussion in Section~\ref{s:elements-of-the-proof} that by the assumption~\eqref{ass:bdd-Ie-Ee} on the initial data we have the `fundamental estimate'
\begin{equation}
\label{ineq:fundamental-estimate}
\frac{\e\tau_\e}2 \int_0^T \int_\R 
\Bigl|\partial_x \sqrt{u_\e(t,x)} \Bigr|^2\gamma_\e(\dd x)\, \dd t
+ \sup_{t\in[0,T]}E_\e(\rho_\e(t)) \leq \frac C\e.
\end{equation}
Here $u_\e$ is the density of $\rho_\e$ with respect to the invariant measure $\gamma_\e$.

\textit{Step 1: 
Concentration for the case of the outer
half-lines.
} Set $O_\ell := (-\infty, x_{c\ell}]$. Recall that $V''\geq \alpha>0$ on $O_\ell$; by Lemma~\ref{l:LSI} we therefore have
\[
\frac\alpha\e \cE(\mu|\gamma_\e,O_\ell) \leq \cR(\mu|\gamma_\e,O_\ell)
\qquad\text{for all }\mu\in \cM_{\geq0}(\R)\text{ and }\e>0.
\]
Then \begin{align*}
\int_0^T \cE(\rho_\e(t)|\gamma_\e,O_\ell)\,\dd t 
&\leq \frac\e\alpha \int_0^T \cR(\rho_\e(t)	|\gamma_\e,O_\ell)\,\dd t \\
&= \frac\e\alpha \int_0^T \frac12 \int_{O_\ell}
  \left|\partial_x \sqrt{\frac{\dd \rho_\e(t)}{\dd \gamma_\e}}\right|^2 \gamma_\e(\dd x)\dd t\\
&\leq \frac\e\alpha \cdot \frac C{\e^2\tau_\e}\qquad \text{ by~\eqref{ineq:fundamental-estimate} }\\
&= \frac C{\alpha \e\tau_\e} \longrightarrow 0 \qquad\text{as }\e\to0.
\end{align*}
Therefore, if $A\subset O_\ell$ with $\dist(A,\{x_a\})>0$, then by Lemma~\ref{l:concentration-LSI},
\[
\int_0^T \rho_\e(t,A)\, \dd t 
\leq \Bigl(\log \frac{\gamma_\e(A)}{\gamma_\e(O_\ell)}\Bigr)^{-1}
  \int_0^T \Bigl[\cE(\rho_\e(t)|\gamma_\e,O_\ell)+1\Bigr]\,\dd t 
\stackrel{\e\to0}\longrightarrow 0.
\]
It follows that $\rho_\e\bONE_{O_\ell}$ concentrates onto $[0,T]\times \{x_a\}$. By a similar argument $\rho_\e\bONE_{[x_{cr},\infty)}$ concentrates onto $[0,T]\times \{x_b\}$. This also implies that $\rho_\e$ is tight on $[0,T]\times \R$.

\medskip

\textit{Step 2: 
Concentration for the case of the whole domain~$\mathbb{R}$.
}
We have proved concentration of $\rho_\e\bONE_{(-\infty,x_{c\ell}]}$  onto $[0,T]\times \{x_a\}$ and of $\rho_\e\bONE_{[x_{cr}, \infty)}$  onto $[0,T]\times \{x_b\}$. What remains is to bridge the gap between $x_{c\ell}$ and $x_{cr}$. 

We write $u^\ell_\e$ for the density of $\rho_\e$ with respect to the left-normalized invariant measure $\gamma_\e^\ell$, i.e.\ $\hat u_\e^\ell = u_\e Z_\e /Z_\e^\ell$. 
We then estimate
\begin{align*}
\int_0^T &\int_\R \Bigl|\partial_x \sqrt{u_\e^\ell}\Bigr|^2(t,x) e^{(V(x_0)-V(x))/\e}\, \dd x \dd t\\
&= Z_\e^\ell \int_0^T \int_\R \Bigl|\partial_x \sqrt{u_\e^\ell}\Bigr|^2(t,x) e^{V(x_0)/\e}\gamma_\e^\ell( \dd x) \dd t\\
&= C\tau_\e Z_\e^\ell \int_0^T \int_\R \Bigl|\partial_x \sqrt{u_\e}\Bigr|^2(t,x) \gamma_\e( \dd x) \dd t\stackrel{\eqref{ineq:fundamental-estimate}} \leq  C\e^{-3/2}.	
\end{align*}
Since $V\leq V(x_0)$ on $[x_a,x_{b+}]$ it follows that 
\[
\int_0^T \int_{x_a}^{x_{b+}} \Bigl|\partial_x \sqrt{u_\e^\ell}\Bigr|^2(t,x) \, \dd x \dd t\leq C\e^{-3/2}.
\]

Applying the generlized Poincar\'e inequality of Lemma~\ref{l:gen-Poincare} to $f(t,x) = \sqrt{u^\ell_\e}$ on $[x_a,x_{b+}]$ we find 
\begin{align*}
\|u_\e^\ell\|_{L^1(0,T;L^\infty(x_a,x_{b+}))}
&= \int_0^T \|u_\e^\ell(t)\|	_{L^\infty(x_a,x_{b+})}\, \dd t\\
&\leq C \biggl[\e^{-3/2} + \int_0^T  \gamma_\e^\ell([x_a,x_{b+}])^{-1} \|u_\e^\ell(t)\|	_{L^1(x_a,x_{b+};\gamma_\e^\ell)}\, \dd t\biggr]\\
&= C \biggl[\e^{-3/2} + \gamma_\e^\ell([x_a,x_{b+}])^{-1} \int_0^T  \rho_\e(t;[x_a,x_{b+}])\, \dd t \biggr]\\
&\leq C\e^{-3/2} \qquad\text{since $\gamma_\e^\ell([x_a,x_{b+}])\to\infty$ by Lemma~\ref{l:basic-props-gamma}.}
\end{align*}

To prove concentration, take an interval $A$ such that $[x_{c\ell},x_{cr}]\subset A\subset \{V\geq \delta\}$ for some $\delta>0$.  Then
\begin{align*}
\int_0^T \rho_\e(t,A)\, \dd t 
&= \int_0^T \int_A u_\e^\ell(t,x) \gamma_\e^\ell(\dd x) \, \dd t \\
&\leq  \gamma_\e^\ell(A)\|u_\e^\ell\|_{L^1(0,T;L^\infty(x_a,x_{b+}))}
 \\
&\leq  C\e^{-3/2} \cdot \frac{e^{-\delta/\e}}{Z^\ell_\e} 
\stackrel{\e\to0}\longrightarrow0.
\end{align*}
Therefore $\rho_\e$ does not charge the region $[0,T]\times A$ in the limit. 

\smallskip

Concluding, $\rho_\e$ concentrates onto $[0,T]\times \{x_a,x_b\}$ as $\e\to0$. It follows that the  limit $\rho_0$ has support contained in $[0,T]\times \{x_a,x_b\}$, and for almost every $t\in [0,T]$, $\rho_0(t,\cdot)$ has mass one on $\R$. This establishes the structure~\eqref{eq:structure-rhozero}, except for the continuity of $z$;  at this stage we only know that $z\in L^\infty(0,T)$ with $0\leq z\leq 1$, and the absolute continuity of $z$ will follow in Step 4 below.

\medskip

\begin{remark}
After completing the proof of compactness outlined in the previous two steps, Andr\'e Schlichting pointed out that by using the Muckenhoupt criterion it is possible to replace the  assumption of convex wells by two monotonicity assumptions, one for each well; see Theorem~3.19 in~\cite{Schlichting12TH} for an example.
\end{remark}

\textit{Step 3: Lower bound on $\mathcal I_\e$.}
From Definition~\ref{def:Ie} and the bound~\eqref{ass:bdd-Ie-Ee} we have for any $b\in C_c^{0,1}(\tmsp)$ 	the estimate
\begin{equation}
\label{ineq:est:lb}
C \geq \cI_\e(\rho_\e,j_\e)
\geq \int_0^T \int_\R \Bigl[ j_\e b - \e\tau_\e \rho_\e 
  \Bigl( \partial_x b - \frac1\e b V' + \frac12 b^2 \Bigr)\Bigr]\, \dd x \dd t.
\end{equation}
Fix $\psi\in C^1([0,T])$ with $\inf \psi > -1$ and $\psi(T)=0$. Define $F_\e:[0,T]\times \R\to\R$ by 
\[
F_\e(t,x) := \log \Bigl(1+\psi(t)\underbrace{(\tfrac12 - \phi_\e(x))}_{=:\wt \phi_\e(x)}\,\Bigr), 
\qquad\text{with $\phi_\e$ given by~\eqref{eqdef:phi_e}}.
\]
\begin{lemma}
$F_\e$ and $\wt\phi_\e$ have  the following properties:
\begin{enumerate}
\item $F_\e\in C^1_b(\tmsp)$ and $\partial_x F_\e\in C^1_c(\tmsp)$;
\item $F_\e(T,x) = 0$ for all $x\in\R$;
\item $\sup_{\e,t,x} |F_\e(t,x)| \leq \max\{	\log(1+\sup\psi), -\log(1+\inf\psi)\}$;
\item $\wt \phi_\e$ converges uniformly on $[0,T]\times B_a$ to $1$ and on $[0,T]\times B_b$ to zero;
\item $F_\e$ converges uniformly on $[0,T]\times B_a$ to $\log(1+\psi(t))$ and on $[0,T]\times B_b$ to zero.
\end{enumerate}
\end{lemma}
\noindent 
These follow directly from Lemma~\ref{l:props-phie-ye}.

\smallskip
We now set  $b_\e(t,x) = 2\partial_x  F_\e (t,x)= 2\psi(t){\wt \phi_\e}'(x)/(1+\psi(t)\wt\phi_\e(x))$ and find that the expression in brackets in~\eqref{ineq:est:lb} equals
\begin{align*}
\partial_x b_\e - \frac1\e b_\e V' + \frac12 b_\e^2 
&= 	 \frac{2\psi}{1+\psi\wt\phi_\e} \Bigl[ {\wt\phi_\e''} - \frac1\e {\wt\phi_\e'}V'\Bigr] \\
&\leftstackrel{\eqref{eqdef:phi_e}} = \frac{2\psi}{1+\psi\wt\phi_\e}  \frac{Z_\e^\ell}{\e\tau_\e} e^{V/\e}\mu_\e.
\end{align*}
By Lemma~\ref{l:props-phi_e} and the concentration of $\rho_\e$ we therefore find that 
\begin{align}
\lim_{\e\to0} \int_0^T \int_\R \e\tau_\e \rho_\e \Bigl( \partial_x b_\e &- \frac1\e b_\e V' + \frac12 b_\e^2 \Bigr)\notag\\
&= \lim_{\e\to0} \int_0^T \int_\R \frac{2\psi}{1+\psi\wt\phi_\e} \rho_\e Z_\e^\ell e^{V/\e}\mu_\e\notag\\
&= 2\int_0^T \int_\R  \frac{\psi(t)}{1+\psi(t)}\, \rho_0(t,\dd x) \chi_a(x)\,\dd t \notag\\
&= 2\int_0^T \frac{\psi(t)}{1+\psi(t)}\, z(t)\, \dd t. 
\label{eq:limit-of-second-term}
\end{align}

\smallskip
We now turn to the first term in~\eqref{ineq:est:lb}. Applying the Definition~\ref{def:continuity-equation} of $\CE$, and the assumption~\eqref{ass:initial-data-rho_e} on the convergence of the initial data, we find
\begin{align}
\int_0^T \int_\R j_\e b_\e
&= -2\int_0^T \int_\R \rho_\e \partial_t F_\e -2 \int_\R \rho_\e(0,\dd x)F_\e(0,x)\notag\\
&= -2\int_0^T \int_\R \rho_\e \frac{\psi'\wt\phi_\e}{1+\psi\wt\phi_\e}  -2 \int_\R \rho_\e(0,\dd x)F_\e(0,x)\notag\\
&\stackrel{\e\to0}\longrightarrow \;
  -2\int_0^T z(t)  \frac{\psi'(t)}{1+\psi(t)} \, dt  -2 z^\circ \log(1+\psi(0)).
  \label{eq:limit-of-first-term}
\end{align}
Writing $f(t) := -\log(1+\psi(t))$ we have $f(T) =0$; combining~\eqref{eq:limit-of-second-term} and~\eqref{eq:limit-of-first-term}, and observing that $\psi/(1+\psi) = e^f-1$, we find 
\[
\liminf_{\e\to0} \cI_\e(\rho_\e,j_\e) \geq \cJ_0(z),
\]
with 
\[
\cJ_0(z) := 2\sup \biggl\{\int_0^T z(t)\Bigl[f'(t) - e^{f(t)}+ 1\Bigr]\, \dd t + z^\circ f(0) \;: \;f\in C^1_b([0,T]), \ f(T)=0 \biggr\}.
\]

\begin{lemma}
\label{l:var-formulation-I0}
Let $z\in L^\infty(0,T)$ with $z\geq0$, and let $z^\circ\geq 0$. Then $\cJ_0(z) = \cK_0(z)$, 
where
\[
\cK_0(z) := \begin{cases}
\ds 2\int_0^T S\bigl(-\overline z'(t)|\overline z(t)\bigr)\, \dd t &
	\text{if $z=\overline z$ a.e. with $\overline z$ non-increasing }\\[-2\jot]
	& \qquad\text{and absolutely continuous, and $\overline z(0)=z^\circ$}\\[2\jot]
+\infty & \text{otherwise.}
\end{cases}
\]
If $\cK_0(z) = 0$, then $z(t) = z^\circ e^{-t}$ for almost all $0\leq t\leq T$.
\end{lemma}

We prove this lemma below, and first finish the proof of Theorem~\ref{th:compactness-lower-bound-untransformed}. Note that since $\cJ_0(z) = \cK_0(z)<\infty$, the function $z$ has an absolutely continuous representative and $z(0) = z^\circ$; this concludes the proof of  part~\ref{th:compactness-lower-bound-untransformed:char-rhozero} of the theorem.

\medskip

\textit{Step 4 of the proof of Theorem~\ref{th:compactness-lower-bound-untransformed}: Convergence of $j_\e$.}
Choose any $\varphi\in C_c^{1,0}(\tmsp)$ with $\varphi=0$ at $t=T$, and set $\Phi(t,x) := \int_0^x \varphi(t,\xi)\, \dd \xi$; note that  $\Phi\in C^1_b(\tmsp)$ and $\partial_x\Phi\in C_c(\tmsp)$. We calculate
\begin{align*}
\int_0^T\int_\R j_\e \varphi 
&= \int_0^T \int_\R j_\e(t,\dd x) \partial_x \Phi(t,x)\, \dd x\dd t\\
&\leftstackrel{\eqref{eq:weak-form-CE}} = - \int_0^T \int_\R \rho_\e(t,\dd x) \partial_t \Phi(t,x)\, \dd x \dd t
  - \int_\R \rho_\e(0,\dd x)\Phi(0,x)\\
&\stackrel{\e\to0}\longrightarrow
  -\int_0^T \int_\R \rho_0(t,\dd x) \partial_t \Phi(t,x)\, \dd x \dd t
  - \int_\R \rho_0^\circ(\dd x)\Phi(0,x)\\
&\leftstackrel{(\ref{ass:initial-data-rho_e},\ref{eq:structure-rhozero})} = -\int_0^T \bigl[ z(t) \partial_t \Phi(t,x_a) + (1-z(t))\partial_t \Phi(t,x_b)\bigr]\, \dd t\\
& \qquad {}
 - z(0)\Phi(0,x_a) - (1-z(0))\Phi(0,x_b)\\
&= \int_0^T z'(t) \bigl(\Phi(t,x_a)-\Phi(t,x_b)\bigr)\, \dd t\\
&= \int_0^T \int_\R (-z'(t)) \varphi(t,x) \bONE_{[x_a,x_b]}(x) \,\dd x \dd t.
\end{align*}
This proves the convergence of part~\ref{th:compactness-lower-bound-untransformed:conv-j}. Finally, with this definition of the limit $j_0$ of $j_\e$, we have 
\[
\cI_0(\rho_0,j_0) = \cJ_0(z) = \cK_0(z),
\]
and this concludes the proof of Theorem~\ref{th:compactness-lower-bound-untransformed}.
\end{proof}

\begin{proof}[Proof of Lemma~\ref{l:var-formulation-I0}]
A closely related statement and its proof are discussed in~\cite[Sec.~3]{PattersonRenger19}; for completeness we give a standalone proof.

\textit{Step 1: If $\cJ_0(z)<\infty$, then $z$ is non-increasing on $[0,T]$.}
Fix $\varphi\in C_c^\infty((0,T))$ with $\varphi\geq0$. Applying the definition of $\cJ_0(z)$ to $f = -\lambda \varphi$ we find 
\[
-\lambda \int_0^T z\varphi' \leq \cJ_0(z) \qquad \text{for all }\lambda>0,
\]
which implies $\int z\varphi'\geq 0$. Since $\varphi\in C_c^\infty((0,T))$ is arbitrary, it follows that the equivalence class $z\in L^\infty$ has a non-increasing representative, and  from now on we write $z$ for this  non-increasing representative.  We also find that $-z'$ is a positive measure on $(0,T)$. By the monotonicity of $z$, the limits of $z$ at $t=0,T$ exist, and if necessary we redefine $z$ to be continuous at $t=0,T$. By construction, $z$ now is non-increasing on $[0,T]$ and $-z'$ is a positive measure on $[0,T]$ without atoms at $t=0,T$.

\medskip
\textit{Step 2: Reformulation and matching initial data.}
Since $-z'$ is a finite measure and $z$ is continuous at $t=0,T$, we can rewrite
\begin{multline*}
	\cJ_0(z) = 2\sup \biggl\{\int_0^T \Bigl[-z'(t)f(t) - z(t)\bigl(e^{f(t)}- 1\bigr)\Bigr]\, \dd t + (z^\circ-z(0)) f(0) \;:\\
	 \;f\in C^1_b([0,T]), \ f(T)=0 \biggr\}.
\end{multline*}
By choosing functions $f$ with $f(0)=\lambda\in \R$ and $\supp f$ a vanishingly small interval close to $t=0$ we find $\cJ_0(z) \geq 2\lambda(z^\circ-z(0))$, and by taking limits $\lambda\to\pm\infty$ it follows that  $z(0) = z^\circ$. 

\medskip
\textit{Step 3: Primal form.}
Still under the assumption that $\cJ_0(z)<\infty$, we recognize $e^f-1$ as the dual $\eta^*(f)$ of the function $\eta(a) := S(a|1)$ (which is equal to $a\log a -a + 1$ for $a>0$). We then use the well-known duality characterization of convex functions of measures (see e.g.~\cite[Lemma~9.4.4]{AmbrosioGigliSavare2008}) to find, writing $\mu(\dd t) := z(t)\dd t$, 
\[
\cJ_0(z) = 
\begin{cases}
\ds\int_0^T \eta\Bigl( \frac{\dd (-z')}{\dd \mu}(t)\Bigr) \mu(\dd t) 
&\text{if }{-z'} \ll \mu\\
+\infty & \text{otherwise},
\end{cases}
\]
and this functional coincides with $\cK_0$ (see e.g.~\cite[Lemma~2.3]{PeletierRossiSavareTse20TR}).
The reverse statement, assuming $\cK_0(z)<\infty$ and showing that $\cJ_0(z) < \infty$, follows directly by Young's inequality for the pair $(\eta,\eta^*)$.

\medskip
\textit{Step 4: Absolute continuity.}
Finally, if $\cJ_0(z)=\cK_0(z)<\infty$, then the superlinearity of $\eta$ implies that $z'\in L^1(0,T)$, and therefore $z$ is absolutely continuous. 

\medskip
\textit{Step 5: Characterization of minimizers.}
If $\cK_0(z)=0$, then $\overline z'(t) =-\overline z(t)$ for almost all $t$, implying that $\overline z(t) = z^\circ e^{-t}$.
\end{proof}

\section{Recovery sequence}
\label{s:upper-bound}

In this section we state and prove 
Theorem~\ref{t:recovery-sequence-untransformed}, which establishes the existence of a recovery sequence for the $\Gamma$-convergence of Theorem~\ref{thm:intro-main-result}.

\subsection{Spatial transformation}

We start by transforming the system by a nonlinear mapping in space, given by the function $y_\e$ defined in Section~\ref{ss:aux-fcns}; this function maps $\R$ with variable~$x$ to $\R$ with variable~$y$, and is inspired by a similar choice in~\cite{ArnrichMielkePeletierSavareVeneroni2012}. This mapping  desingularizes the system. 

We define the transformed versions $\hat \rho_\e$ and $\hat\gamma^\ell_\e$  of $\rho_\e$ and $\gamma^\ell_\e$ by pushing them forward under $y_
\e$,
\begin{subequations}
\label{def:transformation}
\begin{equation}
\label{def:transformation-rho-hat}
\hat\rho_\e := (y_\e)_\# \rho_\e \qquad \text{and}\qquad
\hat \gamma^\ell_\e := (y_\e)_\# \gamma^\ell_\e,
\end{equation}
which implies that the transformed density $\hat u^\ell_\e$ is  given by
\begin{equation}
\label{def:transformation-u-hat}
\hat{u}_\e^\ell(t,y_\e(x)) := u_\e^\ell(t,x).
\end{equation}
We transform $j_\e$ in such a way that the continuity equation is conserved, which leads to the choice
\begin{equation}
\label{def:transformation-j-hat}
\hat \jmath_\e := (y_\e)_\# \bigl(y_\e'j_\e\bigr),
\end{equation}
which has an equivalent formulation in the case of Lebesgue-absolutely-continuous fluxes,
\begin{equation}
\label{def:transformation-j-hat-ac}
\hat \jmath_\e (t,y_\e(x)):= j_\e(t,x).
\end{equation}
\end{subequations}

Indeed, if $(\rho,j)$ satisfies the continuity equation~\eqref{eq:upscaled-CE}, then the transformed pair $(\hat \rho_\e,\hat \jmath_\e)$ satisfies the corresponding continuity equation in the variables $(t,y)$,
\[
\partial _t \hat \rho_\e + \partial_y \hat \jmath_\e = 0
\qquad\text{in}\qquad \R^+\times \R,
\]
which is defined again as in Definition~\ref{def:continuity-equation}, and one can check that $(\rho,j)\in \CE(0,T) \iff (\hat \rho,\hat\jmath)\in \CE(0,T)$. Since $y_\e$ is a diffeomorphism, there is a one-to-one relationship between $(\rho,j)$ and $(\hat \rho,\hat\jmath)$.

\medskip
In terms of $\hat \jmath_\e$ and the density $\hat u_\e^\ell$ the rate function formally takes the simpler form
\begin{equation*}
\mathcal{I}_\e(\rho,j) = \frac{1}{2}\int_0^T\int_\mathbb{R} \frac{1}{\hat{u}_\e^\ell(t,y)}\big|\hat{\jmath}_\e(t,y) + \partial_y \hat{u}_\e^\ell(t,y)\big|^2\,\dd y\dd t.
\end{equation*}
Note how the parameters $\e$ and $\tau_\e$ are absorbed into the density~$\hat{u}_\e^\ell$ and the derivative with respect to the new coordinate $y$. The coordinate transformation~$y_\e$ is the almost the same as in~\cite{ArnrichMielkePeletierSavareVeneroni2012}; the only difference is that we use the left-normalized stationary measure, whereas in the symmetric case one can use the stationary measure normalized in the usual manner. 

This simpler, transformed form is the basis for the construction of the recovery sequence. To make this precise we first define the rescaled versions of $\cI_\e$ and $\cI_0$. 

\begin{definition}[Rescaled functionals]
\label{def:rescaled-functionals}
For given $\rho$ and $j$, define $\hat \rho$ and $\hat \jmath$ as in~\eqref{def:transformation-rho-hat} and~\eqref{def:transformation-j-hat}. We define $\hat E_\e$, $\widehat\cI_\e$, and $\widehat\cI_0$ to be the rescaled versions of $E_\e$, $\cI_\e$, and $\cI_0$,
\begin{align*}
&\hat E_\e: \cP(\R)\to[0,\infty], && \hat E_\e(\hat \rho) := E_\e(\rho),\\
&\widehat\cI_\e: \CE(0,T) \to [0,\infty], & &\widehat \cI_\e(\hat\rho,\hat \jmath) := \cI_\e(\rho,j),\\
&\widehat\cI_0: \CE(0,T) \to [0,\infty], & &\widehat \cI_0(\hat\rho,\hat \jmath) := \cI_0(\rho,j).
\end{align*}
\end{definition}

\noindent
The following lemma is a direct consequence of the definition~\eqref{eqdef:Ie}, the transformation~\eqref{def:transformation}, and part~\ref{l:reg-rho-j:rho} of Lemma~\ref{l:reg-rho-j}.
\begin{lemma}[Dual formulation of $\widehat\cI_\e$]\label{l:pre-limit-RF}
	We have
	\begin{equation}\label{def:GF-NGF:pre-limitRF}
	\widehat{\mathcal{I}}_\e\left(\hat{\rho},\hat{\jmath}\right) =  \sup_{b\in C_c^\infty(\tmsp)}\int_\tmsp \left[\hat \jmath \,b -  \hat{u}_\e^\ell\left(\partial_y b + \frac{1}{2}b^2\right)\right]\,\dd y\dd t,
	\end{equation}
	provided~$\hat{\rho}(t,\cdot)$ is absolutely continuous with respect to~$\hat{\gamma}_\e^\ell$ with density $\hat{u}_\e^\ell(t,\cdot)$; otherwise we set~$\widehat{\mathcal{I}}_\e\left(\hat{\rho},\hat{\jmath}\right)=+\infty$.\qed
\end{lemma}


\medskip
While the left-normalized stationary measure $\gamma^\ell_\e$ in the original variables concentrates onto the set $\{x: V(x)\leq 0\}=\{x_a\}\cup [x_{b-1},x_{b+}]$, under this  transformation the interval $[x_{b-1},x_{b+}]$ collapses onto a point (see also Figure~\ref{fig:coordinate-transformation}):
\begin{lemma}[The measures $\hat{\gamma}_\e^\ell$ concentrate onto $\{\pm1/2\}$]\label{l:transformed-stationary-measure}
	Let a measurable set $A\subset \R$ have positive distance to $\pm 1/2$. 
Then
	\begin{equation*}
	\lim_{\e\to 0} \hat{\gamma}_\e^\ell(A) = 0.
	\end{equation*}
\end{lemma}

\begin{proof}
Fix $0< \delta< V(x_0)$. Since $A$ has positive distance to $\{\pm1/2\}$, by Lemma~\ref{l:props-phie-ye} we have for sufficiently small $\e$ that $V\geq \delta$ on $y_\e^{-1}(A)$. Therefore
\[
\hat \gamma_\e^\ell (A) = \gamma_\e^\ell(y^{-1}_\e(A))
\leq \gamma_\e^\ell(\{x: V(x) \geq \delta\}).
\]
By Lemma~\ref{l:basic-props-gamma}, the right-hand side vanishes in the limit~$\e\to 0$.
\end{proof}

\subsection{
Statement and proof for the transformed system
}

\begin{theorem}[Upper bound in transformed coordinates]\label{thm:upper-bound}
	For any~$(\hat{\rho}_0,\hat{\jmath}_0)\in\CE(0,T)$ such that $\widehat{\mathcal{I}}_0(\hat{\rho}_0,\hat{\jmath}_0)<\infty$, there exist~$(\hat{\rho}_\e,\hat{\jmath}_\e)\in\CE(0,T)$ such that 
\begin{align}
\label{eq:th:recovery:conv-and-bound}
	&(\hat{\rho}_\e,\hat{\jmath}_\e)\xrightarrow{\mathrm{CE}}(\hat{\rho}_0,\hat{\jmath}_0)
	\qquad \text{and}\qquad \sup_{\e>0} \e \hat E_\e(\hat \rho_\e(0)) < \infty,
\end{align}
and that 
\begin{equation}
\label{ineq:limsup}
\limsup_{\e\to 0}\widehat{\mathcal{I}}_\e(\hat{\rho}_\e,\hat{\jmath}_\e)
	\leq   \widehat{\mathcal{I}}_0(\hat{\rho}_0,\hat{\jmath}_0).
\end{equation}

\end{theorem}

\begin{proof}	
Recall that $\tmsp:= [0,T]\times\mathbb{R}$ and set $\tmspz:=[0,T]\times[-1/2,+1/2]$.
If $\widehat \cI_0(\hat\rho_0,\hat\jmath_0)$  is finite, then by combining Definitions~\ref{def:rescaled-functionals} and~\eqref{eqdef:I0} we find that  the pair~$(\hat{\rho}_0,\hat{\jmath}_0)$ is given by
\begin{align}\label{eq:limsup:rho_0}
	\hat{\rho}_0(t,dy) &= \hat{z}_0(t)\delta_{-1/2}(dy) + (1-\hat{z}_0(t))\delta_{+1/2}(dy),\\
	\label{eq:limsup:j_0}
	\hat{\jmath}_0(t,dy) &= \hat{\jmath}_0(t) \bONE_{(-1/2,+1/2)}(y)\, \dd y,
\end{align}
where~$t\mapsto \hat{z}_0(t)$ is  absolutely continuous and~$\hat{\jmath}_0$ satisfies $\hat{\jmath}_0(t)=-\partial_t\hat{z}_0(t) \geq 0$. For the later construction of $(\hat\rho_\e,\hat\jmath_\e)$ we will want to assume that~$\hat z_0$ satisfies the following regularity assumption.

\begin{assumption}\label{assump:z0-is-ct-and-pos}
	The density~$\hat{z}_0:[0,T]\to[0,1]$ satisfies
\begin{equation}
\label{eq:assump:z0-is-ct-and-pos}
\hat{z}_0\in C^2([0,T])\quad \text{and}\quad  \inf_{t\in[0,T]}\hat z_0(t), |\partial_t \hat{z}_0(t)|>0.
\end{equation}
Note that this implies that $\hat\jmath_0$ is bounded away from zero and of class $C^1$.
\end{assumption}

Indeed, we can assume that $\hat z_0$ has this regularity since this set is energy-dense:

\begin{lemma}[Energy-dense approximations]
\label{lemma:upperbound:energy-dense-set-of-densities}
	If~$\widehat{\mathcal{I}}_0(\hat{\rho}_0,\hat{\jmath}_0)$ is finite, then there are densities~$\hat{z}_0^\delta$ satisfying Assumption~\ref{assump:z0-is-ct-and-pos} such that the pair~$(\hat{\rho}_0^\delta,\hat{\jmath}_0^\delta)$ defined via~$\hat{z}_0^\delta$ as in~\eqref{eq:limsup:rho_0} and~\eqref{eq:limsup:j_0} satisfies 
	\begin{equation*}
	(\hat{\rho}_0^\delta,\hat{\jmath}_0^\delta) \stackrel[\delta\to0]{\CE}\longrightarrow (\hat{\rho}_0,\hat{\jmath}_0) \quad\text{and}\quad 
	\limsup_{\delta\to 0}\widehat{\mathcal{I}}_0(\hat{\rho}_0^\delta,\hat{\jmath}_0^\delta) \leq  \widehat{\mathcal{I}}_0(\hat{\rho}_0,\hat{\jmath}_0).
	\end{equation*}
\end{lemma}
By a standard diagonal argument (e.g.~\cite[Rem.~1.29]{Braides02}) we can continue under the assumption that $\hat z_0$ satisfies Assumption~\ref{assump:z0-is-ct-and-pos}. The bound on the energy~\eqref{eq:th:recovery:conv-and-bound} follows from the $\delta$-independent estimate in~\eqref{lemma:limsup:limit-density:eE-bdd} below. From now on we therefore assume that Assumption~\ref{assump:z0-is-ct-and-pos} is satisfied.

\bigskip

The proof of Theorem~\ref{thm:upper-bound} now consists of  three steps.

\bigskip
\emph{Step 1: characterization of $\widehat{\mathcal{I}}_0(\hat{\rho}_0,\hat{\jmath}_0)$.} 
By Lemma~\ref{lemma:variational-problem} the limiting rate function satisfies
\begin{equation}\label{eq:limit-RF-via-density}
	\widehat{\mathcal{I}}_0(\hat{\rho}_0,\hat{\jmath}_0) = 
	\frac{1}{2} \int_{\tmspz} \hat{b}_0^2\,\hat{u}_0\, \dd y \dd t,
\end{equation}
where~$\hat{u}_0:\tmspz\to[0,\infty)$ is the function given by
\begin{equation}\label{eq:limit-density-u0}
	\hat{u}_0(t,y) = \bigl(\tfrac12 -y\bigr) \Bigl(\hat\jmath_0(t)\bigl(y+\tfrac12\bigr) + \hat z_0(t)\bigl(\tfrac12-y\bigr)\Bigr)
\end{equation}
and~$\hat{b}_0:\tmspz\to\mathbb{R}$ is defined by
\begin{equation}\label{eq:limit-function-b0}
	\hat{b}_0(t,y) := \frac{\hat{\jmath}_0(t)+\partial_y\hat{u}_0(t,y)}{\hat{u}_0(t,y)} 
	= \frac{4(\hat \jmath_0(t)-\hat z_0(t))}{\hat\jmath_0(t)(y+\tfrac12) + \hat z_0(t)
	 (\tfrac12 -y)}.
\end{equation}
The second-order polynomial~$\hat{u}_0(t,\cdot)$ is either concave ($\hat{\jmath}_0>\hat{z}_0$), linear ($\hat{\jmath}_0=\hat{z}_0$) or convex ($\hat{\jmath}_0<\hat{z}_0$). These three cases are sketched in Figure~\ref{fig:u_concance_linear_convex}.  Note that under Assumption~\ref{assump:z0-is-ct-and-pos}, $\hat b_0$ and $\partial_y \hat b_0$ are  bounded on $\tmspz$.

\begin{figure}[h!]
\labellist
\pinlabel $y$ at 950 0
\pinlabel $\hat{z}_0(t)$ at 900 500
\pinlabel $\hat{u}_0(t,y)$ at 550 750
\endlabellist
	\centering
	\includegraphics[scale=.18]{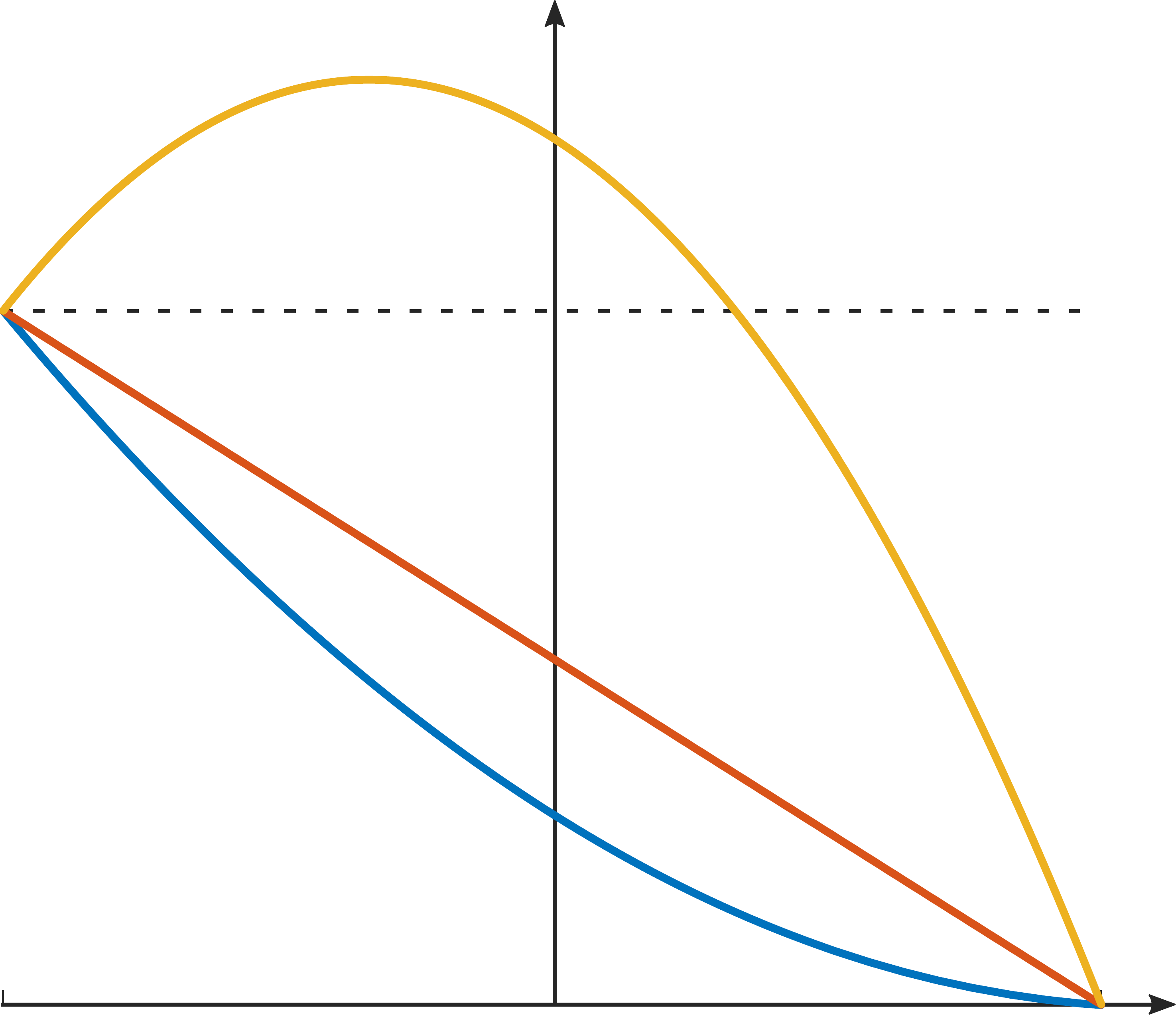}
	\caption{The polynomial $y\mapsto \hat{u}_0(t,y)$ on~$[-1/2,+1/2]$ for the three cases~$\hat{\jmath}_0(t)>\hat{z}_0(t)$ (yellow), $\hat{\jmath}_0(t)=\hat{z}_0(t)$ (red) and~$\hat{\jmath}_0(t)<\hat{z}_0(t)$ (blue). In particular, the function always satisfies $\hat{u}_0(t,-1/2) = \hat{z}_0(t)$ and~$\hat{u}_0(t,+1/2)=0$.}
	\label{fig:u_concance_linear_convex}
\end{figure}

\emph{Step 2: Solve an auxiliary PDE for $\e>0$.}
We define the function~$\hat{u}_\e^\ell:E\to[0,\infty)$ as the weak solution to the auxiliary PDE 
\begin{equation}\label{eq:steps-auxiliary-PDE}
	\hat{g}_\e^\ell \partial_t \hat{u}_\e^\ell = \partial_{yy}\hat{u}_\e^\ell - \partial_y(\hat{b}_0\bONE_{\tmspz} \hat{u}_\e^\ell),
\end{equation}
where~$\hat{g}_\e^\ell\in L^\infty(\R)$ is the Lebesgue density of the left-stationary measure~$\hat{\gamma}_\e^\ell$ from~\eqref{def:transformation-rho-hat}, that is~$\hat{\gamma}_\e^\ell(dy) = \hat{g}_\e^\ell(y)dy$.

This choice is inspired by the observation that if we define the pair~$(\hat{\rho}_\e,\hat{\jmath}_\e)$ by 
\begin{equation}
\label{eq:limsup:rho_eps-and-j_eps}	
	\hat{\rho}_\e(t,dy) := \hat{u}_\e^\ell(t,y)\hat{\gamma}_\e^\ell(dy)\quad\text{and}\quad \hat{\jmath}_\e := -\partial_y\hat{u}_\e^\ell + \hat{b}_0\bONE_{\tmspz} \hat{u}_\e^\ell,
\end{equation}
then by the characterization of weighted $L^2$-norms we have
\begin{align*}
\widehat{\cI}_\e(\hat \rho_\e,\hat\jmath_\e)
&= \sup_{b}\int_\tmsp \left[\bigl(-\partial_y\hat{u}_\e^\ell + \hat{b}_0\bONE_{\tmspz} \hat{u}_\e^\ell\bigr)\,b -  \hat{u}_\e^\ell\left(\partial_y b + \frac{1}{2}b^2\right)\right]\,\dd y\dd t\\
&= \sup_{b}\int_\tmsp \hat{u}_\e^\ell\left[b\,\hat{b}_0\bONE_{\tmspz}  -  \left(\partial_y b + \frac{1}{2}b^2\right)\right]\,\dd y\dd t\\
&=\frac12 \int_\tmspz \hat{u}_\e^\ell \hat b_0^2 ,
\end{align*}
which is an approximation of $\widehat{\cI}_0(\hat \rho_0,\hat \jmath_0)$ as given by~\eqref{eq:limit-RF-via-density}.

\medskip

We choose initial data~$\hat{u}_\e^{\ell,\circ}$ for~\eqref{eq:steps-auxiliary-PDE} that approximate $\hat\rho_0(t=0)$ in the following sense (see Lemma~\ref{l:cond-initial-data-can-be-satisfied} for a proof that such initial data can be found):
\begin{subequations}
\label{cond:initial-data-recovery-u}
\begin{align}
	&0< \hat{u}_\e^{\ell,\circ}(y) \leq 1 \quad\text{for almost all }y\in \R, 
	\label{lemma:limsup:limit-density:u-bdd} \\
	&\hat{u}_\e^{\ell,\circ}(y)\hat{\gamma}_\e^\ell(dy)\longweakto \hat{\rho}_0(0,dy)\quad\text{as}\quad \e\to 0, \label{conv-initial-density}\\
	&\int_\R \hat{u}_\e^{\ell,\circ}(y)\hat{\gamma}_\e^\ell(dy) = 1, \label{lemma:limsup:limit-density:u-mass-one}\\ 
	&\hat u_\e^{\ell,\circ} \quad\text{is constant on $(-\infty,-1/4)$ and on $(1/4,\infty)$}, 
	\label{lemma:limsup:limit-density:uz-constant}\\
	&\sup_{\e>0} \e \hat E_\e\bigl(\hat{u}_\e^{\ell,\circ}\hat{\gamma}_\e^\ell\bigr) \leq |V(x_b)| + 1.
	\label{lemma:limsup:limit-density:eE-bdd}
\end{align}
\end{subequations}
The following lemma gives the relevant properties of  $\hat u_\e^\ell$, $\hat\rho_\e$, and $\hat \jmath_\e$.

\begin{lemma}[Auxiliary PDE]\label{lemma:limsup:auxiliary-PDE}
Assume Assumption~\ref{assump:z0-is-ct-and-pos}.  For any~$\e>0$ and any initial condition~$\hat{u}_\e^{\ell,\circ}$ satisfying~\eqref{cond:initial-data-recovery-u},  there exists a solution~$\hat{u}_\e^\ell$ to the PDE~\eqref{eq:steps-auxiliary-PDE} in the following sense: $\hat{u}_\e^\ell:E\to[0,\infty)$ is such that 	
\begin{gather*}
\hat u_\e^\ell>0 \quad\text{a.e. on }\tmsp,\\
\hat{u}_\e^\ell \in  C(0,T;L^2(\tmspz)),\qquad 
\partial_y \hat{u}_\e^\ell \in L^2(0,T;L^2(\R)), 
\quad\text{and}\\
\hat \gamma_\e^\ell  \hat{u}_\e^\ell(t)   \in \mathcal P(\R) \qquad\text{for all $t$,}
\end{gather*}
and for any~$\varphi\in C_c ^1(\tmsp)$ with $\varphi = 0$ at $t=T$,
\begin{equation}
\label{eq:u-weak-form}
\int_\tmsp\Bigl[ \hat{g}_\e^\ell \hat{u}_\e^\ell \partial_t \varphi 
  + \bigl(-\partial_y \hat{u}_\e^\ell +\hat{b}_0\bONE_{\tmspz}\hat{u}_\e^\ell\bigr) \, \partial_y \varphi\Bigr]\,\dd y\dd t 
+ \int_\R \hat g^\ell_\e(y)\hat u^{\ell,\circ}_\e(y) \varphi(0,y)\, \dd y = 0.
\end{equation}
%
%
%
Define the pair~$(\hat{\rho}_\e,\hat{\jmath}_\e)$ by~\eqref{eq:limsup:rho_eps-and-j_eps}.

\medskip
Then we have 
\begin{enumerate}[label=(\roman*)]
\item \label{item:lemma-limit-density:i}$(\hat{\rho}_\e,\hat{\jmath}_\e)\in \CE(0,T)$ and
\begin{equation}\label{eq:RF-via-density}
\widehat{\mathcal{I}}_\e(\hat{\rho}_\e,\hat{\jmath}_\e) = \frac{1}{2}\int_{\tmspz}\hat{b}_0^2 \hat{u}_\e^\ell\, \dd y \dd t.
\end{equation}
\item \label{item:lemma-limit-density:eE-bound} 
$\sup_{\e>0} \e \hat E_\e(\hat \rho_\e(0,\cdot)) \leq |V(x_b)| + 1.$
\item \label{item:lemma-limit-density:rho-j-converge} 
The pair~$(\hat{\rho}_\e,\hat{\jmath}_\e)$ converges to~$(\hat{\rho}_0,\hat{\jmath}_0)$ in the sense of Definition~\ref{def:converge-in-CE}.
\item \label{item:lemma-limit-density:ii} There exists a function~$\hat{u}_0^\ell\in L^2(0,T;H^1(\Omega))$ such that 
\begin{equation}\label{eq:u_eps-converges-to-u-in-E0}
\hat{u}_\e^\ell\bONE_{\tmspz}  \xrightharpoonup{\e\to 0} \hat{u}_0 \quad\text{weakly in}\;L^2(\tmspz),
\end{equation}
\end{enumerate}
\end{lemma}

\bigskip

\emph{Step 3: Conclude.}
The convergence of $(\hat{\rho}_\e,\hat{\jmath}_\e)$  to~$(\hat{\rho}_0,\hat{\jmath}_0)$ in $\CE(0,T)$ is given by part~\ref{item:lemma-limit-density:rho-j-converge} of Lemma~\ref{lemma:limsup:auxiliary-PDE}. The energy bound~\eqref{eq:th:recovery:conv-and-bound} is satisfied by part~\ref{item:lemma-limit-density:eE-bound}, and note that this bound is independent of the regularity Assumption~\ref{assump:z0-is-ct-and-pos}.

To prove the limsup-bound~\eqref{ineq:limsup}, we observe that 
\begin{align*}
	\lim_{\e\to 0}\widehat{\mathcal{I}}_\e(\hat{\rho}_\e,\hat{\jmath}_\e) &\overset{\eqref{eq:RF-via-density}}{=} \lim_{\e\to 0} \frac{1}{2}\int_{\tmspz}\hat{b}_0^2 \hat{u}_\e^\ell\, \dd y \dd t	\\
	&\overset{\eqref{eq:u_eps-converges-to-u-in-E0}}{=} \frac{1}{2}\int_{\tmspz}\hat{b}_0^2\hat{u}_0\, \dd y \dd t \\
	&\overset{\eqref{eq:limit-RF-via-density}}{=} \widehat{\mathcal{I}}_0(\hat{\rho}_0,\hat{\jmath}_0).
\end{align*}
This concludes the proof of Theorem~\ref{thm:upper-bound}.
\end{proof}

\subsection{Proof of Lemma~\ref{lemma:limsup:auxiliary-PDE}}

We now prove the main Lemma~\ref{lemma:limsup:auxiliary-PDE} used for the proof of Theorem~\ref{thm:upper-bound}.

\medskip

\emph{Step 1: Existence of the solution $\hat u_\e^\ell$.}
Using classical methods such as those in~\cite{Lions69} one finds a function $\hat u_\e^\ell$ with
\[
\hat u_\e^\ell(t)\geq0 \qquad \text{and}\qquad 
\int _\R \hat u_\e^\ell(t,y)\hat \gamma_\e^\ell(\dd y)  = 1\qquad\text{for all }t,
\]
that  satisfies
 the $\e$-independent bounds
\begin{subequations}
\label{props:reg-u}
\begin{align}
&\|\hat u_\e^\ell\|_{C([0,T];L^2(\tmspz))} \leq C,
\label{props:reg-u:CL2}\\	
&\|\partial_y\hat u_\e^\ell\|_{L^2(\tmsp)} \leq C,
\label{props:reg-u:L2H1}\\	
&\|\hat g_\e^\ell \partial_t \hat u_\e^\ell\|_{L^2(0,T; H^{-1}(\R))} \leq C,
\label{props:reg-u:L2Hm1}
\end{align}
\end{subequations}
and solves equation~\eqref{eq:steps-auxiliary-PDE} in the weak form~\eqref{eq:u-weak-form}.

To briefly indicate the main steps in this existence proof, define the function $B(t,y) := \int_{-1/2}^y \hat b_0(t,\tilde y)\,d\tilde y$ and observe that the transformed function $\hat{v}_\e^\ell := e^{-B}\hat{u}_\e^\ell$ satisfies the equation
\begin{equation}\label{eq:limsup:eq-for-v_eps}
	\hat{g}_\e^\ell \partial_t \left(e^{B}\hat{v}_\e^\ell\right) = \partial_y\left(e^{B}\partial_y \hat{v}_\e^\ell\right).
\end{equation}
Applying the usual method of multiplying by the solution $\hat v_\e^\ell$ and integrating we obtain this \emph{a priori} estimate:
\begin{align}
\int_0^T\int_\mathbb{R} e^{B} |\partial_y\hat{v}_\e^\ell|^2\, \dd y \dd t + \sup_{t\in[0,T]}\int_\mathbb{R} e^{B(t)}\hat{g}_\e^\ell \hat{v}_\e^\ell(t)^2\, \dd y 
&\leq \int_\R e^{-B(0)} \hat g_\e^\ell (\hat u_\e^{\ell,\circ})^2 \notag \\
&\leftstackrel{\eqref{lemma:limsup:limit-density:u-bdd}}\leq \|e^{-B(0)}\|_\infty.
\label{eq:limsup:energy-bound-v_eps}
\end{align}
One then constructs by e.g.\ Galerkin approximation a sequence of approximating solutions of~\eqref{eq:limsup:eq-for-v_eps} that satisfy~\eqref{eq:limsup:energy-bound-v_eps}, for which one can extract a subsequence that converges to a limit. Upon transforming back to the function~$\hat u_\e^\ell$ one obtains the weak form~\eqref{eq:u-weak-form} and the bounds~\eqref{props:reg-u:L2H1} and~\eqref{props:reg-u:L2Hm1}.

In order to deduce~\eqref{props:reg-u:CL2} from~\eqref{props:reg-u:L2H1} and~\eqref{props:reg-u:L2Hm1} one applies e.g.~\cite[Th.~5]{Simon87} with the compact embedding $H^1(\tmspz)\hookrightarrow L^2(\tmspz)$. The missing $L^2(\tmspz)$-estimate can be obtained from~\eqref{props:reg-u:L2H1} by applying the generalized Poincar\'e inequality of Lemma~\ref{l:gen-Poincare} to $\mu = \hat \gamma_\e^\ell$ and observing that $\hat\gamma_\e^\ell([-1/2,1/2])\to \infty$ as $\e\to0$. 

\smallskip
By the strong maximum principle and the positivity~\eqref{conv-initial-density} of the initial data the solutions $\hat u_\e^\ell$ are strictly positive, and since $\hat \jmath_\e\in L^2(\tmsp)$ the mass of $\hat \rho_\e(t) = \gamma_\e^\ell \hat u_\e^\ell(t)$ equals the mass of the initial data $\hat \gamma_\e^\ell$, which is one by~\eqref{lemma:limsup:limit-density:u-mass-one}.

\smallskip

Note that by Assumption~\ref{assump:z0-is-ct-and-pos} the function $B$ is not only bounded but also independent of~$\e$, implying that the constants in~\eqref{props:reg-u} also are independent of $\e$.

\medskip

\emph{Step 2: Part~\ref{item:lemma-limit-density:i}, the value of\/ $\widehat \cI_\e(\hat\rho_\e,\hat \jmath_\e)$.} The fact that~$(\hat{\rho}_\e,\hat{\jmath}_\e)\in\CE(0,T)$ follows from the regularity~\eqref{props:reg-u} of $\hat u_\e^\ell$ and from the weak form~\eqref{eq:u-weak-form} of the equation.  The value of $\widehat \cI_\e(\hat\rho_\e,\hat \jmath_\e)$ was already calculated before Lemma~\ref{lemma:limsup:auxiliary-PDE}.

\medskip

\emph{Step 3: Convergence of $(\hat \rho_\e,\hat \jmath_\e)$.}  
By construction (see~\eqref{conv-initial-density}) the initial measures $\hat\rho_\e(0,\cdot)$ converge to $\hat\rho_0(0,\cdot)$. To prove convergence of $(\hat \rho_\e,\hat \jmath_\e)$ we therefore need to show convergence in the continuity equation. 

For any test function~$\varphi\in C_b(\tmsp)$,
\begin{align*}
	\left|\int_\tmsp \varphi \hat{\rho}_\e\right|^2 \;&\leftstackrel{\eqref{eq:limsup:rho_eps-and-j_eps}} = \left|\int_\tmsp\varphi\, e^{B}\hat{g}_\e^\ell \hat{v}_\e^\ell\, \dd y \dd t\right|^2\\
	&\leftstackrel{\mathrm{CS}}{\leq} \left(\int_\tmsp e^{B}\hat{g}_\e^\ell \hat{v}_\e^\ell(t)^2\, \dd y \dd t\right) \left(\int_\tmsp |\varphi e^{B}\hat{g}_\e^\ell| \, \dd y \dd t\right)\\
	&\leftstackrel{\eqref{eq:limsup:energy-bound-v_eps}}{\leq} C \int_{\text{supp}(\varphi)}\hat{g}_\e^\ell(y)\, \dd y.
\end{align*}
Hence for any test function with support outside of~$[0,T]\times \{\pm 1/2\}$, by Lemma~\ref{l:transformed-stationary-measure},
\begin{equation}
\label{conv:hrhoz-to-0}	
	\int_\tmsp \varphi \hat{\rho}_\e \xrightarrow{\e\to 0} 0.
\end{equation}

Take any sequence $\e_k\to0$. By~\eqref{conv:hrhoz-to-0} the family of measures~$\hat{\rho}_{\e_k}$ is tight, and therefore it converges weakly on $[0,T]\times \R$, along a subsequence (denoted the same),  to a measure~$\overline{\rho}_0$ that is concentrated on~$[0,T]\times\{\pm 1/2\}$, and therefore has the form
\[
\overline \rho_0(t,\dd y) = \overline z_0(t) \delta_{-1/2}(\dd y) + (1-\overline z_0(t)) \delta_{1/2}(\dd y)
\]
for some measurable function $\overline z_0:[0,T]\to [0,1]$.
\smallskip

Since the function~$B$ is bounded, we find  that~$\hat{\jmath}_\e= -e^{B}\partial_y\hat{v}_\e^\ell$ is bounded in~$L^2(\tmsp)$, because
\begin{equation*}
	\int_\tmsp |\hat{\jmath}_\e|^2\, \dd y \dd t \leq \|e^B\|_\infty \int_\tmsp e^{B} |\partial_y\hat{v}_\e^\ell|^2\, \dd y \dd t \overset{\eqref{eq:limsup:energy-bound-v_eps}}{\leq} C.
\end{equation*}
Hence, taking another subsequence, the flux $\hat\jmath_{\e_k}$ converges weakly in~$L^2(\tmsp)$  to some~$\overline {\jmath}_0 \in L^2(\tmsp)$. 

Combining these convergence statements of~$\hat{\rho}_{\e_k}$ and~$\hat{\jmath}_{\e_k}$, we find for any test function~$\varphi\in C^1_c(\tmsp)$,
\begin{equation*}
	0 \overset{\mathrm{CE}}{=} \int_\tmsp \bigl[\partial_t \varphi \, \hat{\rho}_{\e_k} +\partial_y\varphi  \, \hat{\jmath}_{\e_k} \bigr] \xrightarrow{k \to \infty} \int_\tmsp \bigl[\partial_t \varphi \, \overline\rho_0 + \partial_y\varphi \,\overline{\jmath}_0\bigr].
\end{equation*}
Therefore  $(\hat \rho_{\e_k},\hat \jmath_{\e_k})$ converges to $(\overline \rho_0,\overline\jmath_0)$ in the sense of $\CE(0,T)$. 

\smallskip
Finally, since~$\overline{\rho}_0$ is concentrated on~$[0,T]\times\{\pm1/2\}$, the limiting flux $\overline \jmath_0$ is piecewise constant in $y$ with jumps only at~$\{\pm1/2\}$, and $\overline{\jmath}_0\in L^2(\tmsp)$ implies that $\overline\jmath_0$  vanishes outside of~$(-1/2,+1/2)$. Therefore, the continuity equation~$0=\partial_t\overline{\rho}_0 + \partial_y\overline{\jmath}_0$ in the distributional sense implies that the flux is given by
\begin{equation}
\label{eq:connection-j0-z0}	
\overline{\jmath}_0(t,\dd y)= -\partial_t\overline{z}_0(t)\bONE_{(-1/2,+1/2)}(y)\,\dd y.
\end{equation}

\medskip

\emph{Step 4: The limit $(\overline \rho_0,\overline\jmath_0)$ is equal to $(\hat\rho_0,\hat\jmath_0)$.} We now show that the limit $\overline z_0$ obtained above coincides with the function $\hat z_0$ that characterizes $\hat \rho_0$ (see~\eqref{eq:limsup:rho_0}). This proves that $(\overline \rho_0,\overline \jmath_0) = (\hat \rho_0,\hat \jmath_0)$ and $\overline u_0 = \hat u_0$ on $\tmspz$. 

By further extracting subsequences we can assume that
\[
\hat u_{\e_k}^\ell \weakto \overline u_0 \quad\text{in }L^2(\tmspz)
\qquad\text{and}\qquad
\partial_y \hat u_{\e_k}^\ell \weakto \partial_y \overline u_0 \quad\text{in }L^2(\tmsp).
\]
By passing to the limit in~\eqref{eq:limsup:rho_eps-and-j_eps} we find that~$\overline{\jmath}_0 = -\partial_y\overline{u}_0 + \hat{b}_0\bONE_{\tmspz}\overline{u}_0$ almost everywhere in $\tmsp$. In combination with~\eqref{eq:connection-j0-z0} this means  that for almost every $t\in [0,T]$, the function~$y\mapsto \overline{u}_0(t,y)$ is a  weak solution of the ODE
\begin{equation}
\label{eq:ODE-u_0}	
	-\partial_y\overline{u}_0(t,y) + \hat{b}_0(t,y)\overline{u}_0(t,y) = -\partial_t\overline{z}_0(t), \quad\text{for }-\frac12 < y < \frac12 .
\end{equation}

This is a first-order ODE in $y$ on the interval $(-1/2,1/2)$, and we show below that $\overline u$ satisfies not one but \emph{two} boundary conditions, at $\pm1/2$:
\begin{align}
	\overline{u}_0(t,-1/2)= \overline{z}_0(t)
	\qquad\text{and}\qquad \overline{u}_0(t,+1/2) = 0 \qquad\text{for a.e. }t.
	\label{bc:u_0}
\end{align}
The solution of~\eqref{eq:ODE-u_0} with left boundary condition $\overline{u}_0(t,-1/2)= \overline{z}_0(t)$ is given by
\begin{equation*}
	\overline{u}_0(t,y) = e^{B(t,y)}\left[\overline{z}_0(t) + \partial_t\overline{z}_0(t) \int_{-1/2}^ye^{-B(t,z)}\,dz\right].
\end{equation*}
Since
\begin{equation*}
	\int_{-1/2}^{+1/2}e^{-B(t,z)}\,dz = -\frac{\hat{z}_0(t)}{\partial_t\hat{z}_0(t)} > 0,
\end{equation*}
the second boundary condition~$\overline{u}_0(t,+1/2)=0$ therefore enforces
\begin{equation}\label{eq:upper-bound:proof-lemma-limiting-density}
	\partial_t \log\overline{z}_0(t) = \partial_t\log\hat{z}_0(t).
\end{equation}
Combined with the convergence assumption on the initial condition~$\hat{\rho}_\e(0,dy)$, which implies~$\overline{z}_0(0)=\hat{z}_0(0)$, it follows that~$\overline{z}_0 = \hat{z}_0$. This unique characterization of the limit $(\overline \rho,\overline \jmath_0)$ also implies that the convergence holds not only along subsequences but in the sense of a full limit $\e\to0$. 

\medskip

\emph{Step 5: Prove the boundary conditions~\eqref{bc:u_0} on $\overline u_0$.}
To prove the left boundary condition in~\eqref{bc:u_0}, let~$U_\delta$ be a small neighborhood around~$-1/2$ of length~$2\delta>0$. Since~$\partial_y\hat{u}_\e^\ell$ is  bounded in~$L^2(\tmsp)$ by~\eqref{props:reg-u:L2H1}, there is an $\alpha\in L^2(0,T)$ such that
\begin{equation*}
	\hat{u}_\e^\ell(t,y)\leq \hat{u}_\e^\ell(t,-1/2) + \alpha(t)\,|y+1/2|^{1/2}\qquad\text{for all $\e$ and a.e. $(t,y)\in \tmspz$}.
\end{equation*}
We can then estimate for any non-negative $\psi\in C([0,T])$,
\begin{multline*}
	\int_0^T\psi(t)\hat{\rho}_\e(t,U_\delta)\, \dd t = \int_0^T\int_{U_\delta} \psi(t)\hat{u}_\e^\ell(t,y)\hat\gamma_\e^\ell(\dd y)\dd t\\
	\leq \hat \gamma_\e^\ell( U_\delta)\int_0^T\psi(t)\hat{u}_\e^\ell(t,-1/2)\, \dd t  + \|\alpha\|_2 \|\psi\|_{L^\infty} \int_{U_\delta} |y+1/2|^{1/2}\hat{\gamma}_\e^\ell(\dd y)\\
	\leq \hat \gamma_\e^\ell( U_\delta) \int_0^T\psi(t)\hat{u}_\e^\ell(t,-1/2)\, \dd t  + C \|\psi\|_{L^\infty} \delta^{1/2} \int_{U_\delta}\hat{g}_\e^\ell(y)\, \dd y.
\end{multline*}
For each~$\delta>0$,~$\int_{U_\delta}\hat{g}_\e^\ell(y)dy$ converges to $1$ as~$\e\to 0$, and
\begin{equation*}
	\lim_{\e\to 0}\int_0^T\psi(t)\hat{\rho}_\e(t,U_\delta)\, \dd t = \int_0^T\psi(t)\overline{z}_0(t)\, \dd t.
\end{equation*}
Therefore,
\begin{equation*}
	\int_0^T\psi(t)\overline{z}_0(t)\, \dd t \leq \liminf_{\e\to 0}\int_0^T\psi(t)\hat{u}_\e^\ell(t,-1/2)\, \dd t + C'\delta^{1/2}.
\end{equation*}
Noting that~$\delta>0$ is arbitrary and repeating the argument for the reversed inequality, we find that
\begin{equation*}
	\int_0^T\psi(t)\overline{z}_0(t)\, \dd t = \lim_{\e\to 0}\int_0^T\psi(t)\hat{u}_\e^\ell(t,-1/2)\, \dd t.
\end{equation*}
Since the trace map $w\in L^2(0,T;H^1(\tmspz))\mapsto w(\cdot,-1/2)\in L^2(0,T)$ is weakly continuous, the sequence of functions $t\mapsto \hat{u}_\e^\ell(t,-1/2)$ converges weakly in $L^2(0,T)$ to the limit $\overline u_0(t,-1/2)$. 
This proves the first boundary condition in~\eqref{bc:u_0}. The argument for the second boundary condition is similar, using that~$\hat{\gamma}_\e^\ell\bigl((1/2-\delta,1/2+\delta)\bigr)\to\infty$ as~$\e\to 0$.

\smallskip
This concludes the proof of Lemma~\ref{lemma:limsup:auxiliary-PDE}.

\subsection{Proof of Lemma~\ref{lemma:upperbound:energy-dense-set-of-densities}}
Approximation results of this type are very common; see e.g.~\cite[Theorem~6.1]{ArnrichMielkePeletierSavareVeneroni2012} or~\cite[Lemma~4.7]{PeletierRenger20TR}. Fix a pair $(\hat\rho_0,\hat\jmath_0)$ with $\widehat\cI_0(\hat\rho_0,\hat\jmath_0)<\infty$, and write $\hat \rho_0$ in terms of the absolutely continuous function $\hat z_0$ as in~\eqref{eq:limsup:rho_0}.

We first approximate $\hat z_0$ by a sequence of more regular functions $\hat z_\eta$, for $\eta\to0$. We do this by first extending $\hat z_0$ to $\R$ by constants:
\[
\hat z_0(t) := \begin{cases}
\hat z_0(0)	& \text{if }t\leq 0\\
\hat z_0(t) & \text{if }0\leq t\leq T\\
\hat z_0(T)  &\text{if } t\geq T.
\end{cases}
\]
The extended function $\hat z_0$ again is non-increasing; we then regularize by convolution  by setting
\[
\hat z_\eta := \alpha_\eta {*} \hat z_0,
\]
where $\alpha_\eta(s) := \eta^{-1}\alpha(s/\eta)$ is a regularizing sequence. 

Then $\hat z_\eta\to\hat z_0$ in $W^{1,1}(\R)$, and therefore the corresponding pair $(\hat\rho_\eta,\hat\jmath_\eta)$ converges in $\CE$ to $(\hat\rho_0,\hat\jmath_0)$
Since the function $S$ in~\eqref{eq:S-fct} is jointly convex in its two arguments, we have
\begin{align*}
\lefteqn{\int_0^T S(-\partial_t \hat z_\eta(t)|\hat z_\eta(t))\, \dd t
\leq \int_\R S(-\partial_t \hat z_\eta(t)|\hat z_\eta(t))\, \dd t}\qquad &
\\
&\leq \int_\R \bigl(\alpha_\eta {*} S(-\partial_t \hat z_0|\hat z_0)\bigr)(t)\, \dd t
= \int_\R S(-\partial_t \hat z_0(t)|\hat z_0(t))\, \dd t\\
&= \int_0^T S(-\partial_t \hat z_0(t)|\hat z_0(t))\, \dd t.
\end{align*}

Next, define $\overline z(t) := 1/2 -t/{4T}$, and note that $\overline z$ and $-\partial_t \overline z$ are bounded away from zero on $[0,T]$. For each $\eta\in(0,1)$, the convex combination
\[
\wt z_\eta(t) := \eta \overline z(t)  + (1-\eta) \hat z_\eta(t), \qquad t\in [0,T].
\]
also satisfies $\inf \wt z_{\eta}$, $\inf (-\partial_t \wt z_{\eta}) > 0$. Again using the convexity of $S$ we find that 
\[
\int_0^T S(-\partial_t \wt z_{\eta}(t)|\wt z_{\eta}(t))\, \dd t
\leq C\eta + (1-\eta) \int_0^T S(-\partial_t \hat z_0(t)|\hat z_0(t))\, \dd t.
\]
Setting $\hat\rho_0^{\eta} (t) = \wt z_{\eta}(t)\delta_{-1/2}  + (1-\wt z_{\eta}(t))\delta_{1/2}$ and defining $\hat \jmath_0^{\eta}$ accordingly, we then have 
\[
\hat \cI_0(\hat \rho_0^{\eta} , \hat \jmath_0^{\eta}) 
\leq C  \eta + (1-\eta) \hat \cI_0(\hat \rho_0 , \hat \jmath_0) .
\]
The sequence $(\hat \rho_0^{\eta} , \hat \jmath_0^{\eta})$ therefore satisfies the claim of Lemma~\ref{lemma:upperbound:energy-dense-set-of-densities}.

\subsection{The initial data in~\eqref{cond:initial-data-recovery-u} can be realized}

In the proof of Theorem~\ref{thm:upper-bound} we postulated a choice of initial data with certain properties. The next lemma shows that it possible to construct such initial data. 

\begin{lemma}
\label{l:cond-initial-data-can-be-satisfied}
For any given $\rho^\circ = z^\circ \delta_{-1/2} + (1-z^\circ)\delta_{1/2}$ it is possible to choose a sequence $\hat u_\e^{\ell,\circ}$ satisfying the requirements~\eqref{cond:initial-data-recovery-u}.
\end{lemma}

\begin{proof}
For instance one may choose
\[
\hat u_\e^{\ell,\circ}(y) := 
\begin{cases}
z_0(0) & \text{if }y\leq -1/4\\
\text{smooth monotonic } & \text{between $-1/4$ and $1/4$}\\
\ds(1-z_0(0) + a_\e)\frac{Z_\e^\ell}{Z_\e} & \text{if } y\geq 1/4,
\end{cases}
\]
where $a_\e\to0$ can be tuned in order to achieve the mass constraint~\eqref{lemma:limsup:limit-density:u-mass-one}.
One can verify that the definitions of $\gamma_\e$ and~$\gamma_\e^\ell$ imply that  $1-z_0+a_\e\leq 1$, and because ${Z_\e^\ell}/{Z_\e}<1$  we have the bound $\|\hat u_\e^{\ell,\circ}\|_\infty \leq 1$. 

To show~\eqref{lemma:limsup:limit-density:eE-bdd} for this choice we can write 
\[
\e \hat E_\e(\hat u_\e^{\ell,\circ} \hat\gamma^\ell_\e) 
= \e \int_\R \eta\Bigl( \hat u_\e^{\ell,\circ} \frac{Z_\e}{Z_\e^\ell}\Bigr) \dd\hat \gamma_\e 
= \e \int_\R \hat u_\e^{\ell,\circ} \frac{Z_\e}{Z_\e^\ell} \log \Bigl( \hat u_\e^{\ell,\circ} \frac{Z_\e}{Z_\e^\ell}\Bigr) \dd\hat \gamma_\e .
\]
Splitting the integral into parts, the integral over $(1/4,\infty)$ equals
\[
\e  (1-z_0(0)+a_\e)\log (1-z_0(0)+a_\e) \hat\gamma_\e  \bigl((1/4,\infty)\bigr) \leq  0.
\]
The integral over the remaining interval $(-\infty,1/4)$ can be bounded from above by 
\[
\e \hat \rho_\e\bigl((-\infty,1/4)\bigr) \log \Bigl(\|u_\e^{\ell,\circ}\|_\infty \frac{Z_\e}{Z_\e^\ell}\Bigr) 
\leq \e \log \frac {Z_\e}{Z_\e^\ell}\stackrel{\eqref{char:asymp-Ze}} \leq |V(x_b)| + 1 \quad\text{for small }\e.
\qedhere
\]
\end{proof}

\subsection{Recovery sequence for the untransformed system}

\begin{theorem}
\label{t:recovery-sequence-untransformed}
Let $V$ satisfy Assumption~\ref{ass:V}. Let $(\rho_0,j_0)\in \CE(0,T)$ satisfy $\cI_0(\rho_0,j_0)< \infty$. Then there exists a sequence $(\rho_\e,j_\e)\in \CE(0,T)$ such that $(\rho_\e,j_\e)\stackrel{\mathrm{CE}}\longrightarrow (\rho_0,j_0)$, $\sup_{\e>0}\e E_\e(\rho_\e(0))<\infty$, and $\cI_\e(\rho_\e,j_\e) \longrightarrow \cI_0(\rho_0,j_0)$.
\end{theorem}

\begin{proof}
Since $\cI_0(\rho_0,j_0)< \infty$, $\rho_0$ and $j_0$ have the structure~\eqref{structure:rho-zero-j-zero} in terms of $z$ and $j$. Define the corresponding $(\hat \rho_0,\hat \jmath_0)$ by 
\begin{align*}
\hat\rho_0(t,\dd x) &:= z(t)\delta_{-1/2}(\dd x) + (1-z(t))\delta_{1/2}(\dd x),\\
\hat \jmath_0(t,\dd x) &:= j(t)\bONE_{[-1/2,1/2]}(x)\, \dd x.
\end{align*}
By construction $\widehat\cI_0(\hat \rho_0,\hat \jmath_0) = \cI_0(\rho_0,j_0)<\infty$, and therefore by Theorem~\ref{thm:upper-bound} there exists a sequence $(\hat\rho_\e,\hat \jmath_\e)$ that converges to $(\hat \rho_0,\hat \jmath_0)$ with $\hat \cI_\e(\hat\rho_\e,\hat \jmath_\e)\longrightarrow \cI_0(\rho_0,j_0)$. 

We define $(\rho_\e,j_\e)$ by back-transforming the relation~\eqref{def:transformation}: 
\[
\rho_\e := \bigl((y_\e)^{-1}\bigr)_\# \hat\rho_\e 
\qquad\text{and}\qquad
j_\e(t,x) := \hat \jmath_\e (t,y_\e(x)).
\]
By definition then $\cI_\e(\rho_\e,j_\e) = \widehat\cI_\e(\hat\rho_\e,\hat\jmath_\e)\longrightarrow \cI_0(\rho_0,j_0)$. The only remaining fact to check is the convergence $(\rho_\e,j_\e)\stackrel{\mathrm{CE}}\longrightarrow (\rho_0,j_0)$.

By Theorem~\ref{thm:upper-bound}, $\cI_\e(\rho_\e,j_\e)$ and $\e E_\e(\rho_\e(0))$ are bounded.
We next verify the convergence~\eqref{ass:initial-data-rho_e} of the initial data. Note that by the properties of push-forwards, 
\begin{equation}
\label{eq:transformation-hatu-u}
\rho_\e (0,\dd x) = u_\e^{\ell,\circ}(x)\gamma_\e^\ell (\dd x)
\qquad\text{with}\qquad
u_\e^{\ell,\circ}(x)  = \hat u_\e^{\ell,\circ}(y_\e(x)).
\end{equation}
\begin{lemma}
\label{l:props-u-e}
\begin{enumerate}
\item \label{l:props-u-e:bdd}
$u_\e^{\ell,\circ}$ is bounded uniformly in $x$ and $\e$;
\item \label{l:props-u-e:const-left}
for small $\e$, on the interval $(-\infty,\tfrac12 (x_a+x_0))$ the function $u_\e^{\ell,\circ}$ is equal to a constant $a_\e$, with $\lim_{\e\to0} a_\e = z_0(0)$;
\item \label{l:props-u-e:const-right}for small $\e$, on the interval $(\tfrac12 (x_0+x_{b-}),\infty)$ the function $u_\e^{\ell,\circ} Z_\e/Z_\e^\ell$ is equal to a constant $b_\e$, with $\lim_{\e\to0} b_\e = 1-z_0(0)$.
\end{enumerate}
\end{lemma}
\noindent
Assuming this lemma for the moment, we calculate for any $\varphi\in C_b(\R)$ that 
\[
\int_{-\infty}^{\tfrac12(x_a+x_0)} \rho_\e(0,\dd x ) \varphi(x) 
= a_\e \int_{-\infty}^{\tfrac12(x_a+x_0)} \gamma_\e^\ell(\dd x ) \varphi(x) 
\stackrel{\e\to0} \longrightarrow z(0)\varphi(x_a).
\]
Similarly, 
\[
\int_{\tfrac12 (x_0+x_{b-})}^\infty  \rho_\e(0,\dd x ) \varphi(x) 
= b_\e  \int_{\tfrac12 (x_0+x_{b-})}^\infty \gamma_\e(\dd x ) \varphi(x) 
\stackrel{\e\to0} \longrightarrow (1-z(0))\varphi(x_b).
\]
Finally, by the uniform boundedness of $u_\e^{\ell,\circ}$ on $\R$, 
\[
\biggl|\int_{\tfrac12(x_a+x_0)}^{\tfrac12 (x_0+x_{b-})} \rho_\e(0,\dd x ) \varphi(x) \biggr|
\leq C\|\varphi\|_\infty \int_{\tfrac12(x_a+x_0)}^{\tfrac12 (x_0+x_{b-})} \gamma_\e^\ell(\dd x)
\stackrel{\e\to0} \longrightarrow 0.
\]
Therefore $\rho_\e(0,\cdot)$ satisfies the convergence condition~\eqref{ass:initial-data-rho_e}.
Theorem~\ref{th:compactness-lower-bound-untransformed} then implies that up to extraction of a subsequence, $(\rho_\e,j_\e)$ converges to a limit $(\overline \rho_0,\overline \jmath_0)$; the only property to check is that $(\overline \rho_0,\overline \jmath_0) = (\rho_0,j_0)$.

Let  $\overline \rho_0(t) = \overline z_0(t) \delta_{x_a} + (1-\overline z_0(t)) \delta_{x_b}$; by~\eqref{ass:initial-data-rho_e} we have  $\overline z_0(0) = z_0(0)$. Recall from Lemma~\ref{l:props-phie-ye} that the function $\phi_\e\circ y_\e^{-1}$ converges uniformly on $\R$ to the function $\idhalf$. We then calculate for any $\psi\in C_b([0,T])$ that 
\begin{multline*}
\int_0^T \psi(t) \int_\R \rho_\e(t,\dd x) \phi_\e(x) \, \dd t 
= \int_0^T \psi(t) \int_\R \hat\rho_\e(t,\dd y) \phi_\e\bigl(y_\e^{-1}(y)\bigr) \, \dd t\\
\longrightarrow \int_0^T \psi(t) \int_\R \hat \rho_0(t,\dd y) \idhalf(y)\, \dd t
= \int_0^T \psi(t) \bigl[ -\tfrac12 z_0(t) + \tfrac12 (1-z_0(t)) \bigr] \, \dd t.
\end{multline*}
On the other hand, since $\phi_\e$ is uniformly bounded and converges to $\mp 1/2$ in neighbourhoods of $x_a$ and $x_b$, we also have
\begin{align*}
\int_0^T \psi(t) \int_\R \rho_\e(t,\dd x) \phi_\e(x) \, \dd t 
&\longrightarrow\int_0^T \psi(t) \bigl[ -\tfrac12 \overline z_0(t) + \tfrac12 (1-\overline z_0(t)) \bigr] \, \dd t.
\end{align*}
Since these two should agree for all $\psi\in C_b([0,T])$, it follows that $\overline z_0 = z_0$ and therefore $\overline \rho_0 = \rho_0$.

Finally, to show that also $\overline \jmath_0 = j_0$, note that both $\overline \jmath_0$ and $j_0$ are of the form $j(t) \bONE_{[x_a,x_b]}$, and since they satisfy the continuity equation with the same measure $\rho$ we have $\partial_y(\overline \jmath_0-j_0) = 0$ in duality with $C_c^{1,0}([0,T]\times \R)$. It follows that $\overline \jmath_0 = j_0$ almost everywhere on $[0,T]\times \R$. 
\end{proof}

We still owe the reader the proof of Lemma~\ref{l:props-u-e}.

\begin{proof}[Proof of Lemma~\ref{l:props-u-e}]
Part~\ref{l:props-u-e:bdd} follows directly from the boundedness of $\hat u_\e^{\ell,\circ}$ (see~\eqref{lemma:limsup:limit-density:u-bdd}) and the transformation~\eqref{eq:transformation-hatu-u}. 
For part~\ref{l:props-u-e:const-left}, recall from~\eqref{lemma:limsup:limit-density:uz-constant} that $\hat u_\e^{\ell,\circ}$ is a constant (say $a_\e$) on the  interval $(-\infty,-1/4)$. Since $\hat u_\e^{\ell,\circ}\hat\gamma_\e^\ell$ converges to $z_0(0)\delta_{-1/2} + (1-z_0(0))\delta_{1/2}$, the constant $a_\e$ converges to $z_0(0)$. 
 Since $y=-1/2$ is an interior point of the interval $(-\infty,-1/4)$, for sufficiently small $\e$ the function $y_\e$ maps the interval $\bigl(-\infty, \tfrac12(x_a+x_0)\bigr)$ into $(-\infty,-1/4)$ (see Lemma~\ref{l:props-phie-ye}) and therefore $u_\e^{\ell,\circ}$ equals $a_\e$ on $\bigl(-\infty, \tfrac12(x_a+x_0)\bigr)$. 

For part~\ref{l:props-u-e:const-right} the argument is very similar, only replacing the left-normalized~$\hat\gamma_\e^\ell$ by the standard normalized $\hat \gamma_\e$.
\end{proof}

\appendix
\section{Auxiliary results}

\subsection{A generalized Poincar\'e inequality}

\begin{lemma}
\label{l:gen-Poincare}
For any $-\infty<a<b<\infty$ and for all bounded non-negative Borel measures $\mu$ on $[a,b]$, we have the following inequality:
\begin{multline*}
\|f\|^2_{L^\infty(a,b)}
\leq C\Bigl(  \|f'\|^2_{L^2(a,b)}+ \mu([a,b])^{-1} \|f\|^2_{L^2(a,b;\mu)} \Bigr),\\
\text{for all }f\in H^1(a,b)\cap L^2(a,b;\mu),
\end{multline*}
with a constant $C>0$ that only depends on $a$ and $b$.
\end{lemma}

\begin{proof}
By density it suffices to prove the inequality for $f\in C^1([a,b])$. For $x,y\in [a,b]$ we have 
\[
f^2(x)  = f^2(y) + 2\int_x^y f'(s) f(s)\, \dd s,
\]
and therefore
\begin{align*}
\mu([a,b])f^2(x) &\leq  \|f\|_{L^2(\mu)}^2  + 2\mu([a,b]) \|f\|_{L^2(a,b)}\|f'\|_{L^2(a,b)}\\
&\leq  \|f\|_{L^2(\mu)}^2  + \mu([a,b])\biggl\{\frac1 \alpha  \|f'\|_{L^2(a,b)}^2 + \alpha(b-a) \|f\|_{L^\infty(a,b)}^2\biggr\}.
\end{align*}
The assertion follows by choosing e.g.\ $\alpha = 1/2(b-a)$.
\end{proof}

\subsection{Laplace's method}
\begin{lemma}[Laplace's method; see e.g.~{\cite[Sec.~7.2]{Olver74}}]\label{lemma:watson}
	Let $f : [a,b] \to \mathbb{R}$ be twice differentiable. 
\begin{enumerate}[(a)]
\item 	\label{char:Laplace:interior}
Suppose that for some $x_i \in (a,b)$, we have $f(x_i) = \inf_{[a,b]} f$. Then
	\begin{equation*}
	\int_a^b e^{-nf(x)}dx = \left[ 1 + o(1)\right] \sqrt{\frac{2\pi}{n f''(x_i)}} e^{-n f(x_i)} \quad \text{as }n \to \infty.
	\end{equation*}
\item \label{char:Laplace:boundary}
If $x_i = a$ or $x_i=b$, then
	\begin{equation*}
	\int_a^b e^{-nf(x)}dx = \left[ 1 + o(1)\right] \frac{1}{2} \sqrt{\frac{2\pi}{n f''(x_i)}} e^{-n f(x_i)} \quad\text{as } n \to \infty.
	\end{equation*}
\end{enumerate}

\end{lemma}

\subsection{Duality characterizations}

\begin{lemma}[Duality characterization of quadratic entropies]
\label{lemma:appendix:dual-of-convex-functions}
	For $X=[0,T]\times \mathbb{R}^d$, $f,g: X \to \mathbb{R}$ measurable with $g > 0$, and any nonnegative Borel measure $\mu$, we have the characterization
	\begin{equation*}
	\int_X \frac{1}{2} \frac{|f(x)|^2}{g(x)} \, \dd \mu(x) = \sup_{\substack{b \in C_c^\infty(X)}} \int_X \left[\left(-\frac{b(x)^2}{2}\right)g(x) + b(x)f(x)\right]\,\dd\mu(x).
	\end{equation*}
\end{lemma}
A proof is given for instance in~\cite[Lemma~3.4]{ArnrichMielkePeletierSavareVeneroni2012}. The representation  there can be further simplified by setting $a=-b^2/2$.

\begin{lemma}[Duality characterization of $S$]\label{lemma:variational-problem}
	The function~$S$ defined in~\eqref{eq:S-fct} has the alternative characterization 
	\begin{equation}\label{eq:lemma:var-probl}
		S(j|z) = \frac{1}{4} \inf_u \int_{-1/2}^{+1/2} \frac{1}{u(y)}\big|j + u'(y)\big|^2\,\dd y,
		\end{equation}
		where the infimum is taken over smooth functions~$u:[-1/2,+1/2]\to[0,\infty)$ satisfying the boundary conditions $u(-1/2)=z$ and~$u(+1/2)=0$ and the positivity requirement $u'(y)>0$ for all $y\in (-1/2,1/2)$.
		The optimal function $u$ is the polynomial
		\begin{equation}\label{eq:lemma-var-probl:optimal-u}
			u(y) = \bigl(\tfrac12 -y\bigr) \Bigl(j\bigl(y+\tfrac12\bigr) + z\bigl(\tfrac12-y\bigr)\Bigr). 
		\end{equation}
\end{lemma}
\begin{proof}
This result is very similar to that in~\cite[Prop.~A1]{LieroMielkePeletierRenger17}, which dealt with the slightly different argument $(j^2+{u'}^2)/u$ with strictly positive boundary conditions for $u$; the sign of $j$ and the degeneracy of $u$ at the boundary $y=1/2$ require some modifications.

If $j=0$, then the integral equals $\int_{-1/2}^{1/2}4 (v')^2$ in terms of  $v=\sqrt u$, for which the optimal function $v$ is linear and the corresponding value of the integral equals $4z$. This proves the identity~\eqref{eq:lemma:var-probl} for the case $j=0$. 

If $j<0$, then one can estimate
\begin{align*}
\limsup_{a\uparrow 1/2}\int_{-1/2}^a \frac{(j + u')^2}{u} \dd y &= 
\limsup_{a\uparrow 1/2}\int_{-1/2}^a\frac{\,j^2  + {u'}^2 + 2ju'}{u} \dd y \\
&\geq
\limsup_{a\uparrow 1/2} \,2j \bigl(\log u(a)- \log z\bigr) = +\infty,
\end{align*}
which establishes the identity~\eqref{eq:lemma:var-probl} for the case $j<0$. In the case $j>0$ but $z=0$, a similar calculation at $y=-1/2$ yields the same conclusion. 

The final case to consider is $j,z>0$. Following the argument of~\cite{LieroMielkePeletierRenger17}, we set $f(u,u') = (j+u')^2/4u$, such that the Euler-Lagrange equation is $-(\partial_{u'} f(u,u'))' + \partial_u f(u,u') = 0$. It follows by differentiating (or applying Noether's theorem) that the Hamiltonian $u'\partial_{u'}f(u,u') - f(u,u') = ({u'}^2 - j^2)/4u$ is constant on $[-1/2,1/2]$, say equal to $\gamma/4$. By differentiating the resulting equation ${u'}^2 = j^2 + \gamma u$ we find that all solutions are second-order polynomials, and by applying the boundary conditions on $u$ we obtain~\eqref{eq:lemma-var-probl:optimal-u} and $\gamma = 4(z-j)$. The identity~\eqref{eq:lemma:var-probl} then follows from a direct calculation. 	
\end{proof}

\subsection{Proof of Lemma~\ref{l:reg-rho-j}}
\label{app:pf:l:reg-rho}
Results of this type are fairly standard; similar arguments can be found in~\cite[Th.~2.3]{DuongLamaczPeletierSharma17} or~\cite[Lemmas~8.4 and~8.5]{GavishNyquistPeletier2019}. Since we could not find a complete result, we provide a proof here. For the length of this proof we use subscripts $t$ to indicate time slices. 

\textit{Step 1: Alternative duality estimate.}
We have
\begin{align}
\notag
\cI_\e(\rho,j)&\geq 
\sup_{f\in C^{1,2}_b(\tmsp)} \int_\R \bigl(\rho_T f_T - \rho_0 f_0\bigr) \\
&\qquad {}- \int_\tmsp \rho_t \Bigl( \partial_t f_t 
  + \e\tau_\e \Bigl(  \partial_{xx}f_t - \frac1\e \partial_x f_t V' + \frac12 (\partial_x f_t)^2\Bigr)\Bigr).
 \label{ineq:duality-f}
\end{align}
One first obtains this inequality for $f\in C^{1,2}_c(\tmsp)$ by substituting $b=\partial_x f$ in~\eqref{defeq:I_e2} and using the narrow continuity of $\rho$ with a standard truncation and regularization in time. For general  $f\in C^{1,2}_b(\tmsp)$ the inequality follows from a regularized truncation in space, using the finiteness of the measure $\rho$ on $\tmsp$.

\textit{Step 2: Dual equation.}
Fix $\varphi\in C_c^{0,1}(\tmsp)$ and $\psi\in C_c^\infty(\R)$. Define $g\in C_b^{1,2}(\tmsp)$ as the solution of the backward parabolic equation
\begin{align*}
&\partial_tg_t + \e\tau_\e \Bigl(\partial_{xx}g_t -\frac1\e \partial_x g_t V' \Bigr) 
  = \frac{\e\tau_\e}2 g_t \Bigl( \frac12 \varphi_t^2 - \partial_x \varphi_t + \frac1\e \varphi_t V'\Bigr)\\
&g_T = e^{\psi/2} \\\
&g \text{ is bounded.}
\end{align*}
Such a solution exists by e.g.~\cite[Th.~1.12]{Friedman64}. By calculating the derivative explicitly, we find that 
\[
\frac \dd{\dd t} \int_\R g_t^2 e^{V/\e}
= 2\e\tau_\e \int_\R e^{-V/\e} \Bigl(\partial_x g_t + \frac12 g \varphi\Bigr)^2 \geq 0,
\]
implying that
\begin{equation}
\label{ineq:int-g}
\int_\R e^\psi e^{-V/\e}  = \int_\R  g_T^2 e^{-V/\e} \geq \int_\R g_0^2  e^{-V/\e} .
\end{equation}
The function $f := 2\log g$ is an element of $C^{1,2}_b(\tmsp)$ and satisfies the equation
\[
\partial_t f_t + \e\tau_e\Bigl( \partial_{xx} f_t + \frac12 (\partial_x f_t)^2 - \frac1\e \partial_x f_t V'\Bigr)
= \e\tau_e \Bigl(\frac12 \varphi_t^2 - \partial_x \varphi_t + \frac1\e \varphi_t V'\Bigr)
\]
with final datum $f_T = \psi$. Substituting in~\eqref{ineq:duality-f} yields
\begin{align*}
\cI_\e(\rho,j) 
&\geq \int_\R \bigl[\rho_T\psi - \rho_0 f_0\bigr]
- \int_\tmsp \rho\e\tau_\e \Bigl(\frac12 \varphi_t^2 - \partial_x \varphi_t + \frac1\e \varphi_t V'\Bigr).
\end{align*}
By reorganizing this inequality and applying the Donsker-Varadhan dual characterization of the relative entropy we find
\begin{align*}
\int_\R \rho_T \psi &- \log \int_\R \gamma_\e e^\psi
  + \e\tau_\e\int_\tmsp \rho\Bigl(\frac12 \varphi_t^2 - \partial_x \varphi_t + \frac1\e \varphi_t V'\Bigr)\\
&\leq \cI_\e(\rho,j) + \int_\R \rho_0 f_0 - \log \int_\R \gamma_\e e^\psi\\
&\leq \cI_\e(\rho,j) + \RelEnt(\rho_0|\gamma_\e) 
  - \log \frac{ \int_\R \gamma_\e e^\psi}{\int_\R \gamma_\e e^{f_0}}
  \;\stackrel{\eqref{ineq:int-g}} \leq \;
  \cI_\e(\rho,j) + \RelEnt(\rho_0|\gamma_\e).
\end{align*}
Taking the supremum over $\psi\in C_c^\infty(\R)$ and $\varphi\in C^{1,2}_c(\tmsp)$, and applying also the dual formulation of the Fisher Information~\cite[Lemma~D.44]{FengKurtz06}, we find
\[
\RelEnt(\rho_T|\gamma_\e) + \e\tau_\e \int_0^T \cR(\rho_t|\gamma_\e,\R)\, \dd t
\leq \cI_\e(\rho,j) + \RelEnt(\rho_0|\gamma_\e).
\]
Summarizing,  $\cI_\e(\rho,j)<\infty$ implies that $\rho$ is absolutely continuous on $\tmsp$ with respect to $\gamma_\e(\dd x)\dd t$, or equivalently to with respect to Lebesgue measure on $\tmsp$; the density $u := \dd\rho/\dd\gamma_\e$ satisfies $\partial_x u\in L^1_{\mathrm{loc}}(\tmsp)$. This proves parts~\ref{l:reg-rho-j:rho} and~\ref{l:reg-rho-j:rho_x} of Lemma~\ref{l:reg-rho-j}.

\smallskip
\textit{Step 3: Regularity of $j$.}
To show that $j\ll \rho$, we use the regularity of $u$ to rewrite
\[
\cI_\e(\rho,j) = \sup_{b\in C_c^1(\tmsp)} \int_\tmsp \Bigl[ b\Bigl(j + \e\tau_\e \rho \frac{\partial_x u}u\Bigr) - \frac{\e\tau_\e}2 \rho b^2\Bigr)\Bigr].
\]
By the dual characterization of $L^2(\tmsp;\rho)$, finiteness of $\cI_\e(\rho,j)$ implies that there exists $v\in L^2(\tmsp;\rho)$ with $j = v\rho$, and we have the estimate
\begin{equation}
\label{ineq:v-square-Ie}
\frac1{2\e\tau_\e} \int_\tmsp \rho \Bigl( v + \e\tau_\e \frac{\partial_x u}u\Bigr)^2 
\leq
\cI_\e(\rho,j).
\end{equation}

\smallskip
\textit{Step 4: Rewriting $\cI_\e$.}
Finally, to show the identity~\eqref{eq:Ie=E+R}, we note that $v\in L^2(\tmsp;\rho)$ implies that the curve $t\mapsto \rho_t$ is absolutely continuous in the Wasserstein sense~\cite[Th.~8.3.1]{AmbrosioGigliSavare2008}. By \cite[Th.~10.4.9]{AmbrosioGigliSavare2008}, the bound on $\partial_x u_t/u_t$ implies that the global Wasserstein slope $|\partial E_\e|(\rho_t)$ is bounded and $\partial_x u_t/u_t$ is an element of the Fr\'echet subdifferential. Finally, by the chain rule~\cite[Sec.~10.1.2]{AmbrosioGigliSavare2008}, we have 
\[
\int_\tmsp v\frac{\partial_x u} u \, \rho = E_\e(\rho(T))-E_\e(\rho(0)).
\]
Expanding the square in~\eqref{ineq:v-square-Ie} establishes~\eqref{eq:Ie=E+R} with an inequality. 

\smallskip
\textit{Step 5. Inverting the argument.}
The argument up to now can be summarized as ``if $\cI_\e(\rho,j)$ is finite, then $\rho$ and $j$ are regular and identity~\eqref{eq:Ie=E+R} holds as an inequality''. Vice versa, if the regularity conditions on $\rho$ and $j$ are satisfied and the right-hand side in~\eqref{eq:Ie=E+R} is finite, then the calculations can be reversed, and we find that $\cI_\e$ is finite and that the inequality is an identity. This concludes the proof.

\bibliography{mikola}
\bibliographystyle{alphainitials}
\end{document}